\documentclass[a4paper]{article}

\usepackage{graphicx,array,delarray}
\usepackage{latexsym}
\usepackage{amscd}
\usepackage{amsfonts, amsmath, amssymb}
\usepackage{color}
\usepackage{hyperref}

\usepackage{amsthm}

\usepackage{multicol}

\input xy
\xyoption{all}

\newtheorem{theorem}{\textbf{Theorem}}

\newtheorem{example}{\textbf{Example}}

\newtheorem{lemma}[theorem]{\textbf{Lemma}}

\newtheorem{proposition}[theorem]{\textbf{Proposition}}

\newtheorem{remark}[theorem]{Remark}

\theoremstyle{definition}
\newtheorem{definition}{Definition}

\title{Associative triple trisystems and standard embeddings}
\author{ {\bf{Ra\'ul Felipe}}\\ \\
        CIMAT \\
        Callej\'on Jalisco s/n Mineral de Valenciana \\
        Guanajuato, Gto, M\'exico.\\
        \texttt{raulf@cimat.mx}
        \\
        \url{https://orcid.org/0000-0002-2809-3232}
        \and
        \textbf{Guillermo Vera de Salas\footnote{Corresponding author.}} \\ \\
        Universidad Rey Juan Carlos, Madrid (Spain) \\
        \texttt{guillermo.vera@urjc.es}
        \\
        \url{https://orcid.org/0000-0002-3517-7223}
        }

\begin{document}

\maketitle

\begin{abstract}

Building on the established theories of Jordan triple disystems and Leibniz triple systems, we introduce and develop the theory of associative triple trisystems, filling a significant gap in the existing framework. We establish the classical relationships between associative, Jordan, and Lie triple systems within the context of trisystems. We present a significant example by equipping the space of matrices with a non-trivial associative dialgebra structure. We conclude defining the concept of di-endomorphisms of any module, which enables the construction of the standard embedding for any associative triple trisystem.

\bigskip

\noindent{\it{2020 Mathematics Subject Classification (MSC2020):}}
$17A30$, $17A60$, $17C99$.


\end{abstract}


\section{Introduction}

In the framework of triple systems and algebras there exist many well-known connections between their associative, Lie and Jordan structures. A classical summary of these relationships is as follows: given an associative algebra $A$ with multiplication denoted by juxtaposition $ab$, one can obtain a Lie algebra $A^{(-)}$ and a Jordan algebra $A^{(+)}$ by defining the new products $[a,b] := ab - ba$ and $a \bullet b := ab + ba$, respectively. Furthermore, $A$ can be naturally regarded as an associative triple system of the first kind (ATS1) by defining the triple product $\{a,b,c\}:=abc$. If $A$ is equipped with an involution $*$, then the modified product $\{a,b,c\} := ab^*c$ gives rise to an associative triple system of the second kind (ATS2) (see \cite{loos}).

In addition to these associative structures, additional triple systems arise. Specifically, given an ATS1 or ATS2 with triple product $\{a,b,c\}$, the product 
$$
    \langle a,b,c\rangle := \{a,b,c\} + \{c,b,a\}
$$
yields a Jordan triple system (JTS), while the triple product
$$
    [a,b,c] := \{a,b,c\} - \{b,a,c\} - \{c,a,b\} + \{c,b,a\}
$$
defines a Lie triple system (LTS). Conversely, any ATS1 and ATS2 can be embedded, in different ways for each case, in an associative algebra via its so-called standard embedding (see \cite{meyberg}).

Beyond these classical constructions, the study of non-anticommutative generalizations of Lie and Jordan algebras has led to the development of Leibniz algebras (see \cite{lod1}) and Jordan dialgebras (see \cite{brem3}, \cite{Kol} and \cite{fel}), both of which naturally arise from the structure of an associative dialgebra.

A fundamental result of Loday in \cite{lod2} states that every Leibniz algebra can be embedded in an associative dialgebra, which is a module equipped with two associative products $\dashv$ and $\vdash$ satisfying specific compatibility conditions. The induced bracket $[a,b]:= a \dashv b - b\vdash a$ endows $D$ with the structure of a (right) Leibniz algebra denoted by $D^{(-)}$. 

Motivated by these ideas, Velásquez and Felipe (2008) introduced the notion of quasi-Jordan algebras, and the notion (right) Jordan dialgebras as later Bremner in \cite{brem3} and Kolesnikov \cite{Kol} described. For an associative dialgebra $D$, the product $a \bullet b := a \dashv b + b\vdash a$ gives rise to a (right) Jordan dialgebra, denoted by $D^{(+)}$.

Building on this framework, in \cite{brem1}, Bremner, Felipe and Sánchez-Ortega introduced the notion of Jordan triple disystems (JTD), which arise as the structure obtained by applying the Kolesnikov-Pozhidaev (KP) algorithm to the identities of a JTS. They also demonstrated that a JTD can be constructed from an associative dialgebra $D$ with products $a \dashv b$ and $a \vdash b$ using the two triple products:
$$
    \{a,b,c\}_1 :=  a\dashv b \dashv c + c \vdash b \vdash a \quad \text{ and } \quad  \{a,b,c\}_2 = a \vdash b \dashv c + c \vdash b \dashv a.
$$
Subsequently, in \cite{brem2}, Bremner and Sánchez-Ortega introduced the concept of a Leibniz triple system (LeibTS), again after applying the KP algorithm to the identities of a LTS. They showed that an associative dialgebra $D$ also gives rise to the triple product
$$
    [a,b,c] := a \dashv b \dashv c - b \vdash a \dashv c-  c \vdash a \dashv b + c \vdash b \vdash a
$$
endowing it with the structure of an LeibTS.

Surprisingly, while the classical theory provides a unified perspective on associative, Jordan and Lie triple systems, a corresponding associative triple trisystem (ATT) theory was missing. Such a theory would allow for a comprehensive treatment of triple trisystems, mirroring the classical framework and establishing a natural connection between JTD's, LeibTS's, and associative dialgebras. Therefore, in this work, we aim to fill this gap by developing the theory of ATT's and studying their standar embedding in associative dialgebras.

Thus, we structure this work around the following sections:

\smallskip
$1$) Describe the KP algorithm, which serves as the fundamental tool for extending the notion of an associative triple system  to the dialgebra framework.

\smallskip
$2$) Review the theories of Jordan triple disystems and Leibniz triple systems and associative triple systems, remembering their definitions.

\smallskip
$3$) Introduce the concept of an associative triple trisystem of the first kind (ATT1) as the result of applying the KP algorithm to the identities of ATS1. We will show how an ATT1 arises from an associative dialgebra and how it naturally leads to the structures of JTD's and LeibTS's, mirroring the classical case.

\smallskip
$4$) Introduce the concept of an associative triple trisystem of the second kind (ATT2) and demonstrate how an ATT2 can be constructed from an associative dialgebra with involution. Moreover, we also prove that any ATT2 gives rise a JTD and a LeibTS. 

In this section, we also present a significant example involving matrices endowed with a non-trivial associative dialgebra structure. This example demonstrates that the newly introduced definitions are indeed non-trivial and closely connected to familiar mathematical objects, such as matrix spaces, although equipped with a structure not previously considered. Moreover, it serves as a rich source for constructing examples of finite-dimensional Leibniz algebras.

\smallskip
$5$) Although the definition of a di-endomorphism of a Jordan dialgebra was introduced in \cite{GK}, we present a more elementary description of a di-endomorphism of a module inspired on the identities obtained after applying the KP algorithm to a particular algebra. This approach enables us to establish a standard embedding for ATT1's and ATT2's. We then prove that every ATT1 and ATT2 can be canonically embedded into an associative dialgebra, giving the analogous standard embeddings presented in \cite{meyberg}.

\smallskip

By developing this framework, we establish a natural extension of classical triple system theory to the setting of dialgebras, providing a unifiying perspective that connects associative, Jordan triple disystems and Leibniz triple systems in a systematic way.

\section{Preliminaries}

In this work, $\phi$ is an scalar ring with $\frac{1}{2} \in \phi$.

We start by recalling the definition of  associative triple system.
Let $A$ be a $\phi$-module. A trilinear map $\{ \, , \,,\,\}:\,\,A\times A\times A\longrightarrow A$ endows
$A$ with the structure of an  \textbf{associative triple system of the first kind} if it satisfies the relations
\begin{equation}\label{tpada1}
\{\{a,b,c\},d,e\}=\{a,\{b,c,d\},e\}=\{a,b,\{c,d,e\}\},
\end{equation}
for all $a,b,c,d,e\in A$,.

On the other hand, if the trilinear map satisfies the relations,
\begin{equation}\label{ats2}
\{\{a,b,c\},d,e\} = \{a, \{d,c,b\},e\} = \{a,b,\{c,d,e\}\},
\end{equation}
for all $a,b,c,d,e \in A$, then $A$ it is an \textbf{associative triple system of the second kind}. For more detail of these structures see  \cite{loos} and \cite{meyberg}.

Next, we expose the KP algorithm. Initially, it was introduced by Kolesnikov \cite{Kol} and Pozhidaev \cite{Pozhid2}.
We recall here the most general version of this algorithm  appeared in \cite{brem1}. Summarizing, it converts a multilinear
polynomial identity of degree $d$ for an $n$-ary operation into $d$ multilinear
identities of degree $d$ for $n$ new $n$-ary operations.

\begin{definition}

\textbf{KP Algorithm}

Part 1:
We consider a multilinear $n$-ary operation, denoted by the symbol
  \begin{equation} \label{operation}
  \{-,-,\dots,-\}
  \qquad
  \text{($n$ arguments)}.
  \end{equation}
Given a multilinear polynomial identity of degree $d$ in this operation,
we describe the application of the algorithm to one monomial in the identity,
and from this the application to the complete identity follows by linearity.
Let $\overline{a_1 a_2 \dots a_d}$
be a multilinear monomial of degree $d$, where the bar denotes some placement of $n$-ary operation symbols.
We introduce $n$ new $n$-ary operations, denoted by the same symbol but distinguished by subscripts:
  \begin{equation} \label{noperations}
  \{-,-,\dots,-\}_1,
  \quad
  \{-,-,\dots,-\}_2,
  \quad
  \dots,
  \quad
  \{-,-,\dots,-\}_n.
  \end{equation}
For each $i \in \{1, 2, \dots, d\}$ we convert the monomial $\overline{a_1 a_2 \dots a_d}$
in the original $n$-ary operation \eqref{operation}
into a new monomial of the same degree $d$ in the $n$ new $n$-ary operations \eqref{noperations},
according to the following rule which is based on the position of $a_i$.
For each occurrence of the original operation symbol in the monomial,
either $a_i$
occurs within one of the $n$ arguments or not, and we have the following
cases:
  \begin{itemize}
  \item
  If $a_i$ occurs within the $j$-th argument then we convert the original
  operation symbol
  $\{\dots\}$ to the $j$-th new operation symbol $\{\dots\}_j$.
  \item
  If $a_i$ does not occur within any of the $n$ arguments, then either
    \begin{itemize}
    \item
    $a_i$ occurs to the left of the original operation symbol, in which case
    we convert $\{\dots\}$ to the first new operation symbol $\{\dots\}_1$, or
    \item
    $a_i$ occurs to the right of the original operation symbol,
    in which case we convert $\{\dots\}$ to the last new operation
    symbol $\{\dots\}_n$.
    \end{itemize}
  \end{itemize}
In this process, we call $a_i$ the central argument of the monomial.

Part 2:
In addition to the identities constructed in Part 1, we also include the following identities
for all $i, j \in \{ 1, 2, \dots, n\}$ with $i \ne j$
and all $k, \ell \in \{ 1, 2, \dots, n\}$:
  \begin{align*}
  &
  \{ a_1, \dots, a_{i-1}, \{ b_1, \cdots, b_n \}_k, a_{i+1}, \dots, a_n \}_j
  \equiv
  \\
  &
  \{ a_1, \dots, a_{i-1}, \{ b_1, \cdots, b_n \}_\ell, a_{i+1}, \dots, a_n \}_j.
  \end{align*}
This identity says that the $n$ new operations are interchangeable in
the $i$-th argument of the $j$-th new operation when $i \ne j$.
\end{definition}

\begin{example}
    The five axioms of an associative dialgebra can be obtained applying the KP algorithm to the associative identity $(ab)c \equiv a(bc)$. The juxtaposition operation produces two new operations $\{ \, , \,\}_1 = \, \dashv$ and $\{ \, , \,\}_2 =\,\vdash$. Since the associative identity has degree $3$, Part $1$ yields three new identities of degree $3$ by making $a,b,c$ in turn be the central argument:
    \begin{align*}
        (a \dashv b) \dashv c &\equiv a \dashv ( b \dashv c)
        \\
        (a \vdash b) \dashv c& \equiv a \vdash (b \dashv c)
        \\
        (a \vdash b) \vdash c&\equiv  a \vdash (b \vdash c)
    \end{align*}
    and Part $2$ produces the following two identities:
    $$
        a \dashv (b \dashv c) \equiv  a \dashv (b \vdash c), \quad \text{ and } \quad (a \dashv b) \vdash c \equiv (a \vdash b) \vdash c.
    $$
\end{example}

\begin{example}
    The defining identities for left-symmetric dialgebras (see \cite{fel1} for more details) can be obtained by applying the KP algorithm to the left-symmetric identity, the which we write in the form $\{a,\{b,c\}\}-\{\{a,b\},c\}\equiv\{b,\{a,c\}\}-\{\{b,a\},c\}$. The original operation produces two new operations $\{ \, , \,\}_{1}$ and $\{ \, , \,\}_{2}$. Since left-symmetric identy has degree $3$, Part $1$ produces three new identities of degree $3$ by making $a,b,c$ in turn the central argument:
    
    \begin{equation*}
    \{a,\{b,c\}_{1}\}_{1}-\{\{a,b\}_{1},c\}_{1}\equiv\{b,\{a,c\}_{1}\}_{2}-\{\{b,a\}_{2},c\}_{1},
    \end{equation*}
    \begin{equation*}
    \{a,\{b,c\}_{1}\}_{2}-\{\{a,b\}_{2},c\}_{1}\equiv\{b,\{a,c\}_{1}\}_{1}-\{\{b,a\}_{1},c\}_{1},
    \end{equation*}
    \begin{equation*}
    \{a,\{b,c\}_{2}\}_{2}-\{\{a,b\}_{2},c\}_{2}\equiv\{b,\{a,c\}_{2}\}_{2}-\{\{b,a\}_{2},c\}_{2},
    \end{equation*}
    and Part $2$ produces the following two identities:
    \begin{equation*}
    \{a,\{b,c\}_{1}\}_{1}\equiv\{a,\{b,c\}_{2}\}_{1},\,\,\{\{a,b\}_{1},c\}_{2}\equiv\{\{a,b\}_{2},c\}_{2}.
    \end{equation*}
    
    Note that the last two identities are the same. Moreover,
    the first three identities constitute the set of axioms
    for a left-symmetric dialgebra proposed for the first author in \cite{fel1} some time ago.
\end{example}

In \cite{brem1} the authors applied the KP algorithm to the identities of a Jordan triple algebras to obtain a new variety of nonassociative triple systems, called Jordan triple disystems.  We now recall this concept here

\begin{definition}
A \textbf{Jordan triple disystem} is a $\phi$-module  with two trilinear operations $\langle \, , \, , \,\rangle_i \colon D \times D \times D \rightarrow D$, with $i \in \{1,2\}$, satisfying the following equalities:
\begin{itemize}
\item[(JTD1)]  $\langle a,b,c\rangle_2 = \langle c,b,a\rangle_2$,
\item[(JTD2)]  $\langle a, \langle b, c, d \rangle_1, e \rangle_1 =
\langle a, \langle b, c, d \rangle_2, e \rangle_1$,
\item[(JTD3)] $\langle a, b, \langle c, d, e \rangle_1 \rangle_1 =
\langle a, b, \langle c, d, e \rangle_2 \rangle_1$,
\item[(JTD4)] $\langle \langle a, b, c \rangle_1, d, e \rangle_2 =
\langle \langle a, b, c \rangle_2, d, e \rangle_2$,
\item[(JTD5)] $\langle \langle e, d, c \rangle_1, b, a \rangle_1=\langle \langle e, b, a \rangle_1, d, c \rangle_1-\langle e, \langle d, a, b \rangle_1, c \rangle_1+\langle e, d, \langle c, b, a \rangle_1 \rangle_1$,
\item[(JTD6)] $\langle \langle e, d, c \rangle_2, b, a \rangle_1=\langle \langle e, b, a \rangle_1, d, c \rangle_2-\langle e, \langle d, a, b \rangle_1, c \rangle_2+\langle e, d, \langle c, b, a \rangle_1 \rangle_2$,
\item[(JTD7)] $\langle a, b, \langle c, d, e \rangle_1 \rangle_1=\langle \langle a, b, c \rangle_1, d, e \rangle_1-\langle c, \langle b, a, d \rangle_2, e \rangle_2+\langle \langle a, b, e \rangle_1, d, c \rangle_1$,
\item[(JTD8)] $\langle a, b, \langle c, d, e \rangle_1 \rangle_2=\langle \langle a, b, c \rangle_2, d, e \rangle_1 -\langle c, \langle b, a, d \rangle_1, e \rangle_2 +\langle \langle a, b, e \rangle_2, d, c \rangle_1$,
\end{itemize}
for all $a,b,c,d,e\in D$.
\end{definition}

\begin{example}
Let $A$ be a differential associative algebra in the sense of Loday \cite{lod2}:
that is, $A$ is an associative algebra with product $a \cdot b$
together with a linear map $d\colon A \to A$ such that $d^2 = 0$ and
$d(a \cdot b) = d(a) \cdot b + a \cdot d(b)$ for all $a, b \in A$.
One endows $A$ with a dialgebra structure by defining
$a \dashv b = a \cdot d(b)$ and $a \vdash b = d(a) \cdot b$.
Then $A$ becomes a Jordan triple disystem by defining
  \[
  \langle a,b,c \rangle _1 = a \cdot d(b) \cdot d(c) + d(c) \cdot d(b) \cdot a,
  \quad
  \langle a,b,c \rangle _2 = d(a) \cdot b \cdot d(c) + d(c) \cdot b \cdot d(a).
  \]
\end{example}

\smallskip
In \cite{brem2} was proposed a ternary version of a Leibniz algebra, called the Leibniz triple system which we remember now

\begin{definition}\cite{brem2}A \textbf{Leibniz triple system} is a $\phi$-module $T$ with a trilinear operation $T \times T \times T \rightarrow T$
denoted by $[ \, , \, ,\, ]$ satisfying the following two axioms:
\begin{itemize}
    \item[(LTSA)] $[a,b, [c,d,e]]- [[a,b,c],d,e] + [[a,b,d],c,e]- [[a,b,e],d,c] + [[a,b,e],c,d]=0$,
    \item[(LTSB)] $[a,[b,c,d],e]-[[a,b,c],d,e] + [[a,c,b],d,e] + [[a,d,b],c, e]
    -[[a,d,c],b,e]=0$,
\end{itemize}
for all $a,b,c,d,e\in T$.
\end{definition}

In the same article, it was proved that if $(D,\vdash,\dashv)$ is an associative dialgebra then the trilinear map
\begin{equation*}
[a,b,c]=a\dashv b \dashv c -b\vdash a\dashv c- c\vdash a \dashv b + c\vdash b\vdash a,
\end{equation*}
equips $D$ of a structure of Leibniz triple system.

\section{Associative triple trisystems of the first kind}

In this section, we apply the KP algorithm to the identities of an associative triple system, allowing us to derive a variety of triple trisystems. These new identities establish the necessary framework for proving classical results in di-structures. Specifically, we show that an associative triple trisystem of the first kind arises from an associative dialgebra (Proposition \ref{th:asstoass}), while a Jordan triple disystem and a Leibniz triple system emerge from an associative triple trisystem of the first kind (Theorem \ref{th:asstojordan} and Theorem \ref{th:asstoleibniz}).

\smallskip
We apply Part $1$ of the KP algorithm to (\ref{tpada1}) and obtain five relations:
\begin{align*}
& \{\{a,b,c\}_{1},d,e\}_{1}\equiv \{a,\{b,c,d\}_{1},e\}_{1}\equiv \{a,b,\{c,d,e\}_{1}\}_{1}, \\
& \{\{a,b,c\}_{2},d,e\}_{1}\equiv \{a,\{b,c,d\}_{1},e\}_{2}\equiv \{a,b,\{c,d,e\}_{1}\}_{2}, \\
& \{\{a,b,c\}_{3},d,e\}_{1}\equiv \{a,\{b,c,d\}_{2},e\}_{2}\equiv \{a,b,\{c,d,e\}_{1}\}_{3}, \\
& \{\{a,b,c\}_{3},d,e\}_{2}\equiv \{a,\{b,c,d\}_{3},e\}_{2}\equiv \{a,b,\{c,d,e\}_{2}\}_{3}, \\
& \{\{a,b,c\}_{3},d,e\}_{3}\equiv \{a,\{b,c,d\}_{3},e\}_{3}\equiv \{a,b,\{c,d,e\}_{3}\}_{3}.
\end{align*}

\smallskip
Part $2$ of the algorithm produces six relations
\begin{align*}
& \{\{a,b,c\}_{1},d,e\}_{2}\equiv \{\{a,b,c\}_{2},d,e\}_{2}\equiv \{\{a,b,c\}_{3},d,e\}_{2}, \\
& \{\{a,b,c\}_{1},d,e\}_{3}\equiv \{\{a,b,c\}_{2},d,e\}_{3}\equiv \{\{a,b,c\}_{3},d,e\}_{3}, \\
& \{a,\{b,c,d\}_{1},e\}_{1}\equiv \{a,\{b,c,d\}_{2},e\}_{1}\equiv \{a,\{b,c,d\}_{3},e\}_{1}, \\
& \{a,\{b,c,d\}_{1},e\}_{3}\equiv \{a,\{b,c,d\}_{2},e\}_{3}\equiv \{a,\{b,c,d\}_{3},e\}_{3}, \\
& \{a,b,\{c,d,e\}_{1}\}_{1}\equiv \{a,b,\{c,d,e\}_{2}\}_{1}\equiv \{a,b,\{c,d,e\}_{3}\}_{1}, \\
& \{a,b,\{c,d,e\}_{1}\}_{2}\equiv \{a,b,\{c,d,e\}_{2}\}_{2}\equiv \{a,b,\{c,d,e\}_{3}\}_{2}.
\end{align*}

Thus, we introduce our first definition.

\begin{definition}\label{cor1}
Let $A$ be a $\phi$-module. Three trilinear maps $\{ \, , \, ,\,\}_{i}:\,\,A\times A\times A\longrightarrow A$, with $i \in \{1,2,3\}$, endow $A$ with a structure of \textbf{associative triple trisystem of the first kind} if these satisfy the following relations:
\begin{align}
&\label{cor2} \{\{a,b,c\}_{1},d,e\}_{1}=\{a,\{b,c,d\}_{1},e\}_{1}= \{a,b,\{c,d,e\}_{1}\}_{1}, \\
&\label{cor3} \{\{a,b,c\}_{2},d,e\}_{1}= \{a,\{b,c,d\}_{1},e\}_{2}= \{a,b,\{c,d,e\}_{1}\}_{2}, \\
&\label{cor4} \{\{a,b,c\}_{3},d,e\}_{1}= \{a,\{b,c,d\}_{2},e\}_{2}= \{a,b,\{c,d,e\}_{1}\}_{3}, \\
&\label{cor5} \{\{a,b,c\}_{3},d,e\}_{2}= \{a,\{b,c,d\}_{3},e\}_{2}= \{a,b,\{c,d,e\}_{2}\}_{3}, \\
&\label{cor6} \{\{a,b,c\}_{3},d,e\}_{3}= \{a,\{b,c,d\}_{3},e\}_{3}= \{a,b,\{c,d,e\}_{3}\}_{3}, \\
&\label{cor7} \{\{a,b,c\}_{1},d,e\}_{2}= \{\{a,b,c\}_{2},d,e\}_{2}= \{\{a,b,c\}_{3},d,e\}_{2}, \\
&\label{cor8} \{\{a,b,c\}_{1},d,e\}_{3}= \{\{a,b,c\}_{2},d,e\}_{3}= \{\{a,b,c\}_{3},d,e\}_{3}, \\
&\label{cor9} \{a,\{b,c,d\}_{1},e\}_{1}= \{a,\{b,c,d\}_{2},e\}_{1}= \{a,\{b,c,d\}_{3},e\}_{1}, \\
&\label{cor10} \{a,\{b,c,d\}_{1},e\}_{3}= \{a,\{b,c,d\}_{2},e\}_{3}= \{a,\{b,c,d\}_{3},e\}_{3}, \\
&\label{cor11} \{a,b,\{c,d,e\}_{1}\}_{1}= \{a,b,\{c,d,e\}_{2}\}_{1}= \{a,b,\{c,d,e\}_{3}\}_{1}, \\
&\label{cor12} \{a,b,\{c,d,e\}_{1}\}_{2}= \{a,b,\{c,d,e\}_{2}\}_{2}= \{a,b,\{c,d,e\}_{3}\}_{2}.
\end{align}
for all $a,b,c,d,e\in A$.
\end{definition}

Let $A$ be an  associative triple trisystem of the first kind. Define $A^{ann}$ as the subspace spanned by the elements of the form $\{a,b,c\}_{1}-\{a,b,c\}_{2}$ and $\{a,b,c\}_{1}-\{a,b,c\}_{3}$ for all $a,b,c \in A$. It is clear that $A$ is an associative triple system of the first kind if and only if $A^{ann}=\{0\}$. Let us suppose that $A=A^{ann}\oplus \mathcal{A}$ as subspaces where $\mathcal{A}$ is closed with respect to the trilinear maps $\{ \,, \,,\,\}_{i}$ for $i\in \{1,2,3\}$, then $\mathcal{A}$ is an  associative triple system of the first kind.

In the following result we prove that from any associative dialgebra $D$, we can define three triple products that endow $D$ of structure of associative triple trisystem of the first kind.

\begin{proposition}\label{th:asstoass}Let $(D,\vdash,\dashv)$ be an associative dialgebra. Define
\begin{align}
&\label{cor13}\{a,b,c\}_{1}=(a\dashv b) \dashv c = a \dashv (b \dashv c), \\
&\label{cor14}\{a,b,c\}_{2}=(a\vdash b)\dashv c = a\vdash (b\dashv c),  \\
&\label{cor15}\{a,b,c\}_{3}=(a\vdash b)\vdash c = a\vdash (b\vdash c).
\end{align}
for all $a,b,c\in D$. Then, these products convert to $D$ in a  associative triple trisystem of the first kind.
\end{proposition}
\begin{proof}Taking into account that both $\vdash$ and $\dashv$ are associative products, we have
(\ref{cor2}) and (\ref{cor6}). Now,
\begin{align}
&\{\{a,b,c\}_{2},d,e\}_{1}=((a\vdash (b \dashv c))\dashv d) \dashv e =(a\vdash ((b\dashv c) \dashv d))\dashv e \nonumber \\
&=\{a,\{b,c,d\}_{1},e\}_{2}=a\vdash (((b\dashv c) \dashv d)\dashv e)=a\vdash ((b\dashv (c \dashv d))\dashv e) \nonumber \\
&=a\vdash (b\dashv ((c \dashv d)\dashv e))=\{a,b,\{c,d,e\}_{1}\}_{2}, \nonumber
\end{align}
for all $a,b,c,d,e\in D$. It shows (\ref{cor3}). Now, we pass to prove (\ref{cor4}):
\begin{align}
&\{\{a,b,c\}_{3},d,e\}_{1}=((a\vdash (b\vdash c))\dashv d)\dashv e=a\vdash (((b\vdash c) \dashv d)\dashv e) \nonumber \\
&=\{a,\{b,c,d\}_{2},e\}_{2}=a\vdash (((b\vdash c)\dashv d)\dashv e)=a\vdash ((b\vdash c)\dashv (d\dashv e)) \nonumber \\
&=a\vdash (b\vdash (c\dashv d\dashv e))=a\vdash (b \vdash (c\dashv (d\dashv e)))=\{a,b,\{c,d,e\}_{1}\}_{3}, \nonumber
\end{align}
for all $a,b,c,d,e\in D$. On other hand
\begin{align}
&\{\{a,b,c\}_{3},d,e\}_{2}=((a\vdash (b\vdash c))\vdash d)\dashv e=(a\vdash ((b\vdash c)\vdash d))\dashv e \nonumber \\
&=\{a,\{b,c,d\}_{3},e\}_{2}=a\vdash (b\vdash ((c\vdash d)\dashv e))=a\vdash (b\vdash \{c,d,e\}_2) \nonumber \\
&=\{a,b,\{c,d,e\}_{2}\}_{3}, \nonumber
\end{align}
for all $a,b,c,d,e\in D$. Hence, we have established (\ref{cor5}).

It remains to prove the axioms derived by the Part 2 of the KP algorithm. Note that
\begin{align}
&\{\{a,b,c\}_{1},d,e\}_{2}=(a\dashv (b\dashv c))\vdash (d \dashv e)=((a\dashv (b\dashv c))\vdash d) \dashv e \nonumber \\
&=((a\vdash (b\dashv c))\vdash d) \dashv e=\{\{a,b,c\}_{2},d,e\}_{2}=(a\vdash ((b\dashv c)\vdash d))\dashv e \nonumber \\
&=(a\vdash ((b\vdash c)\vdash d))\dashv e=(\{a,b,c\}_3\vdash d)\dashv e=\{\{a,b,c\}_{3},d,e\}_{2}, \nonumber
\end{align}
for all $a,b,c,d,e\in D$. These calculations imply that (\ref{cor7}) holds. Now,
\begin{align}
&\{\{a,b,c\}_{1},d,e\}_{3}=((a\dashv (b\dashv c))\vdash d)\vdash e=((a\vdash (b\dashv c))\vdash d)\vdash e \nonumber \\
&=\{\{a,b,c\}_{2},d,e\}_{3}=(a\vdash ((b\dashv c)\vdash d))\vdash e=(a\vdash ((b\vdash c)\vdash d))\vdash e \nonumber \\
&=((a\vdash (b\vdash c))\vdash d)\vdash e=(\{a,b,c\}_1\vdash d)\vdash e=\{\{a,b,c\}_{3},d,e\}_{3}, \nonumber
\end{align}
for all $a,b,c,d,e\in D$ and it implies (\ref{cor8}).

Observe now that
\begin{align}
&\{a,\{b,c,d\}_{1},e\}_{1}=(a\dashv ((b\dashv c)\dashv d))\dashv e=((a\dashv (b\dashv c))\dashv d)\dashv e \nonumber \\
&=((a\dashv (b\vdash c))\dashv d)\dashv e=(a\dashv ((b\vdash c)\dashv d))\dashv e=\{a,\{b,c,d\}_{2},e\}_{1} \nonumber \\
&=(a\dashv ((b\vdash c)\vdash d))\dashv e=(a\dashv \{b,c,d\}_3)\dashv e=\{a,\{b,c,d\}_{3},e\}_{1}, \nonumber
\end{align}
for all $a,b,c,d,e\in D$. Thus, (\ref{cor9}) holds.
\begin{align}
&\{a,\{b,c,d\}_{1},e\}_{3}=a\vdash ((b\dashv (c\dashv d))\vdash e)=a\vdash ((b\vdash (c\dashv d))\vdash e) \nonumber \\
&=\{a,\{b,c,d\}_{2},e\}_{3}=a\vdash (b\vdash ((c\dashv d)\vdash e))=a\vdash (b\vdash ((c\vdash d)\vdash e)) \nonumber \\
&=a\vdash((b\vdash (c\vdash d))\vdash e)=a\vdash (\{b,c,d\}_3\vdash e)=\{a,\{b,c,d\}_{3},e\}_{3}, \nonumber
\end{align}
for all $a,b,c,d,e\in D$. We conclude that (\ref{cor10}) is true.

We would like to prove (\ref{cor11}). We have
\begin{align}
&\{a,b,\{c,d,e\}_{1}\}_{1}=a\dashv (b\dashv (c\dashv(d\dashv e)))=a\dashv (b\dashv (c\vdash(d\dashv e))) \nonumber \\
&=\{a,b,\{c,d,e\}_{2}\}_{1}=a\dashv (b\dashv ((c\vdash d)\dashv e))=a\dashv (b\dashv ((c\vdash d)\vdash e)) \nonumber \\
&=a\dashv (b \dashv \{c,d,e\}_3)=\{a,b,\{c,d,e\}_{3}\}_{1}, \nonumber
\end{align}
for all $a,b,c,d,e\in D$.

Finally, we have
\begin{align}
&\{a,b,\{c,d,e\}_{1}\}_{2}=a\vdash (b\dashv(c\dashv(d\dashv e)))=a\vdash (b\dashv(c\vdash(d\dashv e))) \nonumber \\
&=\{a,b,\{c,d,e\}_{2}\}_{2}=a\vdash (b\dashv(c\vdash(d\dashv e)))=a\vdash (b\dashv((c\vdash d)\dashv e)) \nonumber \\
&=a\vdash (b\dashv((c\vdash d)\vdash e))=a\vdash (b\dashv \{c,d,e\}_3)=\{a,b,\{c,d,e\}_{3}\}_{2}, \nonumber
\end{align}
for all $a,b,c,d,e\in D$. Hence, we conclude that (\ref{cor12}) holds.

This completes the proof of the proposition.
\end{proof}

It is known (see \cite{brem1}) that in an associative dialgebra $(D,\vdash,\dashv)$, the first and second Jordan triple products
$\langle a,b,c\rangle _{1}=a\dashv (b\dashv c)+(c\vdash b)\vdash a$, $\langle a,b,c \rangle_{2}=(a\vdash b)\dashv c+(c\vdash b)\dashv a$ give rises to a Jordan triple disystem. From (\ref{cor13}) and (\ref{cor15}) we get 
\begin{align*}
\langle a,b,c\rangle _{1}=\{a,b,c\}_{1}+\{c,b,a\}_{3}, \quad\langle a,b,c\rangle _{2}=\{a,b,c\}_{2}+\{c,b,a\}_{2}, \nonumber
\end{align*}
for all $a,b,c\in D$. That is, an associative dialgebra  generates a structure of Jordan triple disystem on $D$. This fact holds in the general case that we start from an associative triple trisystem of the first kind. In fact, we have the following theorem.

\begin{theorem}\label{th:asstojordan}
    Let $A$ be an associative triple trisystem of the first kind. Then $A$ is a Jordan triple disystem with respect to
    the two trilinear products defined as
    $$
    \langle a,b,c\rangle _{1}=\{a,b,c\}_{1}+\{c,b,a\}_{3}, \quad \langle a,b,c\rangle _{2}=\{a,b,c\}_{2}+\{c,b,a\}_{2}, 
    $$
    for all $a,b,c\in A$.
\end{theorem}
\begin{proof}
    It is clear that $\langle a,b,c\rangle_{2}=\langle
    c,b,a\rangle_{2}$ for all $a,b,c \in A$. Let from now on $a,b,c,d,e$
    be arbitrary elements of $A$. Using (\ref{cor9}) and (\ref{cor10})
    we have (JTD2):
    \begin{align}
    \langle& a,\langle b,c,d\rangle_{1},e\rangle_{1}=\{a,\langle
    b,c,d\rangle_{1},e\}_{1}+\{e,\langle b,c,d\rangle_{1},a\}_{3}
    \nonumber \\
    &=\{a,\{ b,c,d\}_{1},e\}_{1}+\{a,\{ d,c,b\}_{3},e\}_{1} +\{e,\{ b,c,d\}_{1},a\}_{3}+\{e,\{ d,c,b\}_{3},a\}_{3}
    \nonumber \\
    &=\{a,\{ b,c,d\}_{2},e\}_{1}+\{a,\{ d,c,b\}_{2},e\}_{1} +\{e,\{ b,c,d\}_{2},a\}_{3}+\{e,\{ d,c,b\}_{2},a\}_{3}
    \nonumber \\
    &=\langle a,\{ b,c,d\}_{2},e\rangle_{1}+\langle a,\{ d,c,b\}_{2},e\rangle_{1} = \langle a,\langle b,c,d\rangle_{2},e\rangle_{1}. \nonumber
    \end{align}
    
    On other hand, of (\ref{cor8}) and (\ref{cor11}) (JTD3) follows:
    \begin{align}
    \langle &a,b,\langle c,d,e\rangle_{1}\rangle_{1}=\{a,b,\langle
    c,d,e\rangle_{1}\}_{1}+\{\langle c,d,e\rangle_{1},b,a\}_{3}
    \nonumber \\
    &=\{a,b,\{c,d,e\}_{1}\}_{1}+\{a,b,\{e,d,c\}_{3}\}_{1}  +\{\{c,d,e\}_{1},b,a\}_{3}+\{\{e,d,c\}_{3},b,a\}_{3}
    \nonumber \\
    &=\{a,b,\{c,d,e\}_{2}\}_{1}+\{a,b,\{e,d,c\}_{2}\}_{1} +\{\{c,d,e\}_{2},b,a\}_{3}+\{\{e,d,c\}_{2},b,a\}_{3}
    \nonumber \\
    &=\langle a,b,\{c,d,e\}_{2}\rangle_{1}+\langle
    a,b,\{e,d,c\}_{2}\rangle_{1} =\langle a,b,\langle c,d,e\rangle_{2}\rangle_{1}. \nonumber
    \end{align}
    
    Now, (\ref{cor7}) and (\ref{cor12}) imply (JTD4):
    \begin{align}
    \langle &\langle a,b,c\rangle_{1},d,e \rangle_{2}=\{\langle
    a,b,c\rangle_{1},d,e\}_{2}+\{e,d,\langle a,b,c\rangle_{1}\}_{2}
    \nonumber \\
    &=\{\{a,b,c\}_{1},d,e\}_{2}+\{\{c,b,a\}_{3},d,e\}_{2}  +\{e,d,\{a,b,c\}_{1}\}_{2}+\{e,d,\{c,b,a\}_{3}\}_{2}
    \nonumber \\
    &=\{\{a,b,c\}_{2},d,e\}_{2}+\{\{c,b,a\}_{2},d,e\}_{2} +\{e,d,\{a,b,c\}_{2}\}_{2}+\{e,d,\{c,b,a\}_{2}\}_{2}
    \nonumber \\
    &=\langle\{a,b,c\}_{2},d,e\rangle_{2}+\langle\{c,b,a\}_{2},d,e\rangle_{2} =\langle \langle a,b,c\rangle_{2},d,e\rangle_{2}. \nonumber
    \end{align}
    
    Next, we verify the more intricate axioms of a Jordan triple disystem. To do so, first we expand the triple products $\langle \, , \, , \, \rangle$ appearing in the Jordan triple disystem identities in terms of $\{ \, , \, ,\,\}$. Then, we compare terms with the same arrangement of elements, and finally, we use the associative triple trisystem identities to cancel them. Let us show (JTD5):
    
    \begin{align*}\label{cor20}
    \langle &\langle e,d,c\rangle_{1},b,a\rangle_{1}=\{\langle e,d,c\rangle_{1},b,a\}_{1}+\{a,b,\langle e,d,c\rangle_{1}\}_{3} \nonumber \\
    &=\{\{e,d,c\}_{1},b,a\}_{1}+\{\{c,d,e\}_{3},b,a\}_{1} +\{a,b,\{e,d,c\}_{1}\}_{3}+\{a,b,\{c,d,e\}_{3}\}_{3}
    \\
    &=t_{1,1} + t_{1,2} + t_{1,3} + t_{1,4}
    \end{align*}
    \begin{align*}
    \langle & \langle e,b,a\rangle_{1},d,c\rangle_{1}=\{\langle e,b,a\rangle_{1},d,c\}_{1}+\{c,d,\langle e,b,a\rangle_{1}\}_{3} \nonumber \\
    &=\{\{e,b,a\}_{1},d,c\}_{1}+\{\{a,b,e\}_{3},d,c\}_{1} +\{c,d,\{e,b,a\}_{1}\}_{3}+\{c,d,\{a,b,e\}_{3}\}_{3}
    \\
    &=t_{2,1} + t_{2,2} + t_{2,3} + t_{2,4}
    \end{align*}
    \begin{align*}
    &-\langle e,\langle d,a,b\rangle_{1},c\rangle_{1}=-\{e,\langle d,a,b\rangle_{1},c\}_{1}-\{c,\langle d,a,b\rangle_{1},e\}_{3} \nonumber \\
    &=-\{e,\{d,a,b\}_{1},c\}_{1}-\{e,\{b,a,d\}_{3},c\}_{1}  -\{c,\{d,a,b\}_{1},e\}_{3}-\{c,\{b,a,d\}_{3},e\}_{3}
    \\
    &=t_{3,1} + t_{3,2} + t_{3,3} + t_{3,4}
    \end{align*}
    \begin{align*}
    \langle& e,d,\langle c,b,a\rangle_{1}\rangle_{1}=\{e,d,\langle c,b,a\rangle_{1}\}_{1}+\{\langle c,b,a\rangle_{1},d,e\}_{3} \nonumber \\
    &=\{e,d,\{c,b,a\}_{1}\}_{1}+\{e,d,\{a,b,c\}_{3}\}_{1}+\{\{c,b,a\}_{1},d,e\}_{3}+\{\{a,b,c\}_{3},d,e\}_{3}
    \\
    &=t_{4,1} + t_{4,2} + t_{4,3} + t_{4,4}
    \end{align*}
    Being that
    $$
        \begin{array}{llll}
            t_{1,1} = t_{4,1}, & t_{1,2} = t_{2,3}, & t_{1,3} = t_{2,2}, & t_{1,4} = t_{4,4},
            \\
            t_{2,1} = -t_{3,2}, &t_{2,4} = -t_{3,3}, &t_{3,1} = -t_{4,2}, & t_{3,4} = -t_{4,3}.
        \end{array}
    $$
    then, it follows that
    \begin{equation*}
    \langle \langle e,d,c\rangle_{1},b,a\rangle_{1}=\langle \langle e,b,a\rangle_{1},d,c\rangle_{1}
    -\langle e,\langle d,a,b\rangle_{1},c\rangle_{1}+ \langle e,d,\langle c,b,a\rangle_{1}\rangle_{1}.
    \end{equation*}
    
    We continue with (JTD6):
    \begin{align*}
    \langle& \langle e,d,c\rangle_{2},b,a\rangle_{1}=\{\langle e,d,c\rangle_{2},b,a\}_{1}+\{a,b,\langle e,d,c\rangle_{2}\}_{3} \nonumber \\
    &=\{\{e,d,c\}_{2},b,a\}_{1}+\{\{c,d,e\}_{2},b,a\}_{1} +\{a,b,\{e,d,c\}_{2}\}_{3}+\{a,b,\{c,d,e\}_{2}\}_{3} \nonumber
    \\
    &=t_{1,1} + t_{1,2}+ t_{1,3} + t_{1,4}
    \end{align*}
    \begin{align*}
    \langle& \langle e,b,a\rangle_{1},d,c\rangle_{2}=\{\langle e,b,a\rangle_{1},d,c\}_{2}+\{c,d,\langle e,b,a\rangle_{1}\}_{2} \nonumber \\
    &=\{\{e,b,a\}_{1},d,c\}_{2}+\{\{a,b,e\}_{3},d,c\}_{2} +\{c,d,\{e,b,a\}_{1}\}_{2}+\{c,d,\{a,b,e\}_{3}\}_{2}\nonumber
    \\
    &=t_{2,1} + t_{2,2} +t_{2,3} + t_{2,4}
    \end{align*}
    \begin{align*}
    &-\langle e,\langle d,a,b\rangle_{1},c\rangle_{2}=-\{e,\langle d,a,b\rangle_{1},c\}_{2}-\{c,\langle d,a,b\rangle_{1},e\}_{2} \nonumber \\
    &=-\{e,\{d,a,b\}_{1},c\}_{2}-\{e,\{b,a,d\}_{3},c\}_{2}-\{c,\{d,a,b\}_{1},e\}_{2}-\{c,\{b,a,d\}_{3},e\}_{2}
    \\
    &=t_{3,1} + t_{3,2} + t_{3,3} + t_{3,4}
    \end{align*}
    
    Moreover
    \begin{align*}
    \langle& e,d,\langle c,b,a\rangle_{1}\rangle_{2}=\{e,d,\langle c,b,a\rangle_{1}\}_{2}+\{\langle c,b,a\rangle_{1},d,e\}_{2} \nonumber \\
    &=\{e,d,\{c,b,a\}_{1}\}_{2}+\{e,d,\{a,b,c\}_{3}\}_{2} +\{\{c,b,a\}_{1},d,e\}_{2}+\{\{a,b,c\}_{3},d,e\}_{2}
    \\
    &=t_{4,1} + t_{4,2} + t_{4,3} + t_{4,4}
    \end{align*}
    
    Since
    $$
        \begin{array}{llll}
            t_{1,1} = t_{4,1}, & t_{1,2} = t_{2,3}, & t_{1,3} = t_{2,2}, & t_{1,4} = t_{4,4},
            \\
            t_{2,1} = -t_{3,2}, &t_{2,4} = -t_{3,3}, &t_{3,1} = -t_{4,2}, & t_{3,4} = -t_{4,3}.
        \end{array}
    $$
    Therefore, it follows
    \begin{equation*}
    \langle \langle e,d,c\rangle_{2},b,a\rangle_{1}=\langle \langle e,b,a\rangle_{1},d,c\rangle_{2}
    -\langle e,\langle d,a,b\rangle_{1},c\rangle_{2}+ \langle e,d,\langle c,b,a\rangle_{1}\rangle_{2}.
    \end{equation*}
    
    Also, to prove (JTD7)
    \begin{align*}
    \langle& a,b,\langle c,d,e\rangle_{1}\rangle_{1}=\{a,b,\langle c,d,e\rangle_{1}\}_{1}+\{\langle c,d,e\rangle_{1},b,a\}_{3} \nonumber 
    \\
    &= \{a,b,\{c,d,e\}_{1}\}_{1}+\{a,b,\{e,d,c\}_{3}\}_{1}  +\{\{c,d,e\}_{1},b,a\}_{3}+\{\{e,d,c\}_{3},b,a\}_{3} \nonumber
    \\
    &=t_{1,1} + t_{1,2}+ t_{1,3} + t_{1,4}
    \end{align*}
    \begin{align*}
    \langle&\langle a,b,c\rangle_{1},d,e\rangle_{1}=\{\langle a,b,c\rangle_{1},d,e\}_{1}+\{e,d,\langle a,b,c\rangle_{1}\}_{3} \nonumber \\
    &=\{\{a,b,c\}_{1},d,e\}_{1}+\{\{c,b,a\}_{3},d,e\}_{1} +\{e,d,\{a,b,c\}_{1}\}_{3}+\{e,d,\{c,b,a\}_{3}\}_{3} \nonumber
    \\
    &=t_{2,1} + t_{2,2} +t_{2,3} + t_{2,4}
    \end{align*}
    \begin{align*}
    &-\langle c,\langle b,a,d\rangle_{2},e\rangle_{2}=-\{c,\langle b,a,d\rangle_{2},e\}_{2}-\{e,\langle b,a,d\rangle_{2},c\}_{2} \nonumber 
    \\
    &=-\{c,\{b,a,d\}_{2},e\}_{2}-\{c,\{d,a,b\}_{2},e\}_{2} -\{e,\{b,a,d\}_{2},c\}_{2}-\{e,\{d,a,b\}_{2},c\}_{2}
    \\
    &=t_{3,1} + t_{3,2} + t_{3,3} + t_{3,4}
    \end{align*}
    \begin{align*}
    \langle&\langle a,b,e\rangle_{1},d,c\rangle_{1}=\{\langle a,b,e\rangle_{1},d,c\}_{1}+\{c,d,\langle a,b,e\rangle_{1}\}_{3} \nonumber 
    \\
    &=\{\{a,b,e\}_{1},d,c\}_{1}+\{\{e,b,a\}_{3},d,c\}_{1} +\{c,d,\{a,b,e\}_{1}\}_{3}+\{c,d,\{e,b,a\}_{3}\}_{3}
    \\
    &=t_{4,1} + t_{4,2} + t_{4,3} + t_{4,4}
    \end{align*}
    We can check that
    $$
        \begin{array}{llll}
            t_{1,1} = t_{2,1}, & t_{1,2} = t_{4,1}, & t_{1,3} = t_{4,4}, & t_{1,4} = t_{2,4},
            \\
            t_{2,2} = -t_{3,1}, &t_{2,3} = -t_{3,4}, &t_{3,2} = -t_{4,3}, & t_{3,3} = -t_{4,2}.
        \end{array}
    $$
    Thus, we arrive to the following equality
    \begin{equation*}
        \langle a,b,\langle c,d,e\rangle_{1}\rangle_{1}=\langle\langle a,b,c\rangle_{1},d,e\rangle_{1} -\langle c,\langle b,a,d\rangle_{2},e\rangle_{2}+\langle\langle a,b,e\rangle_{1},d,c\rangle_{1}.
    \end{equation*}
    
    Finally, we check the last axiom (JTD8)
    \begin{align*}
        \langle& a,b,\langle c,d,e\rangle_{1}\rangle_{2}=\{a,b,\langle c,d,e\rangle_{1}\}_{2}+\{\langle c,d,e\rangle_{1},b,a\}_{2} \nonumber \\
        &=\{a,b,\{c,d,e\}_{1}\}_{2}+\{a,b,\{e,d,c\}_{3}\}_{2} + \{\{c,d,e\}_{1},b,a\}_{2}+ \{\{e,d,c\}_{3},b,a\}_{2} \nonumber
        \\
        &=t_{1,1} + t_{1,2}+ t_{1,3} + t_{1,4}
    \end{align*}
    \begin{align*}
        \langle &\langle a,b,c\rangle_{2},d,e\rangle_{1}=\{\langle a,b,c\rangle_{2},d,e\}_{1}+\{e,d,\langle a,b,c\rangle_{2}\}_{3} \nonumber \\
        &=\{\{a,b,c\}_{2},d,e\}_{1}+\{\{c,b,a\}_{2},d,e\}_{1} +\{e,d,\{a,b,c\}_{2}\}_{3}+ \{e,d,\{c,b,a\}_{2}\}_{3} \nonumber
        \\
        &=t_{2,1} + t_{2,2} +t_{2,3} + t_{2,4}
    \end{align*}
    \begin{align*}
        &-\langle c,\langle b,a,d\rangle_{1},e\rangle_{2}=-\{c,\langle b,a,d\rangle_{1},e\}_{2}-\{e,\langle b,a,d\rangle_{1},c\}_{2} \nonumber \\
        &=-\{c,\{b,a,d\}_{1},e\}_{2}-\{c,\{d,a,b\}_{3},e\}_{2} -\{e,\{b,a,d\}_{1},c\}_{2}-\{e,\{d,a,b\}_{3},c\}_{2}
        \\
        &=t_{3,1} + t_{3,2} + t_{3,3} + t_{3,4}
        \end{align*}
        \begin{align*}
        \langle &\langle a,b,e\rangle_{2},d,c\rangle_{1}=\{\langle a,b,e\rangle_{2},d,c\}_{1}+\{c,d,\langle a,b,e\rangle_{2}\}_{3} \nonumber \\
        &=\{\{a,b,e\}_{2},d,c\}_{1}+\{\{e,b,a\}_{2},d,c\}_{1} +\{c,d,\{a,b,e\}_{2}\}_{3}+\{c,d,\{e,b,a\}_{2}\}_{3}
        \\
        &=t_{4,1} + t_{4,2} + t_{4,3} + t_{4,4}
    \end{align*}
    
    Let us note that
    $$
        \begin{array}{llll}
            t_{1,1} = t_{2,1}, & t_{1,2} = t_{4,1}, & t_{1,3} = t_{4,4}, & t_{1,4} = t_{2,4},
            \\
            t_{2,2} = -t_{3,1}, &t_{2,3} = -t_{3,4}, &t_{3,2} = -t_{4,3}, & t_{3,3} = -t_{4,2}.
        \end{array}
    $$
    Hence, we conclude 
    \begin{equation*}
        \langle a,b,\langle c,d,e\rangle_{1}\rangle_{2}=\langle \langle a,b,c\rangle_{2},d,e\rangle_{1}-\langle c,\langle b,a,d\rangle_{1},e\rangle_{2} +\langle \langle a,b,e\rangle_{2},d,c\rangle_{1}.
    \end{equation*}
    
    It terminates the proof of the theorem.
\end{proof}

\smallskip

\begin{theorem} \label{th:asstoleibniz}
    Let $A$ be an associative triple trisystem of the first kind. Then $A$ with the new product
    \begin{equation*}
         [a,b,c]=\{a,b,c\}_{1} - \{b,a,c\}_{2} - \{c,a,b\}_{2} +\{c,b,a\}_{3},
    \end{equation*}
    is a Leibniz triple system.
\end{theorem}

\begin{proof}
    We shall check (LTSA) and (LTSB). We have
    \begin{align*}
    [a,b,[c,d,e]]=&[a,b,\{c,d,e\}_{1} - \{d,c,e\}_{2} - \{e,c,d\}_{2} +\{e,d,c\}_{3}] \nonumber \\
    &=\,\,\,\,\,\{a,b,\{c,d,e\}_{1} - \{d,c,e\}_{2} - \{e,c,d\}_{2} +\{e,d,c\}_{3}\}_{1} \nonumber \\
    &\,\,\,\,\,\,-\{b,a,\{c,d,e\}_{1} - \{d,c,e\}_{2} - \{e,c,d\}_{2} +\{e,d,c\}_{3}\}_{2}  \\
    &\,\,\,\,\,\,-\{\{c,d,e\}_{1} - \{d,c,e\}_{2} - \{e,c,d\}_{2} +\{e,d,c\}_{3},a,b\}_{2} \nonumber \\
    &\,\,\,\,\,\,+\{\{c,d,e\}_{1} - \{d,c,e\}_{2} - \{e,c,d\}_{2} +\{e,d,c\}_{3},b,a\}_{3} \nonumber \\
    &=\,\,t_{1,1}+t_{1,2}+\cdots +t_{1,15}+t_{1,16}. \nonumber
    \end{align*}
    Since, for two arbitrary terms, at least two elements occupy
    different places, it follows of the axioms defining any associative triple trisystem of the first kind that
    none of terms can be canceled (this observation
    is also hold for the calculation of the remaining terms).
    
    \begin{align*} 
    -[[a,b,c],d,e]=&-[\{a,b,c\}_{1} - \{b,a,c\}_{2} - \{c,a,b\}_{2} +\{c,b,a\}_{3},d,e] \nonumber \\
    =&-\{\{a,b,c\}_{1} - \{b,a,c\}_{2} - \{c,a,b\}_{2} +\{c,b,a\}_{3},d,e\}_{1} \nonumber \\
    &+\{d,\{a,b,c\}_{1} - \{b,a,c\}_{2} - \{c,a,b\}_{2} +\{c,b,a\}_{3},e\}_{2}  \\
    &+\{e,\{a,b,c\}_{1} - \{b,a,c\}_{2} - \{c,a,b\}_{2} +\{c,b,a\}_{3},d\}_{2} \nonumber \\
    &-\{e,d,\{a,b,c\}_{1} - \{b,a,c\}_{2} - \{c,a,b\}_{2} +\{c,b,a\}_{3}\}_{3} \nonumber \\
    =& \, t_{2,1}+t_{2,2}+\cdots +t_{2,15}+t_{2,16}. \nonumber
    \end{align*}
    
    \begin{align*} 
    [[a,b,d],c,e]&=[\{a,b,d\}_{1} - \{b,a,d\}_{2} - \{d,a,b\}_{2} +\{d,b,a\}_{3},c,e] \nonumber \\
    &=\,\,\,\,\,\{\{a,b,d\}_{1} - \{b,a,d\}_{2} - \{d,a,b\}_{2} +\{d,b,a\}_{3},c,e\}_{1} \nonumber \\
    &\,\,\,\,\,\,-\{c,\{a,b,d\}_{1} - \{b,a,d\}_{2} - \{d,a,b\}_{2} +\{d,b,a\}_{3},e\}_{2}  \\
    &\,\,\,\,\,\,-\{e,\{a,b,d\}_{1} - \{b,a,d\}_{2} - \{d,a,b\}_{2} +\{d,b,a\}_{3},c\}_{2} \nonumber \\
    &\,\,\,\,\,\,+\{e,c,\{a,b,d\}_{1} - \{b,a,d\}_{2} - \{d,a,b\}_{2} +\{d,b,a\}_{3}\}_{3} \nonumber \\
    &=\,\,t_{3,1}+t_{3,2}+\cdots +t_{3,15}+t_{3,16}. \nonumber
    \end{align*}
    
    \begin{align*} 
    -[[a,b,e],d,c]&=-[\{a,b,e\}_{1} - \{b,a,e\}_{2} - \{e,a,b\}_{2} +\{e,b,a\}_{3},d,c] \nonumber \\
    &=-\{\{a,b,e\}_{1} - \{b,a,e\}_{2} - \{e,a,b\}_{2} +\{e,b,a\}_{3},d,c\}_{1} \nonumber \\
    &+\{d,\{a,b,e\}_{1} - \{b,a,e\}_{2} - \{e,a,b\}_{2} +\{e,b,a\}_{3},c\}_{2}  \\
    &+\{c,\{a,b,e\}_{1} - \{b,a,e\}_{2} - \{e,a,b\}_{2} +\{e,b,a\}_{3},d\}_{2} \nonumber \\
    &-\{c,d,\{a,b,e\}_{1} - \{b,a,e\}_{2} - \{e,a,b\}_{2} +\{e,b,a\}_{3}\}_{3} \nonumber \\
    &=\,\,t_{4,1}+t_{4,2}+\cdots +t_{4,15}+t_{4,16}. \nonumber
    \end{align*}
    
    \begin{align*}
    [[a,b,e],c,d]&=[\{a,b,e\}_{1} - \{b,a,e\}_{2} - \{e,a,b\}_{2} +\{e,b,a\}_{3},c,d] \nonumber \\
    &=\,\,\,\,\,\{\{a,b,e\}_{1} - \{b,a,e\}_{2} - \{e,a,b\}_{2} +\{e,b,a\}_{3},c,d\}_{1} \nonumber \\
    &\,\,\,\,\,\,-\{c,\{a,b,e\}_{1} - \{b,a,e\}_{2} - \{e,a,b\}_{2} +\{e,b,a\}_{3},d\}_{2}  \\
    &\,\,\,\,\,\,-\{d,\{a,b,e\}_{1} - \{b,a,e\}_{2} - \{e,a,b\}_{2} +\{e,b,a\}_{3},c\}_{2} \nonumber \\
    &\,\,\,\,\,\,+\{d,c,\{a,b,e\}_{1} - \{b,a,e\}_{2} - \{e,a,b\}_{2} +\{e,b,a\}_{3}\}_{3} \nonumber \\
    &=\,\,t_{5,1}+t_{5,2}+\cdots +t_{5,15}+t_{5,16}. \nonumber
    \end{align*}
    
    Therefore, using this $t_{i,j}$ notation, the left part (LTSA) is equivalent to the sum of these eighty terms $\sum_{i=1}^5\sum_{j=1}^{16} t_{i,j}$. Which is equal to zero because
    $$
        \begin{array}{llll}
            t_{1,1} = -t_{2,1}, & t_{1,2} = -t_{3,1}, & t_{1,3} = -t_{5,1}, & t_{1,4} = - t_{4,1},
            \\
            t_{1,5} = -t_{2,2}, & t_{1,6} = -t_{3,2}, & t_{1,7} = - t_{5,2}, & t_{1,8} = - t_{4,2},
            \\
            t_{1,9} = -t_{4,15}, & t_{1,10} = -t_{5,15}, & t_{1,11} = - t_{3,15}, & t_{1,12} = -t_{2,15},
            \\
            t_{1,13} = - t_{4,16}, & t_{1,14} = -t_{5,16}, & t_{1,15} = -t_{3,16}, & t_{1,16} = -t_{2,16},
            \\
            t_{2,3} = -t_{3,5}, &t_{2,4} = - t_{3,6}, & t_{2,5} = - t_{3,3}, & t_{2,6} = -t_{3,4},
            \\
            t_{2,7} = -t_{5,13}, & t_{2,8} = -t_{5,14}, & t_{2,9} = -t_{5,3}, & t_{2,10} = - t_{5,4},
            \\
            t_{2,11} = -t_{3,13}, & t_{2,12} = - t_{3,14}, & t_{2,13} = -t_{3,11}, & t_{2,14} = - t_{3,12},
            \\
            t_{3,7} = -t_{4,13}, & t_{3,8} = -t_{4,14},&t_{3,9} = -t_{4,3},&t_{3,10} = -t_{4,4},
            \\
            t_{4,5} = - t_{5,9}, & t_{4,6} = - t_{5,10}, & t_{4,7} = - t_{5,11},& t_{4,8} = -t_{5,12},
            \\
            t_{4,9} = -t_{5,5}, & t_{4,10} = -t_{5,6}, & t_{4,11} = - t_{5,7}, & t_{4,12} = - t_{5,8},
        \end{array}
    $$
    where we used the relations of $A$ being an associative triple trisystem of the first kind. Thus, (LTSA) holds. Next, we prove (LTSB) following the same reasoning.
    \begin{align*}
        [a,[b,c,d],e]&=[a,\{b,c,d\}_{1}-\{c,b,d\}_{2}-\{d,b,c\}_{2}+\{d,c,b\}_3,e]   
        \\
        &=\,\,\,\,\, \{a, \{b,c,d\}_{1}-\{c,b,d\}_{2}-\{d,b,c\}_{2}+\{d,c,b\}_3,e\}_1  
        \\
        &\,\,\,\,\, - \{ \{b,c,d\}_{1}-\{c,b,d\}_{2}-\{d,b,c\}_{2}+\{d,c,b\}_3,a,e\}_2  
        \\
        &\,\,\,\,\, - \{e,a,\{b,c,d\}_{1}-\{c,b,d\}_{2}-\{d,b,c\}_{2}+\{d,c,b\}_3\}_2
        \\
        &\,\,\,\,\, + \{e,\{b,c,d\}_{1}-\{c,b,d\}_{2}-\{d,b,c\}_{2}+\{d,c,b\}_3,a\}_3  
        \\
        &= t_{1,1} + t_{1,2} + \cdots + t_{1,15} + t_{1,16}
    \end{align*}
    \begin{align}
        -[[a,b,c],d,e]=& - [\{a,b,c\}_{1} - \{b,a,c\}_{2} - \{c,a,b\}_{2} +\{c,b,a\}_{3}, d, e] \nonumber
        \\
        =&-\{\{a,b,c\}_{1} - \{b,a,c\}_{2} - \{c,a,b\}_{2} +\{c,b,a\}_{3},d,e\}_1 \nonumber 
        \\
        &+\{d,\{a,b,c\}_{1} - \{b,a,c\}_{2} - \{c,a,b\}_{2} +\{c,b,a\}_{3},e \}_2 \nonumber 
        \\
        &+\{e, \{a,b,c\}_{1} - \{b,a,c\}_{2} - \{c,a,b\}_{2} +\{c,b,a\}_{3}, d\}_2 \nonumber 
        \\
        &- \{e,d,\{a,b,c\}_{1} - \{b,a,c\}_{2} - \{c,a,b\}_{2} +\{c,b,a\}_{3}\}_3 \nonumber
        \\
        =&t_{2,1} + t_{2,2} + \cdots + t_{2,15} + t_{2,16} \nonumber
    \end{align}
    \begin{align*}
        [[a,c,b],d,e]&= [\{a,c,b\}_1 - \{c,a,b\}_2-\{b,a,c\}_2 + \{b,c,a\}_3,d,e]
        \\
        &= \,\,\,\,\,\{\{a,c,b\}_1 - \{c,a,b\}_2-\{b,a,c\}_2 + \{b,c,a\}_3,d,e\}_1 
        \\
        &\,\,\,\,\, -\{d,\{a,c,b\}_1 - \{c,a,b\}_2-\{b,a,c\}_2 + \{b,c,a\}_3,e\}_2  
        \\
        &\,\,\,\,\, -\{e, \{a,c,b\}_1 - \{c,a,b\}_2-\{b,a,c\}_2 + \{b,c,a\}_3, d\}_2 
        \\
        &\,\,\,\,\, +\{e,d,\{a,c,b\}_1 - \{c,a,b\}_2-\{b,a,c\}_2 + \{b,c,a\}_3\}_3 
        \\
        &= t_{3,1} + t_{3,2} + \cdots + t_{3,15} + t_{3,16}
    \end{align*}
    \begin{align*}
        [[a,d,b],c, e]&= [\{a,d,b\}_1 - \{d,a,b\}_2 - \{b,a,d\}_2 + \{b,d,a\}_3,c,e]
        \\
        &=\,\,\,\,\,\{\{a,d,b\}_1 - \{d,a,b\}_2 - \{b,a,d\}_2 + \{b,d,a\}_3,c,e\}_1  
        \\
        &\,\,\,\,\, -\{c, \{a,d,b\}_1 - \{d,a,b\}_2 - \{b,a,d\}_2 + \{b,d,a\}_3, e\}_2  
        \\
        &\,\,\,\,\, -\{e,\{a,d,b\}_1 - \{d,a,b\}_2 - \{b,a,d\}_2 + \{b,d,a\}_3,c\}_2 
        \\
        &\,\,\,\,\, +\{e,c,\{a,d,b\}_1 - \{d,a,b\}_2 - \{b,a,d\}_2 + \{b,d,a\}_3\}_3 
        \\
        &= t_{4,1} + t_{4,2} + \cdots + t_{4,15} + t_{4,16}
    \end{align*}
    \begin{align}
        -[[a,d,c],b,e]=& -[ \{a,d,c\}_1 - \{d,a,c\}_2 - \{c,a,d\}_2 + \{c,d,a\}_3, b, e] \nonumber
        \\
        =&- \{\{a,d,c\}_1 - \{d,a,c\}_2 - \{c,a,d\}_2 + \{c,d,a\}_3,b,e\}_1 \nonumber 
        \\
        &+\{b, \{a,d,c\}_1 - \{d,a,c\}_2 - \{c,a,d\}_2 + \{c,d,a\}_3, e\}_2 \nonumber 
        \\
        &+\{e,\{a,d,c\}_1 - \{d,a,c\}_2 - \{c,a,d\}_2 + \{c,d,a\}_3,b\}_2 \nonumber 
        \\
        &-\{e,b,\{a,d,c\}_1 - \{d,a,c\}_2 - \{c,a,d\}_2 + \{c,d,a\}_3\}_3 \nonumber
        \\
        =& t_{5,1} + t_{5,2} + \cdots + t_{5,15} + t_{5,16}\nonumber
    \end{align}
    As before, the left side of the identity (LTSB) is equivalent to the sum of these eighty terms $\sum_{i=1}^5\sum_{j=1}^{16}t_{i,j}$. In this case, these sum is zero because
    $$
        \begin{array}{lllll}
            t_{1,1} = -t_{2,1} , & t_{1,2} = -t_{3,1}, & t_{1,3} = -t_{4,1}, & t_{1,4} = - t_{5,1},& 
            \\
            t_{1,5} = -t_{5,8}, & t_{1,6} = -t_{4,8}, & t_{1,7} = - t_{3,8}, & t_{1,8} = - t_{2,8}, & 
            \\
            t_{1,9} = -t_{2,9}, & t_{1,10} = -t_{3,9}, & t_{1,11} = - t_{4,9}, & t_{1,12} = -t_{5,9},
            \\
            t_{1,13} = - t_{5,16}, & t_{1,14} = -t_{4,16}, & t_{1,15} = -t_{3,16}, & t_{1,16} = -t_{2,16},
            \\
            t_{2,2} = -t_{3,3}, &t_{2,3} = - t_{3,2}, & t_{2,4} = - t_{4,7}, & t_{2,5} = -t_{4,2},
            \\
            t_{2,6} = -t_{3,7}, & t_{2,7} = -t_{3,6}, & t_{2,10} = -t_{3,11}, & t_{2,11} = - t_{3,10},
            \\
            t_{2,12} = -t_{4,15}, & t_{2,13} = - t_{4,10}, & t_{2,14} = -t_{3,15}, & t_{2,15} = - t_{3,14},
            \\
            t_{3,4} = -t_{5,7}, & t_{3,5} = -t_{5,2},&t_{3,12} = -t_{5,15},&t_{3,13} = -t_{5,10},
            \\
            t_{4,3} = - t_{5,5}, & t_{4,4} = - t_{5,6}, & t_{4,5} = - t_{5,3},& t_{4,6} = -t_{5,4},
            \\
            t_{4,11} = -t_{5,13}, & t_{4,12} = -t_{5,14}, & t_{4,13} = - t_{5,11}, & t_{4,14} = - t_{5,12}.
        \end{array}
    $$
    So (LTSB) holds. Hence, $A$ with the triple product $[ \, , \, ,\,]$ is a Leibniz triple system.
\end{proof}

\section{Associative triple systems of second kind}

In this section, we apply the KP algorithm to the identities
$$
    \{\{a,b,c\},d,e\} \equiv \{a,\{d,c,b\}, e\} \equiv \{a,b,\{c,d,e\}\},
$$
which define an associative triple system of the second kind. We then introduce the definition of an associative triple trisystem of the second kind.

To study these structures, it is natural to consider the role of an involution. For this reason, we recall the notion of an involution (of the second type) in an associative dialgebra. Finally, we prove that an associative dialgebra with involution naturally gives rise to an associative triple trisystem of the second kind in the expected way.

The Part $1$ of KP algorithm yields the following $5$ relations:
\begin{align*}
        \{\{a,b,c\}_{1},d,e\}_{1}=\{a,\{d,c,b\}_{1},e\}_{1}= \{a,b,\{c,d,e\}_{1}\}_{1}, 
        \\
        \{\{a,b,c\}_{2},d,e\}_{1}= \{a,\{d,c,b\}_{3},e\}_{2}= \{a,b,\{c,d,e\}_{1}\}_{2}, 
        \\
        \{\{a,b,c\}_{3},d,e\}_{1}= \{a,\{d,c,b\}_{2},e\}_{2}= \{a,b,\{c,d,e\}_{1}\}_{3}, 
        \\
        \{\{a,b,c\}_{3},d,e\}_{2}= \{a,\{d,c,b\}_{1},e\}_{2}= \{a,b,\{c,d,e\}_{2}\}_{3}, 
        \\
        \{\{a,b,c\}_{3},d,e\}_{3}= \{a,\{d,c,b\}_{3},e\}_{3}= \{a,b,\{c,d,e\}_{3}\}_{3}, 
\end{align*}
the Part $2$ gives the remaining $6$ relations:
\begin{align*}
    \{\{a,b,c\}_{1},d,e\}_{2}= \{\{a,b,c\}_{2},d,e\}_{2}= \{\{a,b,c\}_{3},d,e\}_{2}, 
    \\
    \{\{a,b,c\}_{1},d,e\}_{3}= \{\{a,b,c\}_{2},d,e\}_{3}= \{\{a,b,c\}_{3},d,e\}_{3}, 
    \\
    \{a,\{b,c,d\}_{1},e\}_{1}= \{a,\{b,c,d\}_{2},e\}_{1}= \{a,\{b,c,d\}_{3},e\}_{1}, 
    \\
    \{a,\{b,c,d\}_{1},e\}_{3}= \{a,\{b,c,d\}_{2},e\}_{3}= \{a,\{b,c,d\}_{3},e\}_{3}, 
    \\
    \{a,b,\{c,d,e\}_{1}\}_{1}= \{a,b,\{c,d,e\}_{2}\}_{1}= \{a,b,\{c,d,e\}_{3}\}_{1}, 
    \\
    \{a,b,\{c,d,e\}_{1}\}_{2}= \{a,b,\{c,d,e\}_{2}\}_{2}= \{a,b,\{c,d,e\}_{3}\}_{2}.
\end{align*}
Therefore, we define an associative triple trisystem of the second kind if these eleven relations hold.
\begin{definition}\label{cor1b}
    Let $A$ be a $\phi$-module. Three trilinear maps $\{ \,, \, ,\}_{i}:\,\,A\times A\times A\rightarrow A$, with $i \in \{ 1,2,3\}$,
    endow $A$ with a structure of \textbf{associative triple trisystem of the second kind} if these satisfy the following relations:
    \begin{align}
        &\label{cor2b} \{\{a,b,c\}_{1},d,e\}_{1}=\{a,\{d,c,b\}_{1},e\}_{1}= \{a,b,\{c,d,e\}_{1}\}_{1}, 
        \\
        &\label{cor3b} \{\{a,b,c\}_{2},d,e\}_{1}= \{a,\{d,c,b\}_{3},e\}_{2}= \{a,b,\{c,d,e\}_{1}\}_{2}, 
        \\
        &\label{cor4b} \{\{a,b,c\}_{3},d,e\}_{1}= \{a,\{d,c,b\}_{2},e\}_{2}= \{a,b,\{c,d,e\}_{1}\}_{3}, 
        \\
        &\label{cor5b} \{\{a,b,c\}_{3},d,e\}_{2}= \{a,\{d,c,b\}_{1},e\}_{2}= \{a,b,\{c,d,e\}_{2}\}_{3}, 
        \\
        &\label{cor6b} \{\{a,b,c\}_{3},d,e\}_{3}= \{a,\{d,c,b\}_{3},e\}_{3}= \{a,b,\{c,d,e\}_{3}\}_{3}, 
        \\
        &\label{cor7b} \{\{a,b,c\}_{1},d,e\}_{2}= \{\{a,b,c\}_{2},d,e\}_{2}= \{\{a,b,c\}_{3},d,e\}_{2}, 
        \\
        &\label{cor8b} \{\{a,b,c\}_{1},d,e\}_{3}= \{\{a,b,c\}_{2},d,e\}_{3}= \{\{a,b,c\}_{3},d,e\}_{3}, 
        \\
        &\label{cor9b} \{a,\{b,c,d\}_{1},e\}_{1}= \{a,\{b,c,d\}_{2},e\}_{1}= \{a,\{b,c,d\}_{3},e\}_{1}, 
        \\
        &\label{cor10b} \{a,\{b,c,d\}_{1},e\}_{3}= \{a,\{b,c,d\}_{2},e\}_{3}= \{a,\{b,c,d\}_{3},e\}_{3}, 
        \\
        &\label{cor11b} \{a,b,\{c,d,e\}_{1}\}_{1}= \{a,b,\{c,d,e\}_{2}\}_{1}= \{a,b,\{c,d,e\}_{3}\}_{1}, 
        \\
        &\label{cor12b} \{a,b,\{c,d,e\}_{1}\}_{2}= \{a,b,\{c,d,e\}_{2}\}_{2}= \{a,b,\{c,d,e\}_{3}\}_{2}.
    \end{align}
    for all $a,b,c,d,e\in A$.
\end{definition}

In a similar way as in the first kind, let $A$ be an associative triple trisystem of the second kind. Define $A^{ann}$ as the subspace spanned by the elements of the form $\{a,b,c\}_1 - \{a,b,c\}_2$ and $\{a,b,c\}_1 - \{a,b,c\}_3$ for all $a,b,c \in A$. Then $A$ is an associative triple system of the second kind if and only if $A^{ann} = \{0\}$. Moreover, suppose that $A = A^{ann} \oplus \mathcal{A}$ as subspaces where $\mathcal{A}$ is closed with respect to the trilinear maps $\{ \ , \ , \}_i$ for $i \in \{1,2,3\}$, then $\mathcal{A}$ is an associative triple system of the second kind.

\begin{definition}(\cite{Felipe})
    Let $D$ be an associative dialgebra. A linear map $* \colon D \rightarrow D$ is an \textbf{involution of D}  if $(a^*)^* = a$, $(a \dashv b)^* = b^* \vdash a^*$ and $(a \vdash b)^* = b^* \dashv a^*$, for all $a,b \in D$.
\end{definition}

\begin{proposition}
    Let $D$ be an associative dialgebra and $*$ an involution of $D$. Then $D$ with the following triple products
    $$
        \{a,b,c\}_1 := a \dashv (b^* \dashv c), \quad \{a,b,c\}_2 := a \vdash (b^* \dashv c), \quad \{a,b,c\}_3 := a \vdash (b^* \vdash c),
    $$
    is an associative triple trisystem of the second kind.
\end{proposition}
\begin{proof}
    First, as expected from the behavior of the involution, we observe that when we apply the involution to these triple products, the products $\dashv$ and $\vdash$ are reversed:
    \begin{align*}
        \{a,b,c\}_1 ^* &= (a \dashv (b^* \dashv c))^* = (b^* \dashv c)^* \vdash a^* = (c^* \vdash b) \vdash a^* = \{c^*,b^*,a^*\}_3,
        \\
        \{a,b,c\}_2^* &= (a \vdash (b^* \dashv c))^* = (b^* \dashv c)^* \dashv a^* = (c^* \vdash b)\dashv a^* = \{c^*,b^*,a^*\}_2,
        \\
        \{a,b,c\}_3^* &= (a \vdash (b^* \vdash c))^* = (b^* \vdash c)^* \dashv a^* = (c^* \dashv b) \dashv a^* = \{c^*,b^*,a^*\}_1.
    \end{align*}
    We prove explicitly the first five relations (\ref{cor2b})-(\ref{cor6b}). The remaining six relations (\ref{cor7b})-(\ref{cor12b}) are following from the same computations as in Proposition \ref{th:asstoass}.
    Thus,
    \begin{align*}
        \{\{a,b,c\}_1,d,e\}_1 &=  (a \dashv (b^* \dashv c)) \dashv (d^* \dashv e) = a \dashv ((b^* \dashv c) \dashv (d^* \dashv e)) 
        \\
        &=a\dashv (b^* \dashv (c \dashv (d^* \dashv e))) = a \dashv (b^* \dashv \{c,d,e\}_1) = \{a,b,\{c,d,e\}_1\}_1
        \\
        &= a \dashv ((b^* \dashv c) \vdash (d^* \dashv e)) = a \dashv ((b^* \vdash c) \vdash (d^* \dashv e))
        \\
        &= a \dashv (((b^* \vdash c) \vdash d^*)\dashv e) = a \dashv (\{d,c,b\}_1^* \dashv e) = \{a,\{d,c,b\}_1,e\}_1
    \end{align*}
    \begin{align*}
        \{\{a,b,c\}_2,d,e\}_1 &= ((a \vdash (b^* \dashv c)) \dashv d^*) \dashv e = (a \vdash ((b^* \dashv c) \dashv d^*)) \dashv e
        \\
        &= (a \vdash (d \vdash (c^* \vdash b)  )^*) \dashv e =( a \vdash \{d,c,b\}_3^*) \dashv e = \{a,\{d,c,b\}_3,e\}_2 
        \\
        &= a \vdash (((b^* \dashv c) \dashv d^*) \dashv e)) = a \vdash ( b^* \dashv (c \dashv (d^* \dashv e)))
        \\
        &=a \vdash (b^* \dashv \{c,d,e\}_1) = \{a,b,\{c,d,e\}_1\}_2
    \end{align*}
    \begin{align*}
        \{\{a,b,c\}_3,d,e\}_1 &= (a \vdash (b^* \vdash c)) \dashv (d^* \dashv e) = a \vdash ((b^* \vdash c) \dashv (d^* \dashv e)) 
        \\
        &= a\vdash (\{d,c,b\}_2^*\dashv e ) = \{a,\{d,c,b\}_2,e\}_2
        \\
        &=a \vdash (b^* \vdash (c \dashv (d^* \dashv e))) =a \vdash (b^* \vdash \{c,d,e\}_3) =\{a,b,\{c,d,e\}_1\}_3
    \end{align*}
    \begin{align*}
        \{\{a,b,c\}_3,d,e\}_2 &= ((a \vdash (b^* \vdash c)) \vdash d^* )\dashv e = (a \vdash ((b^* \vdash c) \vdash d^*))\dashv e 
        \\
        &=(a \vdash \{d,c,b\}_1^*)\dashv e = \{a,\{d,c,b\}_1,e\}_2
        \\
        &=(a \vdash (b^* \vdash c)) \vdash (d^* \dashv e) = a \vdash ((b^* \vdash c) \vdash (d^*\dashv e))
        \\
        &=a \vdash (b^* \vdash (c \vdash (d^* \dashv e))) = a \vdash (b^* \vdash \{c,d,e\}_2) = \{a,b,\{c,d,e\}_2\}_3
    \end{align*}
    \begin{align*}
        \{\{a,b,c\}_3,d,e\}_3 &= (a \vdash (b^* \vdash c)) \vdash (d^* \vdash e) = a \vdash ((b^* \vdash c) \vdash (d^* \vdash e))
        \\
        &= a \vdash ((b^* \dashv c) \vdash (d^* \vdash e)) = a \vdash (((b^* \dashv c) \vdash d^*)\vdash e)
        \\
        &= a\vdash (((b^* \dashv c) \dashv d^* ) \vdash e) = a \vdash ( \{d,c,b\}_3^* \vdash e) = \{a,\{d,c,b\}_3,e\}_3
        \\
        &=a \vdash (b^* \vdash (c \vdash (d^* \vdash e))) = a \vdash (b^* \vdash \{c,d,e\}_3) = \{a,b,\{c,d,e\}_3\}_3
    \end{align*}
\end{proof}

\begin{example}
    Let $R$ be a ring with involution $\bar{x}$. For any $A \in M_{m\times n}(R)$ let us introduce the following notation. Let us pick $m_1 <  m$ and $n_1 < n$, then 
    $$
        A = \left( \begin{array}{c|c}f_A &u_A \\ \hline l_A&d_A \end{array}\right) \in M_{m \times n}^{m_1 \times n_1} (R), 
    $$
    where $d_A \in M_{m_1 \times n_1}(R), u_A \in M_{(m-m_1) \times n_1}(R), l_A \in M_{m_1 \times (n-n_1)}$ and $f_A \in M_{(m-m_1) \times (n-n_1)}(R)$. We denote by $*$ the conjugate transpose of any matrix $A$.
    
    Let us pick $A \in M_{m \times n}^{m_1 \times n_1}(R)$ and $B \in M_{n \times p}^{n_1\times p_1 }(R)$. We define:
    $$
        A \dashv B := \left( \begin{array}{c|c}0 & u_A\\ \hline 0 & d_A \end{array}\right)\left( \begin{array}{c|c}0 & u_B\\ \hline 0 & d_B \end{array}\right) = \left( \begin{array}{c|c}0 & u_Ad_B\\ \hline 0 & d_A d_B\end{array}\right)
    $$
    and
    $$
        A \vdash B :=\left( \begin{array}{c|c}0 & 0\\ \hline l_A & d_A \end{array}\right)\left( \begin{array}{c|c}0 & 0\\ \hline l_B & d_B \end{array}\right) = \left( \begin{array}{c|c}0 & 0\\ \hline d_Al_B & d_A d_B\end{array}\right).
    $$
    Observe that
    \begin{align*}
        (A \dashv B)^* &=\left( \begin{array}{c|c}0 & u_Ad_B\\ \hline 0 & d_A d_B\end{array}\right)^*  = \left( \begin{array}{c|c}0 & 0\\ \hline d_B^*u_A^*&  d_B^*d_A^*\end{array}\right) 
        \\
        &=\left( \begin{array}{c|c}0 & 0\\ \hline u_B^* & d_B^*\end{array}\right)\left( \begin{array}{c|c}0 & 0\\ \hline u_A^* & d_A^* \end{array}\right)=B^* \vdash A^*,
    \end{align*}
    and
    \begin{align*}
        (A \vdash B)^* &= \left( \begin{array}{c|c}0 & 0\\ \hline d_Al_B & d_A d_B\end{array}\right) ^* = \left( \begin{array}{c|c}0 & l_B^*d_A^*\\ \hline 0& d_B^*d_A^*\end{array}\right)
        \\
        &=\left( \begin{array}{c|c}0 & l_B^*\\ \hline 0& d_B^*\end{array}\right)\left( \begin{array}{c|c}0 & l_A^*\\ \hline 0& d_A^*\end{array}\right) = B^* \dashv A^*.
    \end{align*}
    
    \begin{enumerate}
        \item If $ m= n = p$ and $m_1 = n_1 = p_1$, then $(M_m^{m_1}(R), \dashv, \vdash)$ is an associative dialgebra. Moreover, $*$ is an involution of $M_m^{m_1}(R)$ as a dialgebra. Indeed, let us consider $A,B,C \in M_{m\times m}^{m_1 \times m_1}(R)$.  First, let us show that $(A \dashv B) \dashv C = A \dashv (B \dashv C) = A \dashv (B \vdash C)$:
        \begin{align*}
            (A \dashv B) \dashv C &= \left(\left(  \begin{array}{c|c}0 & u_A\\ \hline 0& d_A \end{array}  \right) \left(\begin{array}{c|c}0 & u_B \\ \hline 0 & d_B  \end{array}\right)\right) \dashv C  = \left(  \begin{array}{c|c}0 & u_Ad_B \\ \hline 0 & d_Ad_B  \end{array}  \right) \dashv C 
            \\
            &=\left(  \begin{array}{c|c}0 & u_Ad_B \\ \hline 0 & d_Ad_B  \end{array}  \right)\left(  \begin{array}{c|c}0 & u_C\\ \hline 0 & d_C  \end{array}  \right) = \left(  \begin{array}{c|c}0 & u_Ad_Bd_C \\ \hline 0 & d_Ad_Bd_C \end{array}  \right),
        \end{align*}
        \begin{align*}
            A \dashv (B \dashv C) &= A \dashv \left(\left(  \begin{array}{c|c}0 & u_B\\ \hline 0& d_B \end{array}  \right) \left(\begin{array}{c|c}0 & u_C \\ \hline 0 & d_C  \end{array}\right)\right) =  A \dashv \left(\begin{array}{c|c}0 & u_Bd_C \\ \hline 0 & d_Bd_C  \end{array}\right)
            \\
            &=\left(\begin{array}{c|c}0 & u_A \\ \hline 0 & d_A  \end{array}\right)\left(\begin{array}{c|c}0 & u_Bd_C \\ \hline 0 & d_Bd_C  \end{array}\right) = \left(\begin{array}{c|c}0 & u_Ad_Bd_C \\ \hline 0 & d_Ad_Bd_C  \end{array}\right),
        \end{align*}
        \begin{align*}
            A \dashv (B \vdash C) &= A \dashv \left( \left( \begin{array}{c|c}0 & 0\\ \hline l_B  & d_B  \end{array}\right)\left( \begin{array}{c|c}0 & 0\\ \hline l_C & d_C \end{array}\right)\right) = A \dashv \left( \begin{array}{c|c}0 & 0\\ \hline d_Bl_C & d_Bd_C \end{array}\right)
            \\
            &=\left( \begin{array}{c|c}0 & u_A\\ \hline 0 & d_A \end{array}\right)\left( \begin{array}{c|c}0 & 0\\ \hline 0 & d_Bd_C \end{array}\right) = \left( \begin{array}{c|c}0 & u_Ad_Bd_C\\ \hline 0 & d_Ad_Bd_C \end{array}\right).
        \end{align*}
        
        Second, we prove $(A \vdash B) \dashv C = A \vdash (B \dashv C)$:
        \begin{align*}
            (A \vdash B) \dashv C & = \left(\left(  \begin{array}{c|c}0 & 0\\ \hline l_A & d_A \end{array}  \right) \left(\begin{array}{c|c}0 & 0\\ \hline l_B  & d_B  \end{array}\right)\right) \dashv C  = \left(  \begin{array}{c|c}0 & 0\\ \hline d_Al_B  & d_Ad_B  \end{array}  \right) \dashv C 
            \\
            &=\left(  \begin{array}{c|c}0 & 0\\ \hline 0  & d_Ad_B  \end{array}  \right)\left(  \begin{array}{c|c}0 & u_C\\ \hline 0  & d_C  \end{array}  \right) = \left(  \begin{array}{c|c}0 & 0\\ \hline 0  & d_Ad_Bd_C  \end{array}  \right),
        \end{align*}
        \begin{align*}
            A \vdash (B  \dashv C) &= A \vdash \left( \left( \begin{array}{c|c}0 & u_B \\ \hline 0 & d_B  \end{array}\right)\left( \begin{array}{c|c}0 & u_C\\ \hline 0 & d_C \end{array}\right)\right) = A \vdash \left( \begin{array}{c|c}0 & u_B d_C\\ \hline 0 & d_B d_C \end{array}\right) 
            \\
            &= \left( \begin{array}{c|c}0 & 0\\ \hline l_A & d_A \end{array}\right)\left( \begin{array}{c|c}0 & 0\\ \hline 0 & d_B d_C \end{array}\right) = \left( \begin{array}{c|c}0 & 0\\ \hline 0 & d_Ad_B d_C \end{array}\right). 
        \end{align*}
        Finally, we prove $(A \vdash B) \vdash C) = (A \dashv B) \vdash C = A \vdash (B \vdash C)$:
        \begin{align*}
            (A \vdash B ) \vdash C &= \left(\left(  \begin{array}{c|c}0 & 0\\ \hline l_A & d_A \end{array}  \right) \left(\begin{array}{c|c}0 & 0\\ \hline l_B  & d_B  \end{array}\right)\right) \vdash C  = \left(  \begin{array}{c|c}0 & 0\\ \hline d_Al_B  & d_Ad_B  \end{array}  \right) \vdash C 
            \\
            &=\left(  \begin{array}{c|c}0 & 0\\ \hline d_Al_B  & d_Ad_B  \end{array}  \right) \left(  \begin{array}{c|c}0 & 0\\ \hline l_C & d_C \end{array}  \right)  = \left(  \begin{array}{c|c}0 & 0\\ \hline d_Ad_B l_C & d_Ad_B d_C \end{array}  \right) ,
        \end{align*}
        \begin{align*}
            (A \dashv B ) \vdash C &= \left(\left(  \begin{array}{c|c}0 & u_A\\ \hline 0& d_A \end{array}  \right) \left(\begin{array}{c|c}0 & u_B \\ \hline 0 & d_B  \end{array}\right)\right) \vdash C  = \left(  \begin{array}{c|c}0 & u_Ad_B \\ \hline 0 & d_Ad_B  \end{array}  \right) \vdash C 
            \\
            &=\left(  \begin{array}{c|c}0 & 0\\ \hline 0 & d_Ad_B  \end{array}  \right) \left(  \begin{array}{c|c}0 & 0\\ \hline l_C & d_C \end{array}  \right)  = \left(  \begin{array}{c|c}0 & 0\\ \hline d_Ad_B l_C & d_Ad_B d_C \end{array}  \right) ,
        \end{align*}
        \begin{align*}
            A \vdash (B  \vdash C) &= A \vdash \left( \left( \begin{array}{c|c}0 & 0\\ \hline l_B  & d_B  \end{array}\right)\left( \begin{array}{c|c}0 & 0\\ \hline l_C & d_C \end{array}\right)\right) = A \vdash \left( \begin{array}{c|c}0 &0\\ \hline d_B l_C & d_B d_C \end{array}\right) 
            \\
            &= \left( \begin{array}{c|c}0 & 0\\ \hline l_A & d_A \end{array}\right)\left( \begin{array}{c|c}0 & 0\\ \hline d_B l_C & d_B d_C \end{array}\right) = \left( \begin{array}{c|c}0 & 0\\ \hline d_Ad_B l_C & d_Ad_B d_C \end{array}\right) .
        \end{align*}
        
        \item The rectangle matrices $M_{m\times n}^{m_1\times n_1}(R)$ with 
        \begin{align*}
            \{A,B,C\}_1 &:= A \dashv ( B^* \dashv C)= \left( \begin{array}{c|c}0 & u_Ad_B^*d_C\\ \hline 0 & d_A d_B^*d_C\end{array}\right),
            \\
            \{A,B,C\}_2 &:= A \vdash (B^* \dashv C)= \left( \begin{array}{c|c}0 & 0\\ \hline 0 & d_A d_B^*d_C\end{array}\right), \\
            \{A,B,C\}_3& := A \vdash (B^* \vdash C)=\left( \begin{array}{c|c}0 & 0\\ \hline d_Ad_B^*l_C & d_A d_B^*d_C\end{array}\right),
        \end{align*}
        is an associative triple trisystem of the second kind that is not arised from an associative dialgebra. We show (\ref{cor3b}), the rest are similar.
        \begin{align*}
            \{\{A,B,C\}_2,D,E\}_1 &= \left\lbrace \left( \begin{array}{c|c}0 & 0\\ \hline 0 & d_A d_B^*d_C\end{array}\right),D,E\right\rbrace_1 = \left( \begin{array}{c|c}0 & 0\\ \hline 0 & d_A d_B^*d_Cd_D^*d_E\end{array}\right) 
            \\
            \\
            \{A,\{D,C,B\}_3,E\}_2 &= \left\lbrace A, \left( \begin{array}{c|c}0 & 0\\ \hline d_Dd_C^*l_B& d_D d_C^*d_B\end{array}\right), E\right\rbrace_2  
            \\
            &= \left( \begin{array}{c|c}0 & 0\\ \hline 0 & d_A (d_dd_C^*d_B)^*d_E\end{array}\right) = \left( \begin{array}{c|c}0 & 0\\ \hline 0 & d_A d_B^*d_Cd_D^*d_E\end{array}\right) 
            \\
            \\
            \{A,B,\{C,D,E\}_1\}_2 &= \left\lbrace A,B,\left( \begin{array}{c|c}0 & u_Cd_D^*d_E\\ \hline 0 & d_C d_D^*d_E\end{array}\right) \right\rbrace= \left( \begin{array}{c|c}0 & 0\\ \hline 0 & d_Ad_B^*d_C d_D^*d_E\end{array}\right)
        \end{align*}
    \end{enumerate}
\end{example}

We would like to remark that this example is the first on the literature, as far as we know, of non-trivial associative dialgebras on matrices. That allows, in particular, to provide a fruitful source of examples of finite dimensional Leibniz algebras. We just need to consider the Leibniz algebra $M_m^{m_1}(R)^{(-)}$.

\begin{example}
    Suppose $m=n = 2$ and $m_1=n_1 =1$. Then we consider $L=M^{1}_2(\mathbb{C})^{(-)}$, where 
    $$
        [A,B] := A\dashv B - B \vdash A  = \begin{pmatrix}0 & a_{12}b_{22}\\ -b_{22}a_{21} & 0\end{pmatrix}
    $$
    for any $A = \begin{pmatrix}a_{11} & a_{12} \\ a_{21} & a_{22} \end{pmatrix}$ and $B = \begin{pmatrix} b_{11} & b_{12} \\ b_{21} & b_{22}\end{pmatrix}$ in $L$. 
    Then $E_1=\begin{pmatrix}
        1&0\\0&0
    \end{pmatrix}$, $E_2= \begin{pmatrix}
        0&1\\0&0
    \end{pmatrix}$ $E_3= \begin{pmatrix}
        0 & 0\\1&0
    \end{pmatrix}$ and $ X= \begin{pmatrix}
        0&0\\0&1
    \end{pmatrix}$ form a basis that satisfies
    $$
        [E_1,X] = 0 \cdot E_1 =0, \quad [E_2,X] = E_2, \quad [E_3,X] =E_3.
    $$
    Hence, $L \cong \mathcal{L}_{23}(0,1)$ (see \cite{canete}). 

    Furthermore, we can consider the subspace $L'$ with basis $\mathcal{B} = \left\lbrace \begin{pmatrix}
        0&1\\0&0
    \end{pmatrix}, \ \begin{pmatrix}
        0 & 0\\1&0
    \end{pmatrix}, \ \begin{pmatrix}
        0&1\\0&1
    \end{pmatrix} \right\rbrace $ that satisfies 
    $$
        [B_1, B_3] = B_1, \quad [B_2,B_3] = B_2, \quad [B_3,B_3] = B_1.
    $$
    Then $L' \cong \mathrm{AO3}(\alpha_1=1, \alpha_2=0)$ (see \cite{ayupov} and \cite{casas}).

    In general, it allows us to construct a family of non-Lie Leibniz algebras of dimension $m^2$ and $m_1(2m-m_1)$, for any $1\leq m_1 < m$, and its respective subalgebras of lower dimensions.
\end{example}

\begin{theorem}
    Let $A$ be an associative triple trisystem of the second kind. Then the triple products
    $$
        \langle a,b,c\rangle_1 := \{a,b ,c\}_1 + \{c,b ,a\}_3, \quad \langle a,b,c \rangle_2 := \{a,b ,c\}_2 + \{c,b ,a\}_2,
    $$
    endow $A$ with a Jordan triple disystem structure.
\end{theorem}
\begin{proof}
    The axiom (JTD1) is clear from the definition. Since axioms (\ref{cor7})-(\ref{cor12}) are the same as (\ref{cor7b})-(\ref{cor12b}), (JTD2), (JTD3) and (JTD4) also hold. Now, following the notation introduced in Theorem \ref{th:asstojordan}, we show (JTD5) and (JTD6) by noticing, for each case:
    $$
        \begin{array}{llll}
            t_{1,1} = t_{4,1}, & t_{1,2} = t_{2,3}, & t_{1,3} = t_{2,2}, & t_{1,4} = t_{4,4},
            \\
            t_{2,1} = -t_{3,1},& t_{2,4}=-t_{3,4}, &t_{3,2} = -t_{4,2}, &t_{3,3} = - t_{4,3}.
        \end{array}
    $$
    Finally, (JTD7) and (JTD8) are satisfied because, for each case:
    $$
        \begin{array}{llll}
            t_{1,1} = t_{2,1}, & t_{1,2} = t_{4,1}, & t_{1,3} = t_{4,4}, & t_{1,4} = t_{2,4},
            \\
            t_{2,2} = -t_{3,2},& t_{2,3}=-t_{3,3}, &t_{3,1} = -t_{4,3}, &t_{3,4} = - t_{4,2}.
        \end{array}
    $$
\end{proof}

\begin{theorem}
    Let $A$ be an associative triple trisystem of the second kind. Then the triple product
    $$
        [a,b,c] := \{a,b,c\}_1 - \{b,a,c\}_2 - \{c,a,b\}_2 + \{c,b,a\}_3,
    $$
    is a Leibniz triple system.
\end{theorem}
\begin{proof}
    As in the preovious theorem, we refer to the proof of Theorem \ref{th:asstoleibniz}. To prove (LTSA), following the notation introduced in Theorem \ref{th:asstoleibniz}, it is sufficient to see that 
    $$
        \begin{array}{llll}
            t_{1,1} = -t_{2,1}, & t_{1,2} = -t_{3,1}, & t_{1,3} = -t_{5,1}, & t_{1,4} = - t_{4,1},
            \\
            t_{1,5} = -t_{2,2}, & t_{1,6} = -t_{3,2}, & t_{1,7} = - t_{5,2}, & t_{1,8} = - t_{4,2},
            \\
            t_{1,9} = -t_{4,15}, & t_{1,10} = -t_{5,15}, & t_{1,11} = - t_{3,15}, & t_{1,12} = -t_{2,15},
            \\
            t_{1,13} = - t_{4,16}, & t_{1,14} = -t_{5,16}, & t_{1,15} = -t_{3,16}, & t_{1,16} = -t_{2,16},
            \\
            t_{2,3} = -t_{3,8}, &t_{2,4} = - t_{3,7}, & t_{2,5} = - t_{5,14}, &t_{2,6} = -t_{5,13},
            \\
            t_{2,7} = -t_{3,4}, & t_{2,8} = -t_{3,3}, & t_{2,9} = -t_{3,14}, & t_{2,10} = - t_{3,13},
            \\
            t_{2,11} = -t_{5,4}, & t_{2,12} = - t_{5,3}, & t_{2,13} = -t_{3,10}, & t_{2,14} = - t_{3,9},
            \\
            t_{3,5} = -t_{4,14}, & t_{3,6} = -t_{4,13},&t_{3,11} = -t_{4,4},&t_{3,12} = -t_{4,3},
            \\
            t_{4,5} = - t_{5,9}, & t_{4,6} = - t_{5,10}, & t_{4,7} = - t_{5,11},& t_{4,8} = -t_{5,12},
            \\
            t_{4,9} = -t_{5,5}, & t_{4,10} = -t_{5,6}, & t_{4,11} = - t_{5,7}, & t_{4,12} = - t_{5,8}.
        \end{array}
    $$
    where we had to modify the rows five, six, seven and eight. And, to prove (LTSB):
    $$
        \begin{array}{lllll}
            t_{1,1} = -t_{5,1} , & t_{1,2} = -t_{4,1}, & t_{1,3} = -t_{3,1}, & t_{1,4} = - t_{2,1},& 
            \\
            t_{1,5} = -t_{5,5}, & t_{1,6} = -t_{4,5}, & t_{1,7} = - t_{3,5}, & t_{1,8} = - t_{2,5}, & 
            \\
            t_{1,9} = -t_{2,12}, & t_{1,10} = -t_{3,12}, & t_{1,11} = - t_{4,12}, & t_{1,12} = -t_{5,12},
            \\
            t_{1,13} = - t_{2,16}, & t_{1,14} = -t_{3,16}, & t_{1,15} = -t_{4,16}, & t_{1,16} = -t_{5,16},
            \\
            t_{2,2} = -t_{3,3}, &t_{2,3} = - t_{3,2}, & t_{2,4} = - t_{4,6}, & t_{2,6} = -t_{3,7},
            \\
            t_{2,7} = -t_{3,6}, & t_{2,8} = -t_{4,2}, & t_{2,9} = -t_{4,15} &t_{2,10} = -t_{3,11}, 
            \\
            t_{2,11} = -t_{3,10}, & t_{2,13} = - t_{4,11}, & t_{2,14} = -t_{3,15}, & t_{2,15} = - t_{3,14},
            \\
            t_{3,4} = -t_{5,6}, & t_{3,8} = -t_{5,2},& t_{3,9} = -t_{5,15},&t_{3,13} = -t_{5,11},
            \\
            t_{4,3} = - t_{5,8}, &t_{4,4} = - t_{5,7}, & t_{4,7} = - t_{5,4},& t_{4,8} = -t_{5,3},
            \\
            t_{4,9} = -t_{5,14}, & t_{4,10} = -t_{5,13}, & t_{4,13} = - t_{5,10}, & t_{4,14} = - t_{5,9}.
        \end{array}
    $$
\end{proof}

\section{Standard embedding for associative triple trisystems}

The concept of a di-endomorphism of a Jordan dialgebra was introduced in \cite{GK}. However, their definition is formulated in the framework of pseudo-algebras over a Hopf algebra and conformal algebras, making the approach somewhat intricate for those unfamiliar with these concepts. Nevertheless, the role of di-endomorphisms in the construction of the TKK extension of a Jordan dialgebra suggests that they are the appropriate tool for studying the standard embeddings of various di-structures. Therefore, in this work, we provide a more elementary approach to defining di-endomorphisms for any $\phi$-module.

For instance, any associative triple system of the first kind $A$ is embedded into an associative algebra $\mathfrak{M}(A,A) \oplus A$ where $\mathfrak{M}(A,A)$ is the module generated by the elements $(L(x,y), \ R(y,x))$ of $\mathrm{End}(A)\oplus \mathrm{End}(A)^\mathrm{op}$ (here $L(x,y)(z) = \{x,y,z\}$ and $R(y,x)(z) = \{z,y,x\}$ with $x,y,z \in A$) (see \cite{meyberg}). We aim for an analogous construction for any associative triple trisystem. To do so, some sort of di-endomorphisms will play an analogous role as the endomorphisms of $\mathfrak{M}(A,A)$.

\subsection{KP algorithm on endomorphisms of an algebra}

Let $D$ be an algebra. We define on $\mathrm{End}(D) \oplus D$ the bilinear product
$$
    f \cdot g := f \circ g, \quad f \cdot x := f(x) \quad \text{ and } \quad x \cdot y := xy,
$$
for all $f, g \in \mathrm{End}(D)$ and $x,y \in D$. We are interested in the non-trivial relation $(f \cdot g) \cdot x = f \cdot (g \cdot x)$.

Applying the KP algorithm, this new product gives rise two new products, let us denoted by $\dashv$ and $\vdash$, such that the new non-trivial relations of the  associative dialgebra are, for all $f,g \in \mathrm{End}(D)$ and $x \in D$:
\begin{align}
    \label{cor16}( f \dashv g) \dashv x = f \dashv (g \dashv x),
    \\
    \label{cor17}(f \vdash g) \dashv x = f \vdash (g \dashv x),
    \\
    \label{cor18}(f \vdash g) \vdash x = f \vdash (g \vdash x),
    \\
    \label{cor19}f \dashv (g \dashv x) = f \dashv (g \vdash x),
    \\
    \label{cor20}(f \dashv g) \vdash x = (f \vdash g) \vdash x.
\end{align}
Observe that, from any $f \in \mathrm{End}(D)$ we can define two new linear maps $f_1(x) := f \dashv x$ and $f_2(x) := f \vdash x$ for all $x \in D$. We can think that each map is a different evaluation of the same $f$. Therefore, using this approach we can rewrite the previous identities as
$$
    \begin{array}{ll}
        ( f \dashv g)_1(x) = f_1(g_1(x)) = f_1(g_2(x)), &
        (f \vdash g)_1(x) = f_2(g_1(x)),
        \\
        (f \vdash g)_2(x) = f _2(g_2(x)), &
        (f \dashv g)_2(x) = (f \vdash g)_2(x),
    \end{array}
$$
for all $f,g \in \mathrm{End}(D)$ and $x \in D$.

On another hand, the Part 2 of the KP algorithm give us also the following non-trivial identities between elements of $D$ and $\mathrm{End}(D)$
$$
    \begin{array}{lll}
        f \dashv (x \dashv y) = f \dashv (x \vdash y), & (f \dashv x) \vdash y = (f \vdash x) \vdash y, &x \dashv (f \dashv y) = x \dashv (f \vdash y).
    \end{array}
$$
or equivalently,
$$
    f_1(x \dashv y) = f_1 (x \vdash y), \quad f_1(x) \vdash y = f_2(x) \vdash y \quad \text{and} \quad x \dashv f_1( y) = x\dashv f_2(y),
$$
for all $f \in \mathrm{End}(D)$ and $x,y \in D$.

\subsection{Di-endomorphisms of a module}

Based on the ideas presented in the previous section, we give the following definitions:

\begin{definition} 
    Let $D,D'$ be $\phi$-modules. On $\mathrm{Hom}(D,D')\oplus \mathrm{Hom}(D,D')$ we define the following two bilinear maps 
    $$
    \prec, \succ \colon \mathrm{Hom}(D,D')\oplus \mathrm{Hom}(D,D')\times D \rightarrow D',
    $$
    given by
    $$
        \begin{array}{lcl}
            f \prec x := f_1(x) & \text{ and }&f \succ x := f_2(x) 
        \end{array}
    $$
    for any $f = (f_1 , f_2) \in \mathrm{Hom}(D,D')\oplus \mathrm{Hom}(D,D')$ and $x \in D$. The $\phi$-module equipped with these  two bilinear maps, $(\mathrm{Hom}(D,D')\oplus \mathrm{Hom}(D,D'), \prec, \succ)$, will be called the \textbf{di-homomorphisms from $D$ to $D'$} and denoted by $\mathrm{DiHom}(D,D')$.
    
    If $D = D'$, then we denote it by $\mathrm{DiEnd(D)}$ and its elements are referred to as \textbf{di-endomorphisms of $D$}.
\end{definition}
We can define two operations on $\mathrm{DiEnd}(D)$, inspired on the identities obtained after applying the KP algorithm, that gives structure of associative dialgebra, as we prove in the next proposition.
\begin{proposition}
    Let $D$ be a $\phi$-module. Then $\mathrm{DiEnd}(D)$ with the following two bilinear products
    $$
    \begin{array}{lcl}
        f \dashv g := (f_1 \circ g_2 , f_2 \circ g_2) & \text{ and } & f \vdash g := (f_2 \circ g_1 , f_2 \circ g_2),
    \end{array}
    $$
    is an associative dialgebra.
\end{proposition}
\begin{proof}
    Let us consider $f,g,h \in \mathrm{DiEnd}(D)$. Then
    \begin{enumerate}
        \item $(f \dashv g) \dashv h = (f_1 \circ g_2, f_2 \circ g_2) \dashv h = (f_1 \circ g_2 \circ h_2, f_2 \circ g_2 \circ h_2) = f \dashv (g_1 \circ h_2 , g_2 \circ h_2) =f \dashv (g \dashv h)$
        \item $(f \vdash g ) \vdash h = (f_2 \circ g_1, f_2 \circ g_2) \vdash h = (f_2 \circ g_2 \circ h_1, f_2 \circ g_2 \circ h_2)= f \vdash (g_2 \circ h_1, g_2 \circ h_2)=f \vdash (g \vdash h)$
        \item $f \dashv (g \dashv h) = f \dashv (g_1 \circ h_2, g_2 \circ h_2) = (f_1 \circ g_2 \circ h_2, f_2 \circ g_2 \circ h_2) = f \dashv (g_2 \circ h_1 , g_2 \circ h_2)=f \dashv (g \vdash h)$
        \item $(f \vdash g) \dashv h =(f_2 \circ g_1, f_2 \circ g_2) \dashv h = (f_2 \circ g_1 \circ h_2, f_2 \circ g_2 \circ h_2) = f \vdash (g_1 \circ h_2, g_2 \circ h_2) =f \vdash ( g \dashv h)$
        \item $(f \vdash g) \vdash h = (f_2 \circ g_1, f_2 \circ g_2) \vdash h = (f_2 \circ g_2 \circ h_1, f_2 \circ g_2 \circ h_2) = (f_1 \circ g_2, f_2 \circ g_2) \vdash h= (f \dashv g) \vdash h$
    \end{enumerate}
\end{proof}

Additionally, on $\mathrm{DiEnd}(D)$ there exist some relations between the products $\dashv$ and $\vdash$ and the bilinear maps $\prec$ and $\succ$, in the same way as the previous section suggested:
\begin{lemma}\label{lemma:diendomorphism}
    Let $D$ be a $\phi$-module. Then
    \begin{enumerate}
        \item $(f \vdash g) \succ x = (f \dashv g )\succ x = f \succ (g \succ x)$,
        \item $(f \vdash g) \prec x = f \succ (g \prec x)$,
        \item $(f \dashv g) \prec x  = f \prec (g \succ x)$,
    \end{enumerate}
    for all $f,g \in \mathrm{DiEnd}(D)$ and $x \in D$. 
\end{lemma}
\begin{proof}
    It is straightforward from the definition.
    \,
    
    1.
    \begin{align*}
        (f \vdash g) \succ x =(f_2 \circ g_1, f_2 \circ g_2) \succ x = f_2(g_2(x)) &= (f_1 \circ g_2, f_2 \circ g_2) \succ x = (f \dashv g) \succ x
        \\
        &= f\succ g_2(x) = f \succ (g \succ x) 
    \end{align*}
    \
    
    2. $(f \vdash g) \prec x = (f_2 \circ g_1, f_2 \circ g_2) \prec x = f_2(g_1(x))=f\succ g_1(x) = f\succ (g \prec x)$
    and
    \
    
    3. $(f \dashv g) \prec x =(f_1 \circ g_2, f_2 \circ g_2) \prec x = f_1(g_2(x))= f\prec g_2(x) = f \prec (g \succ x)$.
\end{proof}

\begin{remark}
    Notice that in general we do not get the relation $f \prec (g \prec x) = f \prec (g \succ x)$, or equivalently, $f_1(g_1(x)) = f_1(g_2(x))$. Nevertheless, as we will show later in Lemma \ref{lemma:extraidentity}, we will work on a particular subset of di-endomorphisms where this identity holds.
\end{remark}

For any $f = (f_1, f_2) \in \mathrm{DiEnd}(D)$, we set $\bar{f} := (\bar{f}_1, \bar{f}_2) \in \mathrm{End}(D) \oplus \mathrm{End}(D)$ where $(x)\bar{f}_1 := f_1(x)$ and $(x) \bar{f}_2 := f_2(x)$. And we also define  $\prec_\mathrm{op}, \succ_\mathrm{op} \colon D \times \mathrm{End}(D) \oplus\mathrm{End}(D)  \rightarrow D$ as
$$
    x \prec_\mathrm{op} \bar{f} := (x)\bar{f}_2, \quad \text{ and } \quad x \succ_\mathrm{op} \bar{f} :=  (x)\bar{f}_1.
$$
We will call it the \textbf{opposite di-endomorphisms of $D$} and we denote $(\mathrm{End}(D) \oplus\mathrm{End}(D) , \prec_\mathrm{op}, \succ_\mathrm{op})$ by $\mathrm{DiEnd}(D)^\mathrm{op}$. Moreover, for any $\bar{f}, \bar{g} \in \mathrm{DiEnd}(D)^\mathrm{op}$, we define the products 
$$
    \bar{f} \dashv_{\mathrm{op}} \bar{g} := \overline{g \vdash f}, \quad \text{ and } \quad \bar{f} \vdash_\mathrm{op} \bar{g} := \overline{g \dashv f}.
$$
It is straightforward to check that for any $x \in D$ and $f, g \in \mathrm{DiEnd}(D)$ we have
$$
    \begin{array}{lcl}
        x \prec_{\mathrm{op}}(\bar{f} \dashv_\mathrm{op} \bar{g}) = (g \vdash f) \succ x,& & x \succ_{\mathrm{op}}(\bar{f} \dashv_\mathrm{op} \bar{g}) = (g \vdash f) \prec x,
        \\
        x \prec_{\mathrm{op}}(\bar{f} \vdash_\mathrm{op} \bar{g}) = (g \dashv f) \succ x,& & x \succ_{\mathrm{op}}(\bar{f} \vdash_\mathrm{op} \bar{g}) = (g \dashv f) \prec x.
    \end{array}
$$

\subsection{The standard embedding of a triple trisystem of the first kind}

We say that a $\phi$-module $T$ is a \textbf{triple trisystem} if it is equipped with three trilinear maps $\{ \, , \, ,\,\}_i \colon T \times T \times T \rightarrow T$, $i \in \{1,2,3\}$. 

In the following, we define a collection of elements that play a fundamental role in establishing the standard embedding of an associative triple trisystem. To motivate these definitions and notations, let us first recall the relevant background. In Proposition \ref{th:asstoass}, we demonstrated that any associative dialgebra $(D, \dashv, \vdash)$ gives rise to an associative triple trisystem of the first kind by defining the following triple products (\ref{cor13}), (\ref{cor14}) and (\ref{cor15}):
\begin{align*}
    \{x,y,z\}_1 &= x \dashv (y \dashv z) = (x \dashv y) \dashv z = x \dashv (y \vdash z),
    \\
    \{x,y,z\}_2 &= x \vdash (y \dashv z) = (x \vdash y) \dashv z 
    \\
    \{x,y,z\}_3 &= x \vdash ( y \vdash z) = (x \vdash y) \vdash z = (x \dashv y) \vdash z.
\end{align*}
We aim to encapsulate the behavior of these parentheses. For example, considering the following di-endomorphism $x \triangleleft_l y = ( \{x,y, \ \cdot \ \}_1, \ \{ x,y, \ \cdot \ \}_3)$, satisfying
$$
    (x \triangleleft_l y) \prec z = \{x,y,z\}_1 = (x \dashv y) \dashv z,
$$
and
$$
    (x \triangleleft_l y) \succ z = \{x,y,z\}_3 = (x \dashv y) \vdash z.
$$
However, the parentheses can also act on the right. Therefore, we will define the element in $\mathrm{DiEnd}(D)\oplus \mathrm{DiEnd}(D)^\mathrm{op}$
$$
    x \triangleleft y = (\{x,y, \ \cdot \ \}_1, \ \{x,y, \ \cdot \ \}_3) \oplus (\{\ \cdot \ ,x,y\}_2, \ \{ \ \cdot \ , x,y\}_1)
$$
satisfying additionally
$$
    z \prec(x \triangleleft y) = \{z,x,y\}_1 = z  \dashv (x \dashv y)
$$
and
$$
    z \succ (x \triangleleft y) = \{z,x,y\}_2 = z \vdash (x \dashv y).
$$
Thus, in this case, for example, we will denote the endomorphisms appearing in $x \triangleleft y$ as
$$
    x \triangleleft y = (L_1(x,y), \ L_3(x,y)) \oplus (\overline{R_2(x,y)}, \ \overline{R_1(x,y)} ).
$$

Let $T$ be a triple trisystem. We define, for any $x,y \in T$, the following six linear maps $L_i(x,y),\, R_i(x,y) \colon T \rightarrow T$ with $i \in \{1,2,3\}$
$$
    L_i(x,y)(z) := \{x,y,z\}_i \quad \text{ and } \quad  R_i(x,y)(z) := \{z,x,y\}_i,
$$
for all $z \in T$. Let us define the following di-endomorphisms
$$
    \begin{array}{ll}
        L^\triangleleft(x,y) = (L_1(x,y), \, L_3(x,y)), & R^\triangleleft(x,y) = (R_2(x,y), \, R_1(x,y)) ,
        \\
        L^\triangleright(x,y) = (L_2(x,y), \, L_3(x,y)), & R^\triangleright(x,y) = (R_3(x,y), \, R_1(x,y)).
    \end{array}
$$
We consider the $\phi$-modules $\mathfrak{L}(T,T), \, \mathfrak{R}(T,T) \subset \mathrm{DiEnd}(T)$ generated by $L^\mu (x,y)$ and $R^\mu(x,y)$, for all $x,y \in T$ and $\mu \in \{ \triangleleft, \triangleright\}$, respectively.

\begin{lemma} \label{lemma:extraidentity}
    Let $A$ be an associative triple trisystem of any kind. Then 
    $$
        f \prec (g \prec x) = f \prec (g \succ x),
    $$
    for all $f,g \in \mathfrak{L}(A,A) \cup \mathfrak{R}(A,A)$ and $x \in A$.
\end{lemma}
\begin{proof}
    It is straightforward. If $A$ is of the first kind is following from the relations (\ref{cor7}), (\ref{cor8}), (\ref{cor11}) and (\ref{cor12}), and if $A$ is of the second kind from (\ref{cor7b}), (\ref{cor8b}), (\ref{cor11b}) and (\ref{cor12b}).
\end{proof}
In particular, if $A$ is an associative triple trisystem of any kind, the elements of $\mathfrak{L}(A,A) \oplus A$ or $\mathfrak{R}(A,A) \oplus A$ satisfy the five identities (\ref{cor16})-(\ref{cor20}) by identifying the actions $\prec$ and $\succ$ with $\dashv$ and $\vdash$, respectively.

Let $\mathfrak{M}(T,T)$ be the $\phi$-module generated by all elements of the form
$$
    x \,  \triangleleft \, y := L^\triangleleft(x,y) \oplus  \overline{R^\triangleleft(x,y)} \in \mathrm{DiEnd}(T)\oplus \mathrm{DiEnd}(T)^{\mathrm{op}},
$$
and 
$$
    x \, \triangleright \,  y:= L^\triangleright(x,y) \oplus \overline{R^\triangleright(x,y)} \in \mathrm{DiEnd}(T)\oplus \mathrm{DiEnd}(T)^{\mathrm{op}},
$$
with $x, y\in T$.

We consider the left and right actions on $\mathfrak{M}(T,T)$ induced by $\mathrm{DiEnd}(T)$ and $\mathrm{DiEnd}(T)^\mathrm{op}$. Indeed, two bilinear maps acting on the left $\prec , \succ \colon T \times \mathfrak{M} (T,T) \rightarrow T$ are given by 
$$
    z \prec  (x \, \mu \, y) := z \prec_\mathrm{op}\overline{R^\mu(x,y)}, 
$$
and
$$z \succ (x \, \mu \, y) := z \succ_\mathrm{op}\overline{R^\mu (x,y)},
$$
for $\mu \in \{ \triangleleft, \triangleright\}$. 
Similarly,  two bilinear maps acting on the right $\prec, \succ \colon \mathfrak{M} (T,T) \times T \rightarrow T$ are defined as
$$
    (x \, \mu \, y) \prec z := L^\mu(x,y) \prec z,
$$
and
$$
    (x \, \mu \, y) \succ z := L^\mu(x,y) \succ z,
$$
for all $x,y,z \in T$ and $\mu \in \{ \triangleleft , \triangleright \}$.

Finally, in this module, we also consider the products $\dashv$ and $\vdash$ that act component-wise on $\mathfrak{M}(T,T)$, induced by those defined in $\mathrm{DiEnd}$ and $\mathrm{DiEnd}^\mathrm{op}$. Concretely, for $\mu, \nu \in \{ \triangleleft , \, \triangleright \}$, we define
$$
    (x \, \mu \, y) \dashv (z \, \nu \, u) := (L^\mu(x,y) \dashv L^\nu(z,u))\oplus (\overline{R^\mu(x,y)} \dashv_\mathrm{op} \overline{R^\nu(z,u))},
$$
and
$$
    (x \, \mu \, y) \vdash (z \, \nu \, u) := (L^\mu(x,y) \vdash L^\nu(z,u)) \oplus( \overline{R^\mu(x,y)} \vdash_\mathrm{op} \overline{R^\nu(z,u)}).
$$

Nevertheless, throughout the remainder of this manuscript, in order to simplify notation and readability, we will omit the overline and the  subscript $\mathrm{op}$ on the opposite di-endomorphisms.

Some relations of these products are given in the following lemma. 
\begin{lemma}
    Let $T$ be a triple trisystem. Then 
    \begin{multicols}{2}
    \begin{enumerate}
        \item $(x \triangleleft y) \dashv (z \triangleleft u) = (x \triangleleft y) \dashv (z \triangleright u)$,
        \item $(x \triangleright y) \dashv (z \triangleleft u) = (x \triangleright y) \dashv (z \triangleright u)$,
        \item $(x \triangleleft y) \vdash (z \triangleleft u) = (x \triangleright y) \vdash (z \triangleleft u)$,
        \item $(x \triangleleft y) \vdash (z \triangleright u) = (x \triangleright y) \vdash (z \triangleright u)$.
    \end{enumerate}
    \end{multicols}
    for all $x,y,z,u \in T$.
\end{lemma}
\begin{proof}
    The calculation is straightforward. For convenience and since these results will be used later in Proposition \ref{lemma:products}, we explicitly provide the expressions for the left-hand side of the equality:
    \begin{enumerate} 
        \item $(x \triangleleft y) \dashv (z \triangleleft u) = (L^\triangleleft(x,y) \dashv L^\triangleleft(z,u)) \oplus (R^\triangleleft(x,y) \dashv R^\triangleleft(z,u))$:
        
        On the left side 
        $$
            L^\triangleleft(x,y) \dashv L^\triangleleft(z,u) = (L_1(x,y) \circ L_3(z,u), \, L_3(x,y) \circ L_3(z,u)),
        $$
        so
        \begin{align*}
            (L^\triangleleft(x,y) \dashv L^\triangleleft(z,u)) \prec v &= (L_1(x,y) \circ L_3(z,u))(v)  = \{x,y,\{z,u,v\}_3\}_1
            \\
            (L^\triangleleft(x,y) \dashv L^\triangleleft(z,u)) \succ v &= (L_3(x,y) \circ L_3(z,u))(v)  = \{x,y,\{z,u,v\}_3\}_3.
        \end{align*}
        On the right side 
        $$
            R^\triangleleft(x,y) \dashv R^\triangleleft(z,u) = ( R_2(x,y) \circ R_1(z,u),\,  R_1(x,y)\circ R_1(z,u)),
        $$
        and therefore
        \begin{align*}
            v \prec (R^\triangleleft(x,y) \dashv R^\triangleleft(z,u)) &= (v)(R_1(x,y)\circ R_1(z,u))   = \{\{v,x,y\}_1,z,u\}_1
            \\
            v \succ (R^\triangleleft(x,y) \dashv R^\triangleleft(z,u)) &= (v)( R_2(x,y) \circ R_1(z,u))  = \{\{v,x,y\}_2,z,u\}_1
        \end{align*}

        \item $(x \triangleright y) \dashv (z \triangleleft u) = (L^\triangleright(x,y) \dashv L^\triangleleft (z,u)) \oplus (R^\triangleright(x,y) \dashv R^\triangleleft (z,u))$:

        On the left side,
        $$
            L^\triangleright(x,y) \dashv L^\triangleleft (z,u) = (L_2(x,y) \circ L_3(z,u), \ L_3(x,y) \circ L_3(z,u)).
        $$
        Therefore
        \begin{align*}
            (L^\triangleright(x,y) \dashv L^\triangleleft (z,u))\prec v &= (L_2(x,y) \circ L_3(z,u))(v)  = \{x,y,\{z,u,v\}_3\}_2
            \\
            (L^\triangleright(x,y) \dashv L^\triangleleft (z,u))\succ v &= (L_3(x,y) \circ L_3(z,u))(v)  = \{x,y,\{z,u,v\}_3\}_3
        \end{align*}

        On the right side,
        $$
             R^\triangleright(x,y) \dashv  R^\triangleleft (z,u) = (R_3(x,y) \circ R_1(z,u),\  R_1(x,y) \circ R_1(z,u)).
        $$
        Therefore
        \begin{align*}
            v \prec (R^\triangleright(x,y) \dashv  R^\triangleleft (z,u)) &= (v)(R_1(x,y) \circ R_1(z,u))   = \{\{v,x,y\}_1,z,u\}_1
            \\
            v \succ (R^\triangleright(x,y) \dashv  R^\triangleleft (z,u)) &= (v)(R_3(x,y) \circ R_1(z,u))   = \{\{v,x,y\}_3,z,u\}_1
        \end{align*}
        
        \item $(z \triangleleft u) \vdash (x \triangleleft y) = (L^\triangleleft(z,u) \vdash L^\triangleleft(x,y)) \oplus ( R^\triangleleft(z,u) \vdash R^\triangleleft(x,y))$:

        On the left side 
        $$
            L^\triangleleft(z,u) \vdash L^\triangleleft(x,y) = (L_3(z,u) \circ L_1(x,y), \ L_3(z,u) \circ L_3(x,y)).
        $$
        Therefore
        \begin{align*}
            (L^\triangleleft(z,u) \vdash L^\triangleleft(x,y)) \prec v &= (L_3(z,u) \circ L_1(x,y))(v)  = \{z,u,\{x,y,v\}_1\}_3
            \\
            (L^\triangleleft(z,u) \vdash L^\triangleleft(x,y)) \succ v &= (L_3(z,u) \circ L_3(x,y))(v)   = \{z,u,\{x,y,v\}_3\}_3
        \end{align*}
        
        On the right side, 
        $$
            R^\triangleleft(z,u) \vdash  R^\triangleleft(x,y) = ( R_1(z,u) \circ R_2(x,y) ,\  R_1(z,u) \circ R_1(x,y)).
        $$
        Therefore
        \begin{align*}
            v \prec (R^\triangleleft(z,u) \vdash  R^\triangleleft(x,y)) &= (v)(R_1(z,u) \circ R_1(x,y))   = \{\{v,z,u\}_1,x,y\}_1
            \\
            v \succ (R^\triangleleft(z,u) \vdash  R^\triangleleft(x,y)) &= (v)(R_1(z,u) \circ R_2(x,y) )  = \{\{v,z,u\}_1,x,y\}_2
        \end{align*}

        \item $(z \triangleleft u) \vdash (x \triangleright y) = (L^\triangleleft(z,u) \vdash L^\triangleright(x,y)) \oplus ( R^\triangleleft(z,u) \vdash  R^\triangleright(x,y) )$:

        On the left side, 
        $$
            L^\triangleleft(z,u) \vdash L^\triangleright(x,y) = (L_3(z,u) \circ L_2(x,y),\ L_3(z,u)  \circ L_3(x,y)).
        $$
        Therefore
        \begin{align*}
            (L^\triangleleft(z,u) \vdash L^\triangleright(x,y))\prec v &=(L_3(z,u) \circ L_2(x,y))(v)  = \{z,u,\{x,y,v\}_2\}_3
            \\
            (L^\triangleleft(z,u) \vdash L^\triangleright(x,y))\succ v &= (L_3(z,u)  \circ L_3(x,y))(v)   = \{z,u,\{x,y,v\}_3\}_3
        \end{align*}

        On the right side,
        $$
            R^\triangleleft(z,u) \vdash  R^\triangleright(x,y) = ( R_1(z,u)\circ R_3(x,y), \ R_1(z,u)\circ R_1(x,y)).
        $$
        Therefore,
        \begin{align*}
            v \prec (R^\triangleleft(z,u) \vdash  R^\triangleright(x,y)) &= (v)(R_1(z,u)\circ R_1(x,y))   = \{\{v,z,u\}_1,x,y\}_1
            \\
            v \succ (R^\triangleleft(z,u) \vdash  R^\triangleright(x,y)) &= (v)( R_1(z,u)\circ R_3(x,y))   = \{\{v,z,u\}_1,x,y\}_3
        \end{align*}
    \end{enumerate}
\end{proof}

If $A$ is an associative triple trisystem of the first kind, then in the following proposition, we prove that these products are closed in $\mathfrak{M}(A,A)$, endowing it with the structure of an associative dialgebra.

\begin{proposition}\label{lemma:products}
    Let $A$ be an associative triple trisystem of the first kind. Then the products in $\mathfrak{M}(A,A)$ satisfy
    $$
        \begin{array}{ll}
            1. \quad (x \triangleleft y)\dashv (z \triangleleft u ) = x \triangleleft \{y,z,u\}_1 = \{x,y,z\}_1 \triangleleft u = (x \triangleleft y)\dashv (z \triangleright u ),  
            \\
            2. \quad (x \triangleleft y)\vdash (z \triangleleft u ) =\{x,y,z\}_3 \triangleleft u = x \triangleright \{y,z,u\}_2 = (x \triangleright y)\vdash (z \triangleleft u ) ,
            \\
            3. \quad (x \triangleright y)\dashv (z \triangleleft u ) = x \triangleright\{y,z,u\}_1 = \{x,y,z\}_2 \triangleleft u = (x \triangleright y)\dashv (z \triangleright u ),  
            \\
            4. \quad (x \triangleleft y)\vdash (z \triangleright u ) =\{x,y,z\}_3 \triangleright u = x \triangleright \{y,z,u\}_3 = (x \triangleright y)\vdash (z \triangleright u ),
        \end{array}
    $$
    for all $x,y,z,u \in A$. In particular, $\mathfrak{M}(A,A)$ is an associative dialgebra.
\end{proposition}
\begin{proof}
    These relations are straightforward. However, unlike the previous lemma, the structure of an associative triple trisystem is required here:
    \begin{enumerate}
        \item Let us prove that $(x \triangleleft y) \dashv (z \triangleleft u) = x \triangleleft \{y,z,u\}_1 = \{x,y,z\}_1 \triangleleft u$:
        \begin{align*}
            ((x \triangleleft y) \dashv (z \triangleleft u)) \prec v = \{x,y,\{z,u,v\}_3 \}_1  &= \{x, \{y,z,u\}_1,v\}_1 = (x \triangleleft \{y,z,u\}_1) \prec v
            \\
            &=\{\{x,y,z\}_1,u,v\}_1 = (\{x,y,z\}_1 \triangleleft u) \prec v
        \end{align*}
        \begin{align*}
            ((x \triangleleft y) \dashv (z \triangleleft u)) \succ v = \{x,y,\{z,u,v\}_3\}_3 &= \{x,\{y,z,u\}_1,v\}_3 = (x \triangleleft \{y,z,u\}_1) \succ v
            \\
            &= \{\{x,y,z\}_1,u,v\}_3 = (\{x,y,z\}_1 \triangleleft u) \succ v
        \end{align*}
        \begin{align*}
            v \prec ((x \triangleleft y) \dashv (z \triangleleft u)) = \{\{v,x,y\}_1,z,u\}_1 &= \{v,x,\{y,z,u\}_1\}_1 = v \prec (x \triangleleft \{y,z.u\}_1)
            \\
            &= \{v,\{x,y,z\}_1,u\}_1 = v \prec (\{x,y,z\}_1 \triangleleft u)
        \end{align*}
        \begin{align*}
            v \succ ((x \triangleleft y) \dashv (z \triangleleft u)) = \{\{v,x,y\}_2,z,u\}_1 &= \{v,x,\{y,z,u\}_1\}_2 = v \succ (x \triangleleft \{y,z.u\}_1)
            \\
            &=\{v,\{x,y,z\}_1,u\}_2 = v \succ (\{x,y,z\}_1 \triangleleft u)
        \end{align*}

        \item Let us prove $(x \triangleleft y) \vdash (z \triangleleft u) = \{x,y,z\}_3 \triangleleft u = x \triangleright \{y,z,u\}_2$:
        \begin{align*}
            ((x \triangleleft y) \vdash (z \triangleleft u)) \prec v = \{x,y,\{z,u,v\}_1\}_3 &= \{\{x,y,z\}_3,u,v\}_1 = (\{x,y,z\}_3 \triangleleft u) \prec v
            \\
            &=\{x,\{y,z,u\}_2,v\}_2= (x \triangleright \{y,z,u\}_2)\prec v
        \end{align*}
        \begin{align*}
            ((x \triangleleft y) \vdash (z \triangleleft u)) \succ v = \{x,y,\{z,u,v\}_3\}_3 &= \{\{x,y,z\}_3,u,v\}_3 = (\{x,y,z\}_3 \triangleleft u) \succ v
            \\
            &=\{x,\{y,z,u\}_2,v\}_3 = (x \triangleright \{y,z,u\}_2)\succ v
        \end{align*}
        \begin{align*}
            v \prec ((x \triangleleft y) \vdash (z \triangleleft u)) = \{\{v,x,y\}_1,z,u\}_1& = \{v,\{x,y,z\}_3,u\}_1 = v \prec (\{x,y,z\}_3 \triangleleft u) 
            \\
            &= \{v,x,\{y,z,u\}_2\}_1 = v \prec (x \triangleright \{y,z,u\}_2)
        \end{align*}
        \begin{align*}
            v \succ ((x \triangleleft y) \vdash (z \triangleleft u)) = \{\{v,x,y\}_1,z,u\}_2  &= \{v,\{x,y,z\}_3, u\}_2 = v \succ (\{x,y,z\}_3 \triangleleft u) 
            \\
            &=\{v,x,\{y,z,u\}_2\}_3 = v \succ (x \triangleright \{y,z,u\}_2)
        \end{align*}

        \item Let us prove $(x \triangleright y) \dashv (z \triangleleft u) = x \triangleright \{y,z,u\}_1 = \{x,y,z\}_2 \triangleleft u$:
        \begin{align*}
            ((x \triangleright y) \dashv (z \triangleleft u)) \prec v =\{x, y, \{z, u, v\}_3\}_2 &=\{x,\{y,z,u\}_1,v\}_2 = (x \triangleleft \{y,z,u\}_1)\prec v
            \\
            &=\{\{x,y,z\}_2,u,v\}_1 = (\{x,y,z\}_2 \triangleleft u) \prec v
        \end{align*}
        \begin{align*}
            ((x \triangleright y) \dashv (z \triangleleft u)) \succ v = \{x, y, \{z, u, v\}_3\}_3&=\{x,\{y,z,u\}_1,v\}_3 = (x \triangleleft \{y,z,u\}_1)\succ v
            \\
            &= \{\{x,y,z\}_2,u,v\}_3 = ( \{x,y,z\}_2 \triangleleft u) \succ v
        \end{align*}
        \begin{align*}
            v \prec ((x \triangleright y) \dashv (z \triangleleft u)) = \{\{v,x,y\}_1,z,u\}_1&= \{v,x,\{y,z,u\}_1\}_1 = v \prec ( x \triangleright \{y,z,u\}_1)
            \\
            &= \{v,\{x,y,z\}_2,u\}_1 = v \prec (\{x,y,z\}_2\triangleleft u)
        \end{align*}
        \begin{align*}
            v \succ ((x \triangleright y) \dashv (z \triangleleft u)) = \{\{v,x,y\}_3,z,u\}_1 &=\{v,x,\{y,z,u\}_1\}_3 = v \succ ( x \triangleright \{y,z,u\}_1)
            \\
            &=\{v,\{x,y,z\}_2,u\}_2 = v \succ (\{x,y,z\}_2 \triangleleft u)
        \end{align*}
        
        \item Let us prove $(x \triangleleft y) \vdash (z \triangleright u) = \{x,y,z\}_3 \triangleright u = x \triangleright \{y,z,u\}_3$:
        \begin{align*}
            ((x \triangleleft y) \vdash (z \triangleright u)) \prec v = \{x,y,\{z,u,v\}_2\}_3&=\{\{x,y,z\}_3,u,v\}_2 =(\{x,y,z\}_3 \triangleright u)\prec v
            \\
            &= \{x,\{y,z,u\}_3,v\}_2 = (x \triangleright \{y,z,u\}_3) \prec v
        \end{align*}
        \begin{align*}
            ((x \triangleleft y) \vdash (z \triangleright u)) \succ v = \{x,y,\{z,u,v\}_3\}_3&= \{\{x,y,z\}_3,u,v\}_3 = (\{x,y,z\}_3 \triangleright u) \succ v
            \\
            &= \{x,\{y,z,u\}_3,v\}_3 = (x \triangleright \{y,z,u\}_3) \succ v
        \end{align*}
        \begin{align*}
            v \prec ((x \triangleleft y) \vdash (z \triangleright u))  = \{\{v,x,y\}_1,z,u\}_1&=\{v,\{x,y,z\}_3,u\}_1= v\prec (\{x,y,z\}_3 \triangleright u)
            \\
            &= \{v,x,\{y,z,u\}_3\}_1 = v \prec (x \triangleright \{y,z,u\}_3)
        \end{align*}
        \begin{align*}
            v \succ ((x \triangleleft y) \vdash (z \triangleright u))  = \{\{v,x,y\}_1,z,u\}_3&= \{v,\{x,y,z\}_3,u\}_3 = v \succ (\{x,y,z\}_3 \triangleright u)
            \\
            &= \{v,x,\{y,z,u\}_3\}_3 = v \succ (x \triangleright \{y,z,u\}_3)
        \end{align*}
    \end{enumerate}
\end{proof}

\begin{definition}
    Let $A$ be an associative triple trisystem of the first kind. We define the \textbf{standard embedding of $A$} as the dialgebra 
    $$
        \mathrm{U}(A) := \mathfrak{M}(A,A) \oplus A
    $$
    equipped with two bilinear products, $\dashv$ and $\vdash$, defined as follows:
    $$
        \begin{array}{ll}
            x \dashv y := x \triangleleft y,&  x \vdash y := x \triangleright y
            \\
            (x \triangleleft y) \dashv z := (x \triangleleft y)  \prec z =  \{x,y,z\}_1, & 
            (x \triangleleft y) \vdash z := (x \triangleleft y)  \succ z  =\{x,y,z\}_3
            \\
            z \dashv (x \triangleleft y) := z \prec (x \triangleleft y)   = \{z,x,y\}_1, & 
            z \vdash (x \triangleleft y):= z \succ (x \triangleleft y)    = \{z,x,y\}_2
            \\
            (x \triangleright y) \dashv z := (x \triangleright y) \prec z =   \{x,y,z\}_2, & 
            (x \triangleright y) \vdash z := (x \triangleright y) \succ z   = \{x,y,z\}_3
            \\
            z \dashv (x \triangleright y) := z \prec (x \triangleright y)  = \{z,x,y\}_1, & 
            z \vdash (x \triangleright y):= z \succ (x \triangleright y)   =\{z,x,y\}_3
        \end{array}
    $$
    for all $x,y,z \in A$. In $\mathfrak{M}(A,A)$, the products $\dashv$ and $\vdash$  coincide with those previously defined.
\end{definition}

\begin{theorem}
    $(A, \{ \, , \,  ,  \, \}_i)$ is an associative triple trisystem of the first kind if and only if $(\mathrm{U}(A), \dashv, \vdash)$ is an associative dialgebra. Moreover, $\{x,y,z\}_1 =x \dashv (y \dashv z)$, $\{x,y,z\}_2 = x \vdash (y \dashv z)$ and $\{x,y,z\}_3 = x \vdash (y \vdash z)$, for all $x,y,z \in A$.
\end{theorem}
\begin{proof}
    Let us prove the five identities that satisfy an associative dialgebra:
    \begin{itemize}
        \item[(a)] $(\alpha \dashv \beta) \dashv \gamma = \alpha \dashv ( \beta \dashv \gamma)$,
        \item[(b)] $(\alpha \vdash \beta) \vdash \gamma = \alpha \vdash ( \beta \vdash \gamma)$,
        \item[(c)] $\alpha \dashv (\beta \dashv \gamma) = \alpha \dashv ( \beta \vdash \gamma)$,
        \item[(d)] $(\alpha \vdash \beta) \dashv \gamma = \alpha \vdash ( \beta \dashv \gamma)$,
        \item[(e)] $(\alpha \vdash \beta) \vdash \gamma = (\alpha \dashv \beta) \vdash \gamma$,
    \end{itemize}
    for any $\alpha, \beta, \gamma \in \mathrm{U}(A)$. We already know that $\mathfrak{M}(A,A)$ is an associative dialgebra, so we just have to prove the following cases:
    \begin{enumerate}
        \item The case $\alpha, \beta, \gamma \in A$ is clear by construction.

        \item  Any combination of two elements from $A$ and one from $\mathfrak{M}(A,A)$ hold by Proposition \ref{lemma:products}.

        \item Let us consider $\alpha = x \triangleleft y \in \mathfrak{M}(A,A), \beta = z \triangleleft u \in \mathfrak{M}(A,A)$ and $\gamma = v \in A$. Let us show that:
        \begin{enumerate}
            \item $((x \triangleleft y) \dashv (z \triangleleft u)) \dashv v = (x \triangleleft y) \dashv ((z \triangleleft u) \dashv v)$:
            \begin{align*}
                ((x \triangleleft y) \dashv (z \triangleleft u)) \dashv v &= (x \triangleleft \{y,z,u\}_1) \dashv v = \{x,\{y,z,u\}_1,v\}_1
                \\
                (x \triangleleft y) \dashv ((z \triangleleft u) \dashv v) &= (x \triangleleft y) \dashv \{z,u,v\}_1 =\{x,y,\{z,u,v\}_1\}_1
            \end{align*}

            \item $((x \triangleleft y) \vdash (z \triangleleft u)) \vdash v = (x \triangleleft y) \vdash ((z \triangleleft u) \vdash v)$:
            \begin{align*}
                ((x \triangleleft y) \vdash (z \triangleleft u)) \vdash v &= (\{x,y,z\}_3\triangleleft u) \vdash v = \{\{x,y,z\}_3,u,v\}_3
                \\
                (x \triangleleft y) \vdash ((z \triangleleft u) \vdash v) &= (x \triangleleft y) \vdash \{z,u,v\}_3 = \{x,y,\{z,u,v\}_3\}_3
            \end{align*}

            \item $(x \triangleleft y) \dashv ((z \triangleleft u) \dashv v) = (x \triangleleft y) \dashv ((z \triangleleft u) \vdash v)$:
            \begin{align*}
                (x \triangleleft y) \dashv ((z \triangleleft u) \dashv v) &= (x \triangleleft y) \dashv \{z,u,v\}_1 = \{x,y,\{z,u,v\}_1\}_1
                \\
                (x \triangleleft y) \dashv ((z \triangleleft u) \vdash v) &= (x \triangleleft y ) \dashv \{z,u,v\}_3 = \{x,y,\{z,u,v\}_3\}_1
            \end{align*}

            \item $((x \triangleleft y) \vdash (z \triangleleft u)) \dashv v = (x \triangleleft y) \vdash ((z \triangleleft u) \dashv v)$:
            \begin{align*}
                ((x \triangleleft y) \vdash (z \triangleleft u)) \dashv v &= (\{x,y,z\}_3\triangleleft u) \dashv v = \{\{x,y,z\}_3,u,v\}_1
                \\
                (x \triangleleft y) \vdash ((z \triangleleft u) \dashv v) &= (x \triangleleft y) \vdash \{z,u,v\}_1 = \{x,y,\{z,u,v\}_1\}_3
            \end{align*}

            \item $((x \triangleleft y) \vdash (z \triangleleft u)) \vdash v = ((x \triangleleft y) \dashv (z \triangleleft u)) \vdash v$:
            \begin{align*}
                ((x \triangleleft y) \vdash (z \triangleleft u)) \vdash v &= (\{x,y,z\}_3 \triangleleft u) \vdash v = \{\{x,y,z\}_3,u,v\}_3
                \\
                ((x \triangleleft y) \dashv (z \triangleleft u)) \vdash v &= (x \triangleleft \{y,z,u\}_1) \vdash v = \{x,\{y,z,u\}_1,v\}_3
            \end{align*}
        \end{enumerate}

        \item Let us consider $\alpha = x \triangleright y \in \mathfrak{M}(A,A), \beta = z \triangleleft u \in \mathfrak{M}(A,A)$ and $\gamma = v \in A$. Let us show that:
        \begin{enumerate}
            \item $((x\triangleright y) \dashv (z \triangleleft u)) \dashv v = (x \triangleright y) \dashv ( (z\triangleleft u) \dashv v)$:
            \begin{align*}
                ((x\triangleright y) \dashv (z \triangleleft u)) \dashv v &= (x \triangleright \{y,z,u\}_1) \dashv v = \{x,\{y,z,u\}_1,v\}_2
                \\
                (x \triangleright y) \dashv ( (z\triangleleft u) \dashv v) &= (x \triangleright y) \dashv \{z,u,v\}_1 = \{x,y,\{z,u,v\}_1\}_2
            \end{align*}

            \item $((x \triangleright y) \vdash (z \triangleleft u)) \vdash v = (x \triangleright y) \vdash ((z \triangleleft u) \vdash v)$:
            \begin{align*}
                ((x \triangleright y) \vdash (z \triangleleft u)) \vdash v &= (\{x,y,z\}_1\triangleleft u) \vdash  v = \{\{x,y,z\}_3,u,v\}_3
                \\
                (x \triangleright y) \vdash ((z \triangleleft u) \vdash v) &= (x \triangleright y) \vdash \{z,u,v\}_3 = \{x,y,\{z,u,v\}_3\}_3
            \end{align*}

            \item $(x \triangleright y) \dashv ((z \triangleleft u) \dashv v) = (x \triangleright y) \dashv ((z \triangleleft u) \vdash v)$: 
            \begin{align*}
                (x \triangleright y) \dashv ((z \triangleleft u) \dashv v) &= (x \triangleright y) \dashv \{z,u,v\}_1 = \{x,y,\{z,u,v\}_1\}_2
                \\
                (x \triangleright y) \dashv ((z \triangleleft u) \vdash v) &= (x \triangleright y) \dashv \{z,u,v\}_3 = \{x,y,\{z,u,v\}_3\}_2
            \end{align*}

            \item $((x \triangleright y) \vdash (z \triangleleft u)) \dashv v = (x \triangleright y) \vdash ((z \triangleleft u) \dashv v)$:
            \begin{align*}
                ((x \triangleright y) \vdash (z \triangleleft u)) \dashv v &= (\{x,y,z\}_3 \triangleleft u) \dashv v = \{\{x,y,z\}_3,u,v\}_1 
                \\
                (x \triangleright y) \vdash ((z \triangleleft u) \dashv v) &= (x \triangleright y ) \vdash \{z,u,v\}_1 = \{x,y,\{z,u,v\}_1\}_3
            \end{align*}

            \item $((x \triangleright y) \vdash (z \triangleleft u)) \vdash v = ((x \triangleright y) \dashv (z \triangleleft u)) \vdash v$:
            \begin{align*}
                ((x \triangleright y) \vdash (z \triangleleft u)) \vdash v &= (\{x,y,z\}_3 \triangleleft u) \vdash v = \{\{x,y,z\}_3, u,v\}_3
                \\
                ((x \triangleright y) \dashv (z \triangleleft u)) \vdash v &= (x \triangleright \{y,z,u\}_1) \vdash v = \{x,\{y,z,u\}_1,v\}_3
            \end{align*}
        \end{enumerate}

        \item Let us consider $\alpha = x \triangleleft y \in \mathfrak{M}(A,A), \beta = z \triangleright u \in \mathfrak{M}(A,A)$ and $\gamma = v \in A$. Let us show that:
        \begin{enumerate}
            \item $((x \triangleleft y) \dashv (z \triangleright u))\dashv v = ( x \triangleleft y) \dashv ( (z \triangleright u) \dashv v)$:
            \begin{align*}
                ((x \triangleleft y) \dashv (z \triangleright u))\dashv v &= (x\triangleleft \{y,z,u\}_1) \dashv v = \{x,\{y,z,u\}_1,v\}_1
                \\
                ( x \triangleleft y) \dashv ( (z \triangleright u) \dashv v) &= (x \triangleleft y) \dashv \{z,y,v\}_2 = \{x,y,\{z,y,v\}_2\}_1
            \end{align*}
            
            \item $((x \triangleleft y) \vdash (z \triangleright u)) \vdash v = (x \triangleleft y) \vdash ( (z \triangleright u) \vdash v)$:
            \begin{align*}
                ((x \triangleleft y) \vdash (z \triangleright u)) \vdash v & = (\{x,y,z\}_3 \triangleright u) \vdash v = \{\{x,y,z\}_3,u,v\}_3
                \\
                (x \triangleleft y) \vdash ( (z \triangleright u) \vdash v) &= (x \triangleleft y) \vdash \{z,u,v\}_3 = \{x,y,\{z,u,v\}_3\}_3
            \end{align*}
            
            \item $(x \triangleleft y) \dashv ((z \triangleright u) \dashv v) = (x \triangleleft y) \dashv ( (z \triangleright u) \vdash v)$:
            \begin{align*}
                (x \triangleleft y) \dashv ((z \triangleright u) \dashv v) &= (x \triangleleft y) \dashv \{z,u,v\}_2 = \{x,y,\{z,u,v\}_2\}_1
                \\
                (x \triangleleft y) \dashv ( (z \triangleright u) \vdash v) &= (x \triangleleft y) \dashv \{z,u,v\}_3 = \{x,y,\{z,u,v\}_3\}_1
            \end{align*}
            
            \item $((x \triangleleft y) \vdash (z \triangleright u)) \dashv v = (x \triangleleft y) \vdash ( (z \triangleright u) \dashv v)$:
            \begin{align*}
                ((x \triangleleft y) \vdash (z \triangleright u)) \dashv v &= (\{x,y,z\}_3 \triangleright u) \dashv v = \{\{x,y,z\}_3,u,v\}_2
                \\
                (x \triangleleft y) \vdash ( (z \triangleright u) \dashv v) &= (x \triangleleft y) \vdash \{z,u,v\}_2 = \{x,y,\{z,u,v\}_2\}_3
            \end{align*}
            
            \item $((x \triangleleft y) \vdash (z \triangleright u)) \vdash v = ((x \triangleleft y) \dashv (z \triangleright u)) \vdash v$:
            \begin{align*}
                ((x \triangleleft y) \vdash (z \triangleright u)) \vdash v &= (\{x,y,z\}_3 \triangleright u) \vdash v = \{\{x,y,z\}_3,u,v\}_3
                \\
                ((x \triangleleft y) \dashv (z \triangleright u)) \vdash v &= (x \triangleleft \{y,z,u\}_1) \vdash v = \{x,\{y,z,u\}_1,v\}_3
            \end{align*}
        \end{enumerate}

        \item Let us consider $\alpha = x \triangleright y \in \mathfrak{M}(A,A), \beta = z \triangleright u \in \mathfrak{M}(A,A)$ and $\gamma = v \in A$. Let us show that:
        \begin{enumerate}
            \item $((x \triangleright y) \dashv (z \triangleright u)) \dashv v = (x \triangleright y) \dashv ( (z \triangleright u) \dashv v)$:
            \begin{align*}
                ((x \triangleright y) \dashv (z \triangleright u)) \dashv v &= (x \triangleright \{y,z,u\}_1)\dashv v = \{x,\{y,z,u\}_1,v\}_2
                \\
                (x \triangleright y) \dashv ( (z \triangleright u) \dashv v) &= ( x \triangleright y) \dashv \{z,u,v\}_2 = \{x,y,\{z,u,v\}_2\}_2
            \end{align*}
            
            \item $((x \triangleright y) \vdash (z \triangleright u)) \vdash v = (x \triangleright y) \vdash ( (z \triangleright u) \vdash v)$:
            \begin{align*}
                ((x \triangleright y) \vdash (z \triangleright u)) \vdash v &= ( x \triangleright \{y,z,u\}_3) \vdash v = \{x,\{y,z,u\}_3,v\}_3
                \\
                (x \triangleright y) \vdash ( (z \triangleright u) \vdash v) &=(x \triangleright y) \vdash \{z,u,v\}_3 = \{x,y,\{z,u,v\}_3\}_3
            \end{align*}
            
            \item $(x \triangleright y) \dashv ((z \triangleright u) \dashv v) = (x \triangleright y) \dashv ( (z \triangleright u) \vdash v)$:
            \begin{align*}
                (x \triangleright y) \dashv ((z \triangleright u) \dashv v) &= (x \triangleright y) \dashv \{z,u,v\}_2 = \{x,y,\{z,u,v\}_2\}_2
                \\
                (x \triangleright y) \dashv ( (z \triangleright u) \vdash v) &= (x \triangleright y) \dashv  \{z,u,v\}_3 =\{x,y,\{z,u,v\}_3\}_2
            \end{align*}
            
            \item $((x \triangleright y) \vdash (z \triangleright u)) \dashv v = (x \triangleright y) \vdash ( (z \triangleright u) \dashv v)$:
            \begin{align*}
                ((x \triangleright y) \vdash (z \triangleright u)) \dashv v &= (\{x,y,z\}_3\triangleright u) \dashv v = \{\{x,y,z\}_3,u,v\}_2
                \\
                (x \triangleright y) \vdash ( (z \triangleright u) \dashv v) &= (x \triangleright y) \vdash  \{z,u,v\}_2 = \{x,y,\{z,u,v\}_2\}_3
            \end{align*}
            
            \item $((x \triangleright y) \vdash (z \triangleright u)) \vdash v = ((x \triangleright y) \dashv (z \triangleright u)) \vdash v$:
            \begin{align*}
                ((x \triangleright y) \vdash (z \triangleright u)) \vdash v &= (\{x,y,z\}_3 \triangleright u) \vdash v = \{\{x,y,z\}_3,u,v\}_3
                \\
                ((x \triangleright y) \dashv (z \triangleright u)) \vdash v &= (x \triangleright \{y,z,u\}_1) \vdash v = \{x,\{y,z,u\}_1,v\}_3
            \end{align*}
        \end{enumerate}
        
        \item Let us consider $\alpha = x \triangleleft y \in \mathfrak{M}(A,A), \beta = z \in A$ and $\gamma = u \triangleleft v \in \mathfrak{M}(A,A)$. Let us show that:
        \begin{enumerate}
            \item $((x \triangleleft y) \dashv z) \dashv (u \triangleleft v) = (x \triangleleft y) \dashv ( z \dashv (u \triangleleft v))$:
            \begin{align*}
                ((x \triangleleft y) \dashv z) \dashv (u \triangleleft v) &= \{x,y,z\}_1 \dashv ( u \triangleleft v) = \{\{x,y,z\}_1,u,v\}_1
                \\
                (x \triangleleft y) \dashv ( z \dashv (u \triangleleft v)) &=(x \triangleleft y) \dashv  \{z,u,v\}_1 = \{x,y,\{z,u,v\}_1\}_1
            \end{align*}
            
            \item $((x \triangleleft y) \vdash z) \vdash (u \triangleleft v) = (x \triangleleft y) \vdash ( z \vdash (u \triangleleft v))$:
            \begin{align*}
                ((x \triangleleft y) \vdash z) \vdash (u \triangleleft v) &= \{x,y,z\}_3 \vdash ( u \triangleleft v) = \{\{x,y,z\}_3,u,v\}_2
                \\
                (x \triangleleft y) \vdash ( z \vdash (u \triangleleft v)) & = (x \triangleleft y) \vdash \{z,u,v\}_2 = \{x,y,\{z,u,v\}_2\}_3
            \end{align*}
            
            \item $(x \triangleleft y) \dashv (z \dashv (u \triangleleft v)) = (x \triangleleft y) \dashv ( z \vdash (u \triangleleft v))$:
            \begin{align*}
                (x \triangleleft y) \dashv (z \dashv (u \triangleleft v)) &= (x \triangleleft y) \dashv \{z,u,v\}_1 = \{x,y,\{z,u,v\}_1\}_1 
                \\
                (x \triangleleft y) \dashv ( z \vdash (u \triangleleft v)) &= (x \triangleleft y) \dashv \{z,u,v\}_2 = \{x,y,\{z,u,v\}_2\}_1
            \end{align*}
            
            \item $((x \triangleleft y) \vdash z) \dashv (u \triangleleft v) = (x \triangleleft y) \vdash ( z \dashv (u \triangleleft v))$:
            \begin{align*}
                ((x \triangleleft y) \vdash z) \dashv (u \triangleleft v) &= \{x,y,z\}_3 \dashv ( u \triangleleft v) = \{\{x,y,z\}_3,u,v\}_1
                \\
                (x \triangleleft y) \vdash ( z \dashv (u \triangleleft v))&= ( x\triangleleft y) \vdash \{z,u,v\}_1 = \{x,y,\{z,u,v\}_1\}_3
            \end{align*}
            
            \item $((x \triangleleft y) \vdash z) \vdash (u \triangleleft v) = ((x \triangleleft y) \dashv z) \vdash (u \triangleleft v)$:
            \begin{align*}
                ((x \triangleleft y) \vdash z) \vdash (u \triangleleft v) &= \{x,y,z\}_3 \vdash (u \triangleleft v) = \{\{x,y,z\}_3,u,v\}_2
                \\
                ((x \triangleleft y) \dashv z) \vdash (u \triangleleft v) &= \{x,y,z\}_1 \vdash (u \triangleleft v) = \{\{x,y,z\}_1,u,v\}_2
            \end{align*}
        \end{enumerate}

        \item Let us consider $\alpha = x \triangleright y \in \mathfrak{M}(A,A), \beta = z \in A$ and $\gamma = u \triangleleft v \in \mathfrak{M}(A,A)$. Let us show that:
        \begin{enumerate}
            \item $((x \triangleright y) \dashv z) \dashv (u \triangleleft v) = (x \triangleright y) \dashv ( z \dashv (u \triangleleft v))$:
            \begin{align*}
                ((x \triangleright y) \dashv z) \dashv (u \triangleleft v) & = \{x,y,z\}_2 \dashv ( u \triangleleft v) = \{\{x,y,z\}_2 ,u,v\}_1
                \\
                (x \triangleright y) \dashv ( z \dashv (u \triangleleft v)) &= (x \triangleright y) \dashv \{z,u,v\}_1 = \{x,y,\{z,u,v\}_1\}_2
            \end{align*}
            
            \item $((x \triangleright y) \vdash z) \vdash (u \triangleleft v) = (x \triangleright y) \vdash ( z \vdash (u \triangleleft v))$:
            \begin{align*}
                ((x \triangleright y) \vdash z) \vdash (u \triangleleft v) &=\{x,y,z\}_3 \vdash (u \triangleleft v) = \{\{x,y,z\}_3,u,v\}_2
                \\
                (x \triangleright y) \vdash ( z \vdash (u \triangleleft v)) &= (x \triangleright y) \vdash \{z,u,v\}_2 = \{x,y,\{z,u,v\}_2\}_3
            \end{align*}
            
            \item $(x \triangleright y) \dashv (z \dashv (u \triangleleft v)) = (x \triangleright y) \dashv ( z \vdash (u \triangleleft v))$:
            \begin{align*}
                (x \triangleright y) \dashv (z \dashv (u \triangleleft v))  &= (x \triangleright y) \dashv \{z,u,v\}_1 = \{x,y,\{z,u,v\}_1\}_2 
                \\
                (x \triangleright y) \dashv ( z \vdash (u \triangleleft v)) &= (x \triangleright y) \dashv \{z,u,v\}_2 = \{x,y,\{z,u,v\}_2\}_2
            \end{align*}
            
            \item $((x \triangleright y) \vdash z) \dashv (u \triangleleft v) = (x \triangleright y) \vdash ( z \dashv (u \triangleleft v))$:
            \begin{align*}
                ((x \triangleright y) \vdash z) \dashv (u \triangleleft v)  &= \{x,y,z\}_3 \dashv (u \triangleleft v) = \{\{x,y,z\}_3,u,v\}_1
                \\
                (x \triangleright y) \vdash ( z \dashv (u \triangleleft v)) & = (x \triangleright y) \vdash \{z,u,v\}_1 = \{x,y,\{z,u,v\}_1\}_3
            \end{align*}
            
            \item $((x \triangleright y) \vdash z) \vdash (u \triangleleft v) = ((x \triangleright y) \dashv z) \vdash (u \triangleleft v)$:
            \begin{align*}
                ((x \triangleright y) \vdash z) \vdash (u \triangleleft v) &= \{x,y,z\}_3 \vdash (u \triangleleft v) = \{\{x,y,z\}_3,u,v\}_2 
                \\
                 ((x \triangleright y) \dashv z) \vdash (u \triangleleft v) &= \{x,y,z\}_2 \vdash (u \triangleleft v) = \{\{x,y,z\}_2,u,v\}_2
            \end{align*}
        \end{enumerate}

        \item Let us consider $\alpha = x \triangleleft y \in \mathfrak{M}(A,A), \beta = z \in A$ and $\gamma = u \triangleright v \in \mathfrak{M}(A,A)$. Let us show that:
        \begin{enumerate}
            \item $((x \triangleleft y) \dashv z) \dashv (u \triangleright v) = (x \triangleleft y) \dashv ( z \dashv (u \triangleright v))$:
            \begin{align*}
                ((x \triangleleft y) \dashv z) \dashv (u \triangleright v) &= \{x,y,z\}_1 \dashv (u \triangleright v) = \{\{x,y,z\}_1,u,v\}_1
                \\
                (x \triangleleft y) \dashv ( z \dashv (u \triangleright v)) & = (x \triangleleft y) \dashv \{z,u,v\}_1 = \{x,y,\{z,u,v\}_1\}_1
            \end{align*}
            
            \item $((x \triangleleft y) \vdash z) \vdash (u \triangleright v) = (x \triangleleft y) \vdash ( z \vdash (u \triangleright v))$:
            \begin{align*}
                ((x \triangleleft y) \vdash z) \vdash (u \triangleright v) &= \{x,y,z\}_3 \vdash ( u \triangleright v) = \{\{x,y,z\}_3,u,v\}_3
                \\
                (x \triangleleft y) \vdash ( z \vdash (u \triangleright v)) &= (x \triangleleft y) \vdash \{z,u,v\}_3 = \{x,y,\{z,u,v\}_3\}_3
            \end{align*}
            
            \item $(x \triangleleft y) \dashv (z \dashv (u \triangleright v)) = (x \triangleleft y) \dashv ( z \vdash (u \triangleright v))$:
            \begin{align*}
                (x \triangleleft y) \dashv (z \dashv (u \triangleright v)) &= (x \triangleleft y) \dashv \{z,u,v\}_1 = \{x,y,\{z,u,v\}_1\}_1 
                \\
                (x \triangleleft y) \dashv ( z \vdash (u \triangleright v)) &= (x \triangleleft y) \dashv \{z,u,v\}_3 = \{x,y,\{z,u,v\}_3\}_1
            \end{align*}
            
            \item $((x \triangleleft y) \vdash z) \dashv (u \triangleright v) = (x \triangleleft y) \vdash ( z \dashv (u \triangleright v))$:
            \begin{align*}
                ((x \triangleleft y) \vdash z) \dashv (u \triangleright v) &= \{x,y,z\}_3 \dashv ( u \triangleright v) = \{\{x,y,z\}_3,u,v\}_1
                \\
                (x \triangleleft y) \vdash ( z \dashv (u \triangleright v)) &= (x \triangleleft y) \vdash \{z,u,v\}_1 = \{x,y,\{z,u,v\}_1\}_3
            \end{align*}
            
            \item $((x \triangleleft y) \vdash z) \vdash (u \triangleright v) = ((x \triangleleft y) \dashv z) \vdash (u \triangleright v)$:
            \begin{align*}
                ((x \triangleleft y) \vdash z) \vdash (u \triangleright v) &= \{x,y,z\}_3 \vdash (u \triangleright v) = \{\{x,y,z\}_3,u,v\}_3
                \\
                ((x \triangleleft y) \dashv z) \vdash (u \triangleright v) &= \{x,y,z\}_1 \vdash (u \triangleright v) = \{\{x,y,z\}_1,u,v\}_3
            \end{align*}
        \end{enumerate}

        \item Let us consider $\alpha = x \triangleright y \in \mathfrak{M}(A,A), \beta = z \in A$ and $\gamma = u \triangleright v \in \mathfrak{M}(A,A)$. Let us show that:
        \begin{enumerate}
            \item $((x \triangleright y) \dashv z) \dashv (u \triangleright v) = (x \triangleright y) \dashv ( z \dashv (u \triangleright v))$:
            \begin{align*}
                ((x \triangleright y) \dashv z) \dashv (u \triangleright v) & =  \{x,y,z\}_2 \dashv ( u \triangleright v) = \{\{x,y,z\}_2,u,v\}_1
                \\
                (x \triangleright y) \dashv ( z \dashv (u \triangleright v)) &= (x \triangleright y) \dashv \{z,u,v\}_1 = \{x,y,\{z,u,v\}_1\}_2
            \end{align*}
            
            \item $((x \triangleright y) \vdash z) \vdash (u \triangleright v) = (x \triangleright y) \vdash ( z \vdash (u \triangleright v))$:
            \begin{align*}
                ((x \triangleright y) \vdash z) \vdash (u \triangleright v)  &= \{x,y,z\}_3 \vdash (u \triangleright v) = \{\{x,y,z\}_3,u,v\}_3
                \\
                (x \triangleright y) \vdash ( z \vdash (u \triangleright v)) &= (x \triangleright y) \vdash \{z,u,v\}_3 = \{x,y,\{z,u,v\}_3\}_3
            \end{align*}
            
            \item $(x \triangleright y) \dashv (z \dashv (u \triangleright v)) = (x \triangleright y) \dashv ( z \vdash (u \triangleright v))$:
            \begin{align*}
                (x \triangleright y) \dashv (z \dashv (u \triangleright v)) &= ( x\triangleright y)\dashv \{z,u,v\}_1 = \{x,y,\{z,u,v\}_1\}_2
                \\
                (x \triangleright y) \dashv ( z \vdash (u \triangleright v)) &= ( x\triangleright y)\dashv\{z,u,v\}_3 = \{x,y,\{z,u,v\}_3\}_2
            \end{align*}
            
            \item $((x \triangleright y) \vdash z) \dashv (u \triangleright v) = (x \triangleright y) \vdash ( z \dashv (u \triangleright v))$:
            \begin{align*}
                ((x \triangleright y) \vdash z) \dashv (u \triangleright v) &= \{x,y,z\}_3 \dashv (u \triangleright v) = \{\{x,y,z\}_3,u,v\}_1
                \\
                (x \triangleright y) \vdash ( z \dashv (u \triangleright v)) &= (x \triangleright y) \vdash \{z,u,v\}_1 = \{x,y,\{z,u,v\}_1\}_3
            \end{align*}
            
            \item $((x \triangleright y) \vdash z) \vdash (u \triangleright v) = ((x \triangleright y) \dashv z) \vdash (u \triangleright v)$:
            \begin{align*}
                ((x \triangleright y) \vdash z) \vdash (u \triangleright v) &= \{x,y,z\}_3 \vdash (u \triangleright v) = \{\{x,y,z\}_3,u,v\}_3
                \\
                ((x \triangleright y) \dashv z) \vdash (u \triangleright v) &= \{x,y,z\}_2 \vdash ( u \triangleright v) = \{\{x,y,z\}_2,u,v\}_3
            \end{align*}
        \end{enumerate}
    
        \item Let us consider $\alpha = x  \in A, \beta = y \triangleleft z \in \mathfrak{M}(A,A)$ and $\gamma = u \triangleleft v \in \mathfrak{M}(A,A)$. Let us show that:
        \begin{enumerate}
            \item $(x \dashv (y \triangleleft z)) \dashv (u \triangleleft v) = x \dashv ( (y \triangleleft z) \dashv (u \triangleleft v))$:
            \begin{align*}
                (x \dashv (y \triangleleft z)) \dashv (u \triangleleft v) &= \{x,y,z\}_1 \dashv (u \triangleleft v) = \{\{x,y,z\}_1,u,v\}_1
                \\
                x \dashv ( (y \triangleleft z) \dashv (u \triangleleft v)) & = x \dashv (y \triangleleft \{z,u,v\}_1) = \{x,y,\{z,u,v\}_1\}_1
            \end{align*}
            
            \item $(x \vdash (y \triangleleft z)) \vdash (u \triangleleft v) = x \vdash ( (y \triangleleft z) \vdash (u \triangleleft v))$:
            \begin{align*}
                (x \vdash (y \triangleleft z)) \vdash (u \triangleleft v) & = \{x,y,z\}_2 \vdash (u \triangleleft v) = \{\{x,y,z\}_2,u,v\}_2
                \\
                x \vdash ( (y \triangleleft z) \vdash (u \triangleleft v)) &= x \vdash (\{y,z,u\}_3 \triangleleft v) = \{x,\{y,z,u\}_3,v\}_2
            \end{align*}
            
            \item $x \dashv ((y \triangleleft z) \dashv (u \triangleleft v)) = x \dashv ( (y \triangleleft z) \vdash (u \triangleleft v))$:
            \begin{align*}
                x \dashv ((y \triangleleft z) \dashv (u \triangleleft v)) &= x \dashv (y \triangleleft \{z,u,v\}_1) = \{x,y,\{z,u,v\}_1\}_1
                \\
                x \dashv ( (y \triangleleft z) \vdash (u \triangleleft v)) &= x \dashv (\{y,z,u\}_1 \triangleleft v) = \{x,\{y,z,u\}_1,v\}_1
            \end{align*}
            
            \item $(x \vdash (y \triangleleft z)) \dashv (u \triangleleft v) = x \vdash ( (y \triangleleft z) \dashv (u \triangleleft v))$:
            \begin{align*}
                (x \vdash (y \triangleleft z)) \dashv (u \triangleleft v) &= \{x,y,z\}_2 \dashv (u \triangleleft v) = \{\{x,y,z\}_2,u,v\}_1
                \\
                x \vdash ( (y \triangleleft z) \dashv (u \triangleleft v)) & = x \vdash (y \triangleleft \{z,u,v\}_1) = \{x,y,\{z,u,v\}_1\}_2
            \end{align*}
            
            \item $(x \vdash (y \triangleleft z)) \vdash (u \triangleleft v) = (x \dashv (y \triangleleft z)) \vdash (u \triangleleft v)$:
            \begin{align*}
                (x \vdash (y \triangleleft z)) \vdash (u \triangleleft v) &= \{x,y,z\}_2 \vdash (u \triangleleft v) =\{\{x,y,z\}_2,u,v\}_2
                \\
                (x \dashv (y \triangleleft z)) \vdash (u \triangleleft v) &=\{x,y,z\}_1 \vdash (u \triangleleft v) = \{\{x,y,z\}_1,u,v\}_2
            \end{align*}
        \end{enumerate}
    
        \item Let us consider $\alpha = x  \in A, \beta = y \triangleright z \in \mathfrak{M}(A,A)$ and $\gamma = u \triangleleft v \in \mathfrak{M}(A,A)$. Let us show that:
        \begin{enumerate}
            \item $(x \dashv (y \triangleright z)) \dashv (u \triangleleft v) = x \dashv ( (y \triangleright z) \dashv (u \triangleleft v))$:
            \begin{align*}
                (x \dashv (y \triangleright z)) \dashv (u \triangleleft v) &= \{x,y,z\}_1 \dashv (u \triangleleft v) = \{\{x,y,z\}_1,u,v\}_1
                \\
                x \dashv ( (y \triangleright z) \dashv (u \triangleleft v)) &= x \dashv (y \triangleright \{z,u,v\}_1) = \{x,y,\{z,u,v\}_1\}_1
            \end{align*}
            
            \item $(x \vdash (y \triangleright z)) \vdash (u \triangleleft v) = x \vdash ( (y \triangleright z) \vdash (u \triangleleft v))$:
            \begin{align*}
                (x \vdash (y \triangleright z)) \vdash (u \triangleleft v) &= \{x,y,z\}_3 \vdash (u \triangleleft v) = \{\{x,y,z\}_3,u,v\}_2
                \\
                x \vdash ( (y \triangleright z) \vdash (u \triangleleft v)) &= x \vdash (\{y,z,u\}_3\triangleleft v) = \{x,\{y,z,u\}_3,v\}_2
            \end{align*}
            
            \item $x \dashv ((y \triangleright z) \dashv (u \triangleleft v)) = x \dashv ( (y \triangleright z) \vdash (u \triangleleft v))$:
            \begin{align*}
                x \dashv ((y \triangleright z) \dashv (u \triangleleft v)) &= x \dashv (y \triangleright \{z,u,v\}_1) = \{x,y,\{z,u,v\}_1\}_1
                \\
                x \dashv ( (y \triangleright z) \vdash (u \triangleleft v)) &= x \dashv ( \{y,z,u\}_3 \triangleleft v) = \{x,\{y,z,u\}_3,v\}_1
            \end{align*}
            
            \item $(x \vdash (y \triangleright z)) \dashv (u \triangleleft v) = x \vdash ( (y \triangleright z) \dashv (u \triangleleft v))$:
            \begin{align*}
                (x \vdash (y \triangleright z)) \dashv (u \triangleleft v)  &= \{x,y,z\}_3 \dashv (u \triangleleft v) = \{\{x,y,z\}_3,u,v\}_1
                \\
                x \vdash ( (y \triangleright z) \dashv (u \triangleleft v)) &= x \vdash (y \triangleright \{z,u,v\}_1) = \{x,y,\{z,u,v\}_1\}_3
            \end{align*}
            
            \item $(x \vdash (y \triangleright z)) \vdash (u \triangleleft v) = (x \dashv (y \triangleright z)) \vdash (u \triangleleft v)$:
            \begin{align*}
                (x \vdash (y \triangleright z)) \vdash (u \triangleleft v) & = \{x,y,z\}_3 \vdash ( u \triangleleft v) = \{\{x,y,z\}_3,u,v\}_2
                \\
                (x \dashv (y \triangleright z)) \vdash (u \triangleleft v)& = \{x,y,z\}_1 \vdash (u \triangleleft v) = \{\{x,y,z\}_1,u,v\}_2
            \end{align*}
        \end{enumerate}
    
        \item Let us consider $\alpha = x  \in A, \beta = y \triangleleft z \in \mathfrak{M}(A,A)$ and $\gamma = u \triangleright v \in \mathfrak{M}(A,A)$. Let us show that:
        \begin{enumerate}
            \item $(x \dashv (y \triangleleft z)) \dashv (u \triangleright v) = x \dashv ( (y \triangleleft z) \dashv (u \triangleright v))$:
            \begin{align*}
                (x \dashv (y \triangleleft z)) \dashv (u \triangleright v) &= \{x,y,z\}_1 \dashv (u \triangleright v) = \{\{x,y,z\}_1,u,v\}_1
                \\
                x \dashv ( (y \triangleleft z) \dashv (u \triangleright v)) &= x \dashv (y \triangleleft \{z,u,v\}_1) = \{x,y,\{z,u,v\}_1\}_1
            \end{align*}
            
            \item $(x \vdash (y \triangleleft z)) \vdash (u \triangleright v) = x \vdash ( (y \triangleleft z) \vdash (u \triangleright v))$:
            \begin{align*}
                (x \vdash (y \triangleleft z)) \vdash (u \triangleright v) &=\{x,y,z\}_2 \vdash (u \triangleright v) = \{\{x,y,z\}_2,u,v\}_3
                \\
                x \vdash ( (y \triangleleft z) \vdash (u \triangleright v)) & = x \vdash ( \{y,z,u\}_3 \triangleright v) = \{x,\{y,z,u\}_3,v\}_3
            \end{align*}
            
            \item $x \dashv ((y \triangleleft z) \dashv (u \triangleright v)) = x \dashv ( (y \triangleleft z) \vdash (u \triangleright v))$:
            \begin{align*}
                x \dashv ((y \triangleleft z) \dashv (u \triangleright v)) &= x \dashv ( y \triangleleft \{z,u,v\}_1) = \{x,y,\{z,u,v\}_1\}_1
                \\
                x \dashv ( (y \triangleleft z) \vdash (u \triangleright v)) & = x \dashv (\{y,z,u\}_3\triangleright v) = \{x,\{y,z,u\}_3,v\}_1
            \end{align*}
            
            \item $(x \vdash (y \triangleleft z)) \dashv (u \triangleright v) = x \vdash ( (y \triangleleft z) \dashv (u \triangleright v))$:
            \begin{align*}
                (x \vdash (y \triangleleft z)) \dashv (u \triangleright v) & = \{x,y,z\}_2 \dashv ( u \triangleright v) = \{\{x,y,z\}_2,u,v\}_1
                \\
                x \vdash ( (y \triangleleft z) \dashv (u \triangleright v)) & = x \vdash ( y \triangleleft \{z,u,v\}_1) = \{x,y,\{z,u,v\}_1\}_2
            \end{align*}
            
            \item $(x \vdash (y \triangleleft z)) \vdash (u \triangleright v) = (x \dashv (y \triangleleft z)) \vdash (u \triangleright v)$:
            \begin{align*}
                (x \vdash (y \triangleleft z)) \vdash (u \triangleright v) &= \{x,y,z\}_2 \vdash (u \triangleright v) = \{\{x,y,z\}_2,u,v\}_3
                \\
                (x \dashv (y \triangleleft z)) \vdash (u \triangleright v) & = \{x,y,z\}_1 \vdash ( u \triangleright v) = \{\{x,y,z\}_1,u,v\}_3
            \end{align*}
        \end{enumerate}
        
        \item Let us consider $\alpha = x  \in A, \beta = y \triangleright z \in \mathfrak{M}(A,A)$ and $\gamma = u \triangleright v \in \mathfrak{M}(A,A)$. Let us show that:
        \begin{enumerate}
            \item $(x \dashv (y \triangleright z)) \dashv (u \triangleright v) = x \dashv ( (y \triangleright z) \dashv (u \triangleright v))$:
            \begin{align*}
                (x \dashv (y \triangleright z)) \dashv (u \triangleright v) &=\{x,y,z\}_1 \dashv (u \triangleright v) = \{\{x,y,z\}_1,u,v\}_1
                \\
                x \dashv ( (y \triangleright z) \dashv (u \triangleright v)) &= x \dashv( y \triangleright \{z,u,v\}_1) = \{x,y,\{z,u,v\}_1\}_1
            \end{align*}
            
            \item $(x \vdash (y \triangleright z)) \vdash (u \triangleright v) = x \vdash ( (y \triangleright z) \vdash (u \triangleright v))$:
            \begin{align*}
                (x \vdash (y \triangleright z)) \vdash (u \triangleright v) &= \{x,y,z\}_3 \vdash ( u \triangleright v) = \{\{x,y,z\}_3,u,v\}_3
                \\
                x \vdash ( (y \triangleright z) \vdash (u \triangleright v)) &= x \vdash (\{y,z,u\}_3\triangleright v) = \{x,\{y,z,u\}_3,v\}_3
            \end{align*}
            
            \item $x \dashv ((y \triangleright z) \dashv (u \triangleright v)) = x \dashv ( (y \triangleright z) \vdash (u \triangleright v))$:
            \begin{align*}
                x \dashv ((y \triangleright z) \dashv (u \triangleright v)) &= x \dashv (y \triangleright \{z,u,v\}_1) = \{x,y,\{z,u,v\}_1\}_1
                \\
                x \dashv ( (y \triangleright z) \vdash (u \triangleright v)) &=x \dashv (\{y,z,u\}_3 \triangleright v) = \{x,\{y,z,u\}_3,v\}_1
            \end{align*}
            
            \item $(x \vdash (y \triangleright z)) \dashv (u \triangleright v) = x \vdash ( (y \triangleright z) \dashv (u \triangleright v))$:
            \begin{align*}
                (x \vdash (y \triangleright z)) \dashv (u \triangleright v) &= \{x,y,z\}_3 \dashv (u \triangleright v) = \{\{x,y,z\}_3,u,v\}_1
                \\
                x \vdash ( (y \triangleright z) \dashv (u \triangleright v)) &= x \vdash ( y \triangleright \{z,u,v\}_1) = \{x,y,\{z,u,v\}_1\}_3
            \end{align*}
            
            \item $(x \vdash (y \triangleright z)) \vdash (u \triangleright v) = (x \dashv (y \triangleright z)) \vdash (u \triangleright v)$:
            \begin{align*}
                (x \vdash (y \triangleright z)) \vdash (u \triangleright v) &= \{x,y,z\}_3 \vdash (u \triangleright v) = \{\{x,y,z\}_3,u,v\}_3
                \\
                (x \dashv (y \triangleright z)) \vdash (u \triangleright v) &= \{x,y,z\}_1 \vdash (u \triangleright v) = \{\{x,y,z\}_1,u,v\}_3
            \end{align*}
        \end{enumerate}
    \end{enumerate}
    Reciprocally,  if $\mathrm{U}(A)$ is an associative dialgebra, we recover the relations
    \begin{multicols}{2}
        \begin{itemize}
            \item (\ref{cor2}) from 3(a) and 7(a),
            \item (\ref{cor3}) from 4(a) and 8(a),
            \item (\ref{cor4}) from Proposition \ref{lemma:products}(2),
            \item (\ref{cor5}) from 5(d) and 12(b),
            \item (\ref{cor6}) from 3(b) and 6(b),
            \item (\ref{cor7}) from 7(e) and 8(e),
            \item (\ref{cor8}) from 9(e) and 10(e),
            \item (\ref{cor9}) from Proposition \ref{lemma:products} (2-3),
            \item (\ref{cor10}) from Proposition \ref{lemma:products}(1-2),
            \item (\ref{cor11}) from 3(c) and 5(c),
            \item (\ref{cor12}) from 4(c) and 6(c).
        \end{itemize}
    \end{multicols}
\end{proof}

\subsection{Standard embedding for associative triple trisystems of the second kind}

Let $A$ be an associative triple trisystem of the second kind. We recall the $\phi$-modules $\mathfrak{L}(A,A), \mathfrak{R}(A,A) \subset \mathrm{DiEnd}(A)$  generated by the di-endomorphisms
$$
    L^\triangleleft(x,y) = (L_1(x,y), L_3(x,y)), \quad L^\triangleright (x,y) = (L_2(x,y), L_3(x,y)),
$$
and 
$$
    R^\triangleleft (x,y) = (R_2(x,y),R_1(x,y)), \quad  R^\triangleright(x,y) = (R_3(x,y),R_1(x,y)),
$$
for any $x,y \in A$, respectively.
\begin{remark}\label{rem:asssecond}
    Observe that, for any $x,y \in A$, $\lambda \in \mathfrak{L}(A,A)$ and $\rho \in \mathfrak{R}(A,A)$, we have the following relations
    $$
        \lambda \dashv L^\triangleleft(z,u) = \lambda \dashv L^\triangleright(z,u), \quad L^\triangleleft(x,y) \vdash \lambda = L^\triangleright(x,y) \vdash \lambda,
    $$
    and
    $$
        \rho \dashv R^\triangleleft(z,u) = \rho \dashv R^\triangleright(z,u), \quad R^\triangleleft(x,y) \vdash \rho = R^\triangleright(x,y) \vdash \rho.
    $$
\end{remark}
In the next proposition, we prove that $\mathfrak{L}(A,A)$ and $\mathfrak{R}(A,A)^\mathrm{op}$ are indeed closed under the products $\dashv$ and $\vdash$, so they are subdialgebras of $\mathrm{DiEnd}(A)$ and $\mathrm{DiEnd}(A)^\mathrm{op}$ respectively.

\begin{proposition}\label{prop:asssecondkind}
    Let $A$ be an associative triple trisystem of the second kind. Then the products in $\mathfrak{L}(A,A)$ and $\mathfrak{R}(A,A)^\mathrm{op}$ satisfy
    \begin{enumerate}
        \item $L^\triangleleft(x,y) \dashv L^\triangleleft(z,u) = L^\triangleleft(x,\{u,z,y\}_i) = L^\triangleleft(\{x,y,z\}_1,u)$, 
        
        \item $L^\triangleright(x,y) \dashv L^\triangleleft(z,u) = L^\triangleright(x,\{u,z,y\}_3) =L^\triangleleft(\{x,y,z\}_2,u)$,
        
        \item $L^\triangleleft(x,y) \vdash L^\triangleleft (z,u) = L^\triangleright(x,\{u,z,y\}_2)= L^\triangleleft(\{x,y,z\}_3,u)$,
        
        \item $L^\triangleleft(x,y) \vdash L^\triangleright(z,u) = L^\triangleright(x,\{u,z,y\}_1) = L^\triangleright(\{x,y,z\}_i,u)$,
        
        \item $R^\triangleleft(x,y) \dashv R^\triangleleft(z,u) = R^\triangleleft(\{z,y,x\}_3,u)= R^\triangleleft(x,\{y,z,u\}_i)$,
        
        \item $   R^\triangleright(x,y) \dashv R^\triangleleft(z,u)= R^\triangleleft(\{z,y,x\}_2,u) = R^\triangleright(x,\{y,z,u\}_1)$,
        
        \item $R^\triangleleft (x,y) \vdash R^\triangleleft (z,u) = R^\triangleleft(\{z,y,x\}_1,u) = R^\triangleright(x,\{y,z,u\}_2)$,
        
        \item $R^\triangleleft (x,y) \vdash R^\triangleright(z,u) =R^\triangleright(\{z,y,x\}_i,u) = R^\triangleright(x,\{y,z,u\}_3)$,
    \end{enumerate}
    for every $x,y,z,u \in A$ and $i \in \{1,2,3\}$. In particular, they both are associative dialgebras.
\end{proposition}
\begin{proof}
    In order to simplify the proof, we denote by $\equiv$ when we evaluate the di-homomorphism. For instance, for any $x,y,z \in A$,
    $$
        L^\triangleleft(x,y) \equiv (L^\triangleleft(x,y) \prec z, L^\triangleleft(x,y) \succ z ) = (\{x,y,z\}_1,\{x,y,z\}_3)
    $$
    and
    $$
        R^\triangleright(x,y) \equiv (z \succ R^\triangleright(x,y), z \prec R^\triangleright(x,y))= ( \{z,x,y\}_3,\{z,x,y\}_1).
    $$
    Given $x,y,z,u,v \in A$,
    \begin{enumerate}
        \item 
        \begin{align*}
            L^\triangleleft(x,y) \dashv L^\triangleleft(z,u) &= (L_1(x,y) \circ L_3(z,u), L_3(x,y) \circ L_3(z,u))
            \\
            &\equiv (\{x,y,\{z,u,v\}_3\}_1, \{x,y,\{z,u,v\}_3\}_3)
            \\
            &=(\{x,\{u,z,y\}_i,v\}_1,\{x,\{u,z,y\}_i,v\}_3) \equiv L^\triangleleft(x,\{u,z,y\}_i)
            \\
            &=(\{\{x,y,z\}_1,u,v\}_1,\{\{x,y,z\}_1,u,v\}_3) \equiv L^\triangleleft(\{x,y,z\}_1,u)
        \end{align*}

        \item 
        \begin{align*}
            L^\triangleright(x,y) \dashv L^\triangleleft(z,u) &= (L_2(x,y) \circ L_3(z,u), L_3(x,y) \circ L_3(z,u))
            \\
            &\equiv(\{x,y,\{z,u,v\}_3\}_2, \{x,y,\{z,u,v\}_3\}_3)
            \\
            &=( \{x,\{u,z,y\}_3,v\}_2,\{x,\{u,z,y\}_3,v\}_3) \equiv L^\triangleright(x,\{u,z,y\}_3) 
            \\
            &=(\{\{x,y,z\}_2,u,v\}_1, \{\{x,y,z\}_2,u,v\}_3) \equiv L^\triangleleft(\{x,y,z\}_2,u)
        \end{align*}
        
        \item 
        \begin{align*}
            L^\triangleleft(x,y) \vdash L^\triangleleft (z,u) &= (L_3(x,y) \circ L_1(z,u), L_3(x,y) \circ L_3(z,u))
            \\
            &\equiv(\{x,y,\{z,u,v\}_1\}_3, \{x,y,\{z,u,v\}_3\}_3)
            \\
            &=(\{x,\{u,z,y\}_2,v\}_2, \{x,\{u,z,y\}_2,v\}_3) \equiv L^\triangleright(x,\{u,z,y\}_2)
            \\
            &=(\{\{x,y,z\}_3,u,v\}_1, \{\{x,y,z\}_3,u,v\}_3) \equiv L^\triangleleft(\{x,y,z\}_3,u)
        \end{align*}

        \item 
        \begin{align*}
            L^\triangleleft(x,y) \vdash L^\triangleright(z,u) &=  (L_3(x,y) \circ L_2(z,u), L_3(x,y) \circ L_3(z,u))
            \\
            &\equiv(\{x,y,\{z,u,v\}_2\}_3, \{x,y,\{z,u,v\}_3\}_3)
            \\
            &=(\{x,\{u,z,y\}_1,v\}_2, \{x,\{u,z,y\}_1,v\}_3) \equiv L^\triangleright(x,\{u,z,y\}_1)
            \\
            &=(\{\{x,y,z\}_i,u,v\}_2,\{\{x,y,z\}_i,u,v\}_3) \equiv L^\triangleright(\{x,y,z\}_i,u)
        \end{align*}
        
        \item 
        \begin{align*}
            R^\triangleleft(x,y)\dashv R^\triangleleft(z,u) &= (R_2(x,y) \circ R_1(z,u), R_1(x,y) \circ R_1(z,u))
            \\
            &\equiv ( \{\{v,x,y\}_2,z,u\}_1, \{\{v,x,y\}_1,z,u\}_1)
            \\
            &=(\{v,\{z,y,x\}_3,u\}_2, \{v,\{z,y,x\}_3,u\}_1) \equiv R^\triangleleft(\{z,y,x\}_3,u)
            \\
            &=(\{v,x,\{y,z,u\}_i\}_2, \{v,x,\{y,z,u\}_i\}_1 ) \equiv R^\triangleleft(x,\{y,z,u\}_i)
        \end{align*}

        \item 
        \begin{align*}
            R^\triangleright(x,y) \dashv R^\triangleleft(z,u) &= (R_3(x,y) \circ R_1(z,u), R_1(x,y) \circ R_1(z,u))
            \\
            &\equiv(\{\{v,x,y\}_3,z,u\}_1, \{\{v,x,y\}_1,z,u\}_1)
            \\
            &=( \{v,\{z,y,x\}_2,u\}_2, \{v,\{z,y,x\}_2,u\}_1) \equiv R^\triangleleft(\{z,y,x\}_2,u)
            \\
            &=(\{v,x,\{y,z,u\}_1\}_3, \{v,x,\{y,z,u\}_1\}_1) \equiv R^\triangleright(x,\{y,z,u\}_1)
        \end{align*}

        \item 
        \begin{align*}
            R^\triangleleft (x,y) \vdash R^\triangleleft (z,u) &= (R_1(x,y)\circ R_2(z,u), R_1(x,y) \circ R_1(z,u))
            \\
            &\equiv (\{\{v,x,y\}_1,z,u\}_2, \{\{v,x,y\}_1,z,u\}_1)
            \\
            &=( \{v,\{z,y,x\}_1,y\}_2,\{v,\{z,y,x\}_1,u\}_1) \equiv R^\triangleleft(\{z,y,x\}_1,u)
            \\
            &=(\{v,x,\{y,z,u\}_2\}_3, \{v,x,\{y,z,u\}_2\}_1) \equiv R^\triangleright(x,\{y,z,u\}_2)
        \end{align*}

        \item 
        \begin{align*}
            R^\triangleleft (x,y) \vdash R^\triangleright(z,u) &= (R_1(x,y)\circ R_3(z,u), R_1(x,y) \circ R_1(z,u))
            \\
            & \equiv (\{\{v,x,y\}_1,z,u\}_3, \{\{v,x,y\}_1,z,u\}_1)
            \\
            &=( \{v,\{z,y,x\}_i,u\}_3, \{v,\{z,y,x\}_i,u\}_1) \equiv R^\triangleright(\{z,y,x\}_i,u)
            \\
            &= (\{v,x,\{y,z,u\}_3\}_3, \{v,x,\{y,z,u\}_3\}_1) \equiv R^\triangleright(x,\{y,z,u\}_3)
        \end{align*}
    \end{enumerate}
\end{proof}

\begin{proposition}
    Let $A$ be an associative triple trisystem of the second kind. Then the map $*$ defined on $\mathfrak{L}(A,A)$ and $\mathfrak{R}(A,A)^\mathrm{op}$ by
    $$
        (L^\triangleleft(x,y))^* := L^\triangleright(y,x), \quad (L^\triangleright(x,y))^* := L^\triangleleft(y,x),
    $$
    and
    $$
        (R^\triangleleft(x,y))^* := R^\triangleright(y,x), \quad (R^\triangleright(x,y))^* := R^\triangleleft(y,x),
    $$
    for all $x,y \in A$, is an involution.
\end{proposition}
\begin{proof}
    It follows from the previous proposition.
\end{proof}
We define the following extra bilinear maps $\prec^*, \succ^* \colon A \times \mathfrak{L}(A,A) \rightarrow A$ as
$$
    x \prec^* L := L^* \succ x, \quad \text{ and } \quad x \succ^* L := L^* \prec x.
$$
and $\prec^*, \succ^* \colon \mathfrak{R}(A,A)^\mathrm{op} \times A \rightarrow A$ as
$$
    R \prec^* x := x \succ R^*, \quad \text{ and } \quad R\succ^* x := x \prec R^*.
$$
We take an isomorphic copy of $A$, denoted by $\bar{A}$. 

\begin{theorem} 
    Let $(A, \{ \, , \, , \, \}_i)$ be an associative triple trisystem of the second kind. Then $\mathfrak{L}(A,A) \oplus A \oplus \bar{A} \oplus \mathfrak{R}(A,A)^\mathrm{op}$ with the products given by the matrix multiplication
    \begin{align*}
        \begin{pmatrix} \lambda_1 & x_1 \\ \bar{y}_1 & \rho_1\end{pmatrix} \dashv \begin{pmatrix} \lambda_2 & x_2 \\ \bar{y}_2 & \rho_2\end{pmatrix} &= \begin{pmatrix} \lambda_1  \dashv \lambda_2 + L^\triangleleft (x_1,y_2)& \lambda_1 \prec x_2 + x_1 \prec \rho_2 \\ \overline{y_1\prec^* \lambda_2} + \overline{\rho_1 \prec^* y_2}& R^\triangleleft(y_1,x_2) + \rho_1 \dashv \rho_2\end{pmatrix},
    \end{align*}
    and
    \begin{align*}
        \begin{pmatrix} \lambda_1 & x_1 \\ \bar{y}_1 & \rho_1\end{pmatrix} \vdash \begin{pmatrix} \lambda_2 & x_2 \\ \bar{y}_2 & \rho_2\end{pmatrix} &= \begin{pmatrix} \lambda_1 \vdash \lambda_2 + L^\triangleright (x_1,y_2)& \lambda_1 \succ x_2 + x_1 \succ \rho_2 \\ \overline{y_1\succ^* \lambda_2} +\overline{\rho_1 \succ^* y_2}& R^\triangleright(y_1,x_2) + \rho_1 \vdash \rho_2\end{pmatrix} ,
    \end{align*}
    for every $\lambda_1, \lambda_2 \in \mathfrak{L}(A)$, $x_1,x_2 \in A$, $\bar{y}_1,\bar{y}_2\in \bar{A}$ and $\rho_1, \rho_2 \in \mathfrak{R}(A)$, is an associative dialgebra. Moreover, 
    $$
        \begin{pmatrix}\lambda & x \\ \bar{y} & \rho \end{pmatrix}^\star := \begin{pmatrix} \lambda ^* & y \\ \bar{x} & \rho^* \end{pmatrix}
    $$
    for all $\lambda \in \mathfrak{L}(A,A)$, $x \in A$, $\bar{y} \in \overline{A}$ and $\rho \in \mathfrak{R}(A,A)$, is an involution of $\mathfrak{L}(A,A) \oplus A \oplus \bar{A} \oplus \mathfrak{R}(A,A)^\mathrm{op}$. And, $\{x,y,z\}_1 = x \dashv (y^\star \dashv z)$, $\{x,y,z\}_2 = x \vdash (y^\star \dashv z)$ and $\{x,y,z\}_3 = x \vdash (y^\star \vdash z)$,  for all $x,y,z \in A$.
\end{theorem}
\begin{proof}
    The proof is straightforward. Nevertheless, we give here the proof with details. First, we prove the axioms $(X \dashv Y) \dashv Z = X \dashv (Y \dashv Z) = X \dashv ( Y \vdash Z)$ for all $X, Y,Z \in \mathfrak{L}(A) \oplus A \oplus \bar{A} \oplus \mathfrak{R}(A)$. To do so, we compute the three different matrices and we compare their entries.
    \begin{align*}
        &\left(\begin{pmatrix} \lambda_1 & x_1 \\ \bar{y}_1 & \rho_1\end{pmatrix} \dashv \begin{pmatrix} \lambda_2 & x_2 \\ \bar{y}_2 & \rho_2\end{pmatrix}\right) \dashv \begin{pmatrix} \lambda_3 & x_3 \\ \bar{y}_3 & \rho_3\end{pmatrix}
        \\
        &=\begin{pmatrix} \lambda_1  \dashv \lambda_2 + L^\triangleleft (x_1,y_2)& \lambda_1 \prec x_2 + x_1 \prec \rho_2 \\ \overline{ \lambda_2^* \succ y_1} + \overline{y_2 \succ \rho_1^* }& R^\triangleleft(y_1,x_2) + \rho_1 \dashv \rho_2\end{pmatrix} \dashv \begin{pmatrix} \lambda_3 & x_3 \\ \bar{y}_3 & \rho_3\end{pmatrix} = \begin{pmatrix} \alpha_{11} & \alpha_{12}\\ \alpha_{21} & \alpha_{22}\end{pmatrix}
    \end{align*}
    where
    \begin{align*}
        \alpha_{11}&=(\lambda_1\dashv \lambda_2) \dashv \lambda_3 + L^\triangleleft(x_1,y_2)\dashv \lambda_3 + L^\triangleleft(\lambda_1 \prec x_2,y_3) + L^\triangleleft(x_1 \prec \rho_2,y_3)
        \\
        \alpha_{12}&=(\lambda_1 \dashv \lambda_2 )\prec x_3 + L^\triangleleft(x_1,y_2)\prec x_3 + (\lambda_1 \prec x_2)\prec \rho_3 + (x_1 \prec \rho_2) \prec \rho_3
        \\
        \alpha_{21}&= \overline{ \lambda_3^* \succ(\lambda_2^* \succ y_1)} + \overline{ \lambda_3^*\succ (y_2\succ \rho_1^*)} + \overline{ y_3 \succ R^\triangleright(x_2,y_1)}+ \overline{ y_3\succ (\rho_2^* \vdash \rho_1^*)}
        \\
        \alpha_{22}&=R^\triangleleft( \lambda_2^* \succ y_1,x_3) + R^\triangleleft(y_2\succ \rho_1^*,x_3) + R^\triangleleft(y_1,x_2)\dashv \rho_3 + (\rho_1 \dashv \rho_2) \dashv \rho_3
    \end{align*}
    
    \begin{align*}
        &\begin{pmatrix} \lambda_1 & x_1 \\ \bar{y}_1 & \rho_1\end{pmatrix} \dashv \left(\begin{pmatrix} \lambda_2 & x_2 \\ \bar{y}_2 & \rho_2\end{pmatrix} \dashv \begin{pmatrix} \lambda_3 & x_3 \\ \bar{y}_3 & \rho_3\end{pmatrix}\right)
        \\
        &=\begin{pmatrix} \lambda_1 & x_1 \\ \bar{y}_1 & \rho_1\end{pmatrix} \dashv \begin{pmatrix} \lambda_2  \dashv \lambda_3 + L^\triangleleft (x_2,y_3)& \lambda_2 \prec x_3 + x_2 \prec \rho_3 \\ \overline{\lambda_3^* \succ y_2} + \overline{y_3 \succ \rho_2^*}& R^\triangleleft(y_2,x_3) + \rho_2 \dashv \rho_3\end{pmatrix} = \begin{pmatrix} \beta_{11} & \beta_{12}\\ \beta_{21} & \beta_{22}\end{pmatrix},
    \end{align*}
    where
    \begin{align*}
        \beta_{11}&=\lambda_1\dashv (\lambda_2 \dashv \lambda_3) + \lambda_1 \dashv L^\triangleleft(x_2,y_3)+ L^\triangleleft(x_1,\lambda_3^*\succ y_2) + L^\triangleleft(x_1,y_3 \succ \rho_2^*)
        \\
        \beta_{12}&=\lambda_1 \prec (\lambda_2 \prec x_3 )+ \lambda_1\prec (x_2 \prec \rho_3) + x_1 \prec R^\triangleleft(y_2,x_3) + x_1 \prec (\rho_2 \dashv \rho_3)
        \\
        \beta_{21}&= \overline{(\lambda_3^* \vdash \lambda_2^*) \succ y_1} + \overline{L^\triangleright(y_3,x_2)\succ y_1}+ \overline{(\lambda_3^* \succ y_2 ) \succ  \rho_1^*} + \overline{(y_3 \succ \rho_2^*) \succ \rho_1^*}
        \\
        \beta_{22}&= R^\triangleleft(y_1,\lambda_2\prec x_3) + R^\triangleleft(y_1,x_2 \prec \rho_3) + \rho_1 \dashv R^\triangleleft(y_2,x_3) + \rho_1 \dashv (\rho_2 \dashv \rho_3)
    \end{align*}
    
    \begin{align*}
        &\begin{pmatrix} \lambda_1 & x_1 \\ \bar{y}_1 & \rho_1\end{pmatrix} \dashv \left(\begin{pmatrix} \lambda_2 & x_2 \\ \bar{y}_2 & \rho_2\end{pmatrix} \vdash \begin{pmatrix} \lambda_3 & x_3 \\ \bar{y}_3 & \rho_3\end{pmatrix}\right)
        \\
        &=\begin{pmatrix} \lambda_1 & x_1 \\ \bar{y}_1 & \rho_1\end{pmatrix} \dashv \begin{pmatrix} \lambda_2  \vdash \lambda_3 + L^\triangleright (x_2,y_3)& \lambda_2 \succ x_3 + x_2 \succ \rho_3 \\ \overline{\lambda_3^* \prec y_2} + \overline{y_3 \prec \rho_2^*}& R^\triangleright(y_2,x_3) + \rho_2 \vdash \rho_3\end{pmatrix}  = \begin{pmatrix} \gamma_{11} & \gamma_{12}\\ \gamma_{21} & \gamma_{22}\end{pmatrix},
    \end{align*}
    where
    \begin{align*}
        \gamma_{11}&=\lambda_1 \dashv (\lambda_2 \vdash \lambda_3) + \lambda_1 \dashv L^\triangleright(x_2,y_3) + L^\triangleleft(x_1,\lambda_3^*\prec y_2) + L^\triangleleft(x_1,y_3 \prec \rho_2^*)
        \\
        \gamma_{12}&=\lambda_1 \prec (\lambda_2 \succ x_3) + \lambda_1 \prec (x_2 \succ \rho_3) + x_1 \prec R^\triangleright(y_2,x_3) + x_1 \prec (\rho_2 \vdash \rho_3)
        \\
        \gamma_{21}&= \overline{(\lambda_3^* \dashv \lambda_2^*)\succ y_1}+ \overline{L^\triangleleft(y_3,x_2) \succ y_1} + \overline{(\lambda_3^* \prec y_2) \succ \rho_1^*} + \overline{(y_3 \prec \rho_2^*) \succ \rho_1^*}
        \\
        \gamma_{22}&= R^\triangleleft(y_1,\lambda_2 \succ x_3) + R^\triangleleft(y_1,x_2 \succ \rho_3) + \rho_1 \dashv R^\triangleright(y_2,x_3) + \rho_1 \dashv ( \rho_2 \vdash \rho_3)
    \end{align*}
    In each entry of the matrices there is a linear combination of four different terms. To prove that these three matrices are equal, entry by entry, we will see that term by term are equal.
    \begin{enumerate}
        \item Let us prove that $\alpha_{11} = \beta_{11} = \gamma_{11}$:
        \begin{enumerate}
            \item $(\lambda_1 \dashv \lambda_2) \dashv \lambda_3 = \lambda_1 \dashv (\lambda_2 \dashv \lambda_3) = \lambda_1 \dashv (\lambda_2 \vdash \lambda_3)$. Because $\mathfrak{L}(A)$ is an associative dialgebra.

            \item $L^\triangleleft(x_1,y_2) \dashv \lambda_3=L^\triangleleft(x_1,\lambda_3^* \succ y_2) = L^\triangleleft(x_1,\lambda_3^* \prec y_2)$. If $\lambda_3 = L^\mu(z,u)$ with $z,u \in A$ and $\mu \in \{\triangleleft, \triangleright\}$, then, by Remark \ref{rem:asssecond} and Proposition \ref{prop:asssecondkind}(1), $L^\triangleleft (x_1,y_2) \dashv \lambda_3 = L^\triangleleft(x_1,\{u,z,y_2\}_i)$ for any $i \in \{1,2,3\}$. Hence, we can conclude that $L^\triangleleft (x_1,y_2) \dashv \lambda_3= L^\triangleleft(x_1,\lambda_3^*\succ y_2)= L^\triangleleft(x_1,\lambda_3^*\prec y_2)$ for all $\lambda_3 \in \mathfrak{L}(A)$.

            \item $L^\triangleleft(\lambda_1 \prec x_2,y_3) = \lambda_1 \dashv L^\triangleleft(x_2,y_3) = \lambda_1 \dashv L^\triangleright(x_2,y_3)$. By Remark \ref{rem:asssecond}, $\lambda_1 \dashv L^\triangleleft(x_2,y_3) = \lambda_1 \dashv L^\triangleright(x_2,y_3)$. If $\lambda_1 = L^\triangleleft(z,u)$ with $z,u \in A$, then, by Proposition \ref{prop:asssecondkind}(1), $\lambda_1 \dashv L^\triangleleft(x_2,y_3)  = L^\triangleleft(\{z,u,x_2\}_1,y_3) = L^\triangleleft ( \lambda_1\prec x_2,y_3)$. Otherwise, if $\lambda_1 = L^\triangleright(z,u)$ then, by Proposition \ref{prop:asssecondkind}(2), $\lambda_1 \dashv L^\triangleleft(x_2,y_3)  =L^\triangleleft(\{z,u,x_2\}_2,y_3)= L^\triangleleft ( \lambda_1\prec x_2,y_3)$. Therefore, for all $\lambda_1 \in \mathfrak{L}(A)$ we can assure that $\lambda_1 \dashv L^\triangleleft(x_2,y_3) = L^\triangleleft ( \lambda_1\prec x_2,y_3)$.

            \item $L^\triangleleft(x_1 \prec \rho_2, y_3) = L^\triangleleft (x_1, y_3 \succ \rho_2^*) = L^\triangleleft(x_1,y_3 \prec \rho_2^*)$. For both cases $\rho_2 = R^\triangleleft(z,u)$ or $\rho_2 = R^\triangleright(z,u)$ with $z,u \in A$, we have that $x_1 \prec \rho_2 = \{x_1,z,u\}_1$. Thus, by Proposition \ref{prop:asssecondkind}(1), for any $i \in \{1,2,3\}$, we have $L^\triangleleft(x_1 \prec \rho_2, y_3) = L^\triangleleft(x_1,\{y_3,u,z\}_i) = L^\triangleleft(x_1, y_3 \prec \rho_2^*) = L^\triangleleft(x_1, y_3 \succ \rho_2^*)$.
        \end{enumerate}

        \item Let us prove that $\alpha_{12} = \beta_{12} = \gamma_{12}$:
        \begin{enumerate}
            \item $(\lambda_1 \dashv \lambda_2) \prec x_3 = \lambda_1 \prec (\lambda_2 \prec x_3) = \lambda_1 \prec (\lambda_2 \succ x_3)$. By Lemma \ref{lemma:diendomorphism}, $(\lambda_1 \dashv \lambda_2) \prec x_3 =  \lambda_1 \prec (\lambda_2 \succ x_3)$. Let us pick $\lambda_1 = L^\mu (y,z)$ and $\lambda_2 = L^\nu(u,v)$, with $y,z,u,v \in A$ and $\mu, \nu \in \{\triangleleft, \triangleright\}$. If $\mu = \nu = \triangleleft$, then 
            $$
                \lambda_1\prec (\lambda_2 \prec x_3) = \{y,z,\{u,v,x_3\}_1\}_1 = \{y,z,\{u,v,x_3\}_3\}_1 = \lambda_1\prec (\lambda_2 \succ x_3).
            $$
            If $\mu = \triangleleft$ and $\nu = \triangleright$, then
            $$
                \lambda_1\prec (\lambda_2 \prec x_3) = \{y,z,\{u,v,x_3\}_2\}_1 = \{y,z,\{u,v,x_3\}_3\}_1 = \lambda_1\prec (\lambda_2 \succ x_3).
            $$
            If $\mu = \triangleright$ and $\nu = \triangleleft$, then
            $$
                \lambda_1\prec (\lambda_2 \prec x_3) = \{y,z,\{u,v,x_3\}_1\}_2 = \{y,z,\{u,v,x_3\}_3\}_2 = \lambda_1\prec (\lambda_2 \succ x_3).
            $$
            If $\mu = \nu = \triangleright$, then 
            $$
                \lambda_1\prec (\lambda_2 \prec x_3) = \{y,z,\{u,v,x_3\}_2\}_2 = \{y,z,\{u,v,x_3\}_3\}_2 = \lambda_1\prec (\lambda_2 \succ x_3).
            $$
            
            \item $L^\triangleleft(x_1,y_2)\prec x_3 = x_1 \prec R^\triangleleft(y_2,x_3) = x_1 \prec R^\triangleright(y_2,x_3)$. By the definition of the di-endomorphisms, all of them are equal to $\{x_1,y_2,x_3\}_1$.

            \item $(\lambda_1 \prec x_2) \prec \rho_3 = \lambda_1 \prec (x_2 \prec \rho_3) = \lambda_1 \prec (x_2 \succ \rho_3)$. Let us pick $\lambda_1 = L^\mu (y,z)$ and $\rho_3 = R^\nu(u,v)$ with $y,z,u,v \in A$ and $\mu, \nu \in \{\triangleleft, \triangleright\}$. If $\mu = \nu = \triangleleft$, then
            \begin{align*}
                (\lambda_1 \prec x_2) \prec \rho_3 = \{\{y,z,x_2\}_1,u,v\}_1 &=\{y,z,\{x_2,u,v\}_1\}_1 = \lambda_1 \prec (x_2 \prec \rho_3)
                \\
                &= \{y,z,\{x_2,u,v\}_2\}_1= \lambda_1 \prec (x_2 \succ \rho_3).
            \end{align*}
            If $\mu = \triangleleft$ and $\nu = \triangleright$, then
            \begin{align*}
                (\lambda_1 \prec x_2) \prec \rho_3 = \{\{y,z,x_2\}_1,u,v\}_1 &=\{y,z,\{x_2,u,v\}_1\}_1 = \lambda_1 \prec (x_2 \prec \rho_3)
                \\
                &= \{y,z,\{x_2,u,v\}_3\}_1= \lambda_1 \prec (x_2 \succ \rho_3).
            \end{align*}
            If $\mu = \triangleright$ and $\nu = \triangleleft$, then
            \begin{align*}
                (\lambda_1 \prec x_2) \prec \rho_3 = \{\{y,z,x_2\}_2,u,v\}_1 &=\{y,z,\{x_2,u,v\}_1\}_2 = \lambda_1 \prec (x_2 \prec \rho_3)
                \\
                &= \{y,z,\{x_2,u,v\}_2\}_2= \lambda_1 \prec (x_2 \succ \rho_3).
            \end{align*}
            If $\mu = \nu = \triangleright$, then
            \begin{align*}
                (\lambda_1 \prec x_2) \prec \rho_3 = \{\{y,z,x_2\}_2,u,v\}_1 &=\{y,z,\{x_2,u,v\}_1\}_2 = \lambda_1 \prec (x_2 \prec \rho_3)
                \\
                &= \{y,z,\{x_2,u,v\}_3\}_2= \lambda_1 \prec (x_2 \succ \rho_3).
            \end{align*}

            \item $(x_1 \prec \rho_2 ) \prec \rho_3 = x_1 \prec (\rho_2 \dashv \rho_3) = x_1 \prec (\rho_2 \vdash \rho_3)$. Let us pick $\rho_2 = R^\mu(y,z)$ and $\rho_3 = R^\nu(u,v)$ with $y,z,u,v \in A$ and $\mu, \nu \in \{ \triangleleft, \triangleright\}$. Therefore
            \begin{align*}
                (x_1 \prec \rho_2 ) \prec \rho_3 = \{\{x_1,y,z\}_1,u,v\}_1 =x_1 \prec (\rho_2 \dashv \rho_3) = x_1 \prec (\rho_2 \vdash \rho_3).
            \end{align*}
        \end{enumerate}

        \item Let us prove that $\alpha_{21} = \beta_{21} = \gamma_{21}$:
        \begin{enumerate}
            \item $\lambda_3^* \succ (\lambda_2^* \succ y_1) = (\lambda_3^* \vdash \lambda_2^*) \succ y_1 = (\lambda_3^* \dashv \lambda_2^*) \succ y_1$. By Lemma \ref{lemma:diendomorphism}(1).
            
            \item $\lambda_3^* \succ (y_2 \succ \rho_1^*) = (\lambda_3^* \succ y_2) \succ \rho_1^* = (\lambda_3^* \prec y_2) \succ \rho_1^*$. Let us pick $\lambda_3^* = L^\mu (x,z)$ and $\rho_1^* = R^\nu(u,v)$, with $x,z,u,v \in A$ and $\mu, \nu \in\{ \triangleleft, \triangleright\}$. If $\mu = \nu = \triangleleft$, then
            \begin{align*}
                \lambda_3^* \succ (y_2 \succ \rho_1^*) = \{x,z,\{y_2,u,v\}_2 \}_3 &=\{\{x,z,y_2 \}_3 ,u,v\}_2 = (\lambda_3^* \succ y_2) \succ \rho_1^*
                \\
                &=\{\{x,z,y_2 \}_1 ,u,v\}_2 =(\lambda_3^* \prec y_2) \succ \rho_1^*.
            \end{align*}
            
            If $\mu = \triangleleft$ and $\nu = \triangleright$, then
            \begin{align*}
                \lambda_3^* \succ (y_2 \succ \rho_1^*) = \{x,z,\{y_2,u,v\}_3 \}_3 &=\{\{x,z,y_2 \}_3 ,u,v\}_3 = (\lambda_3^* \succ y_2) \succ \rho_1^*
                \\
                &=\{\{x,z,y_2 \}_1 ,u,v\}_3 =(\lambda_3^* \prec y_2) \succ \rho_1^*.
            \end{align*}

            If $\mu = \triangleright$ and $\nu = \triangleleft$, then
            \begin{align*}
                \lambda_3^* \succ (y_2 \succ \rho_1^*) = \{x,z,\{y_2,u,v\}_2 \}_3 &=\{\{x,z,y_2 \}_3 ,u,v\}_2 = (\lambda_3^* \succ y_2) \succ \rho_1^*
                \\
                &=\{\{x,z,y_2 \}_2 ,u,v\}_2 =(\lambda_3^* \prec y_2) \succ \rho_1^*.
            \end{align*}

            If $\mu = \nu = \triangleright$, then
            \begin{align*}
                \lambda_3^* \succ (y_2 \succ \rho_1^*) = \{x,z,\{y_2,u,v\}_3 \}_3 &=\{\{x,z,y_2 \}_3 ,u,v\}_3 = (\lambda_3^* \succ y_2) \succ \rho_1^*
                \\
                &=\{\{x,z,y_2 \}_2 ,u,v\}_3 =(\lambda_3^* \prec y_2) \succ \rho_1^*.
            \end{align*}
            
            \item $y_3 \succ R^\triangleright(x_2,y_1) = L^\triangleright(y_3,x_2) \succ y_1 = L^\triangleleft(y_3,x_2) \succ y_1$. By the definition of the di-endomorphisms, all of them are equal to $\{y_3,x_2,y_1\}_3$.
            
            \item $ y_3 \succ ( \rho_2^* \vdash \rho_1^*) = (y_3 \succ \rho_2^*) \succ \rho_1^* = (y_3 \prec \rho_2^*) \succ \rho_1^*$. By Lemma \ref{lemma:diendomorphism}, $y_3 \succ (\rho_2^* \vdash \rho_1^*) = (y_3 \prec \rho_2^*) \succ \rho_1^*$. Let us pick $\rho_2^* = R^\mu(x,z)$ and $\rho_1^* = R^\nu (u,v)$, with $x,z,u,v \in A$ and $\mu, \nu \in \{ \triangleleft, \triangleright\}$. If $\mu = \nu = \triangleleft$, then
            $$
                (y_3 \succ \rho_2^*) \succ \rho_1^* = \{\{y_3,x,z\}_2,u,v\}_2 = \{\{y_3,x,z\}_1,u,v\}_2 = (y_3 \prec \rho_2^*) \succ \rho_1^*.
            $$
            If $\mu = \triangleleft$ and $\nu = \triangleright$, then
            $$
                (y_3 \succ \rho_2^*) \succ \rho_1^* = \{\{y_3,x,z\}_2,u,v\}_3 = \{\{y_3,x,z\}_1,u,v\}_3 = (y_3 \prec \rho_2^*) \succ \rho_1^*.
            $$
            If $\mu = \triangleright$ and $\nu = \triangleleft$, then
            $$
                (y_3 \succ \rho_2^*) \succ \rho_1^* = \{\{y_3,x,z\}_3,u,v\}_2 = \{\{y_3,x,z\}_1,u,v\}_2 = (y_3 \prec \rho_2^*) \succ \rho_1^*.
            $$
            If $\mu = \nu = \triangleright$, then
            $$
                (y_3 \succ \rho_2^*) \succ \rho_1^* = \{\{y_3,x,z\}_3,u,v\}_3 = \{\{y_3,x,z\}_1,u,v\}_3 = (y_3 \prec \rho_2^*) \succ \rho_1^*.
            $$
        \end{enumerate}

        \item Let us prove that $\alpha_{22} = \beta_{22} = \gamma_{22}$:
        \begin{enumerate}
            \item $R^\triangleleft(\lambda_2^* \succ y_1,x_3) = R^\triangleleft(y_1, \lambda_2 \prec x_3) = R^\triangleleft(y_1, \lambda_2 \succ x_3)$. For both cases $\lambda_2 = L^\triangleleft(z,u)$ or $\lambda_2 = L^\triangleright(z,u)$ with $z,u \in A$, we have $\lambda_2^* \succ y_1 = \{u,z,y_1\}_3$. Thus, by Proposition \ref{prop:asssecondkind}(5), for any $i \in \{1,2,3\}$, we have $R^\triangleleft(\lambda_2^* \succ y_1, x_3) = R^\triangleleft(y_1,\{z,u,x_3\}_i) = R^\triangleleft(y_1, \lambda_2 \prec x_3) = R^\triangleleft(y_1, \lambda_2 \succ x_3)$.
            
            \item $R^\triangleleft(y_2 \succ \rho_1^*,x_3) = \rho_1 \dashv R^\triangleleft(y_2,x_3) = \rho_1 \dashv R^\triangleright(y_2,x_3)$. By Remark $\ref{rem:asssecond}$, $\rho_1 \dashv R^\triangleleft(y_2,x_3) = \rho_1 \dashv R^\triangleright(y_2,x_3)$. If $\rho_1 = R^\triangleleft(z,u)$ with $z,u \in A$, then, by Proposition \ref{prop:asssecondkind}(5), $\rho_1 \dashv R^\triangleleft(y_2,x_3) = R^\triangleleft(\{y_2,u,z\}_3,x_3) = R^\triangleleft(y_2  \succ \rho_1^*,x_3)$. Otherwise, if $\rho_1 = R^\triangleright(z,u)$, then by Proposition \ref{prop:asssecondkind}(6), $\rho_1 \dashv R^\triangleleft(y_2,x_3) = R^\triangleleft(\{y_2,u,z\}_2,x_3) = R^\triangleleft(y_2  \succ \rho_1^*,x_3)$. Therefore, for any $\rho_1 \in \mathfrak{R}(A)$ we can assure that $\rho_1 \dashv R^\triangleleft(y_2,x_3) = R^\triangleleft (y_2 \succ \rho_1^*,x_3)$.
            
            \item $R^\triangleleft(y_1,x_2) \dashv \rho_3 = R^\triangleleft(y_1, x_2 \prec \rho_3) = R^\triangleleft(y_1, x_2\succ \rho_3)$. If $\rho_3 = R^\mu(z,u)$ with $z, u \in A$ and $\mu \in \{ \triangleleft , \triangleright \}$, then, by Remark \ref{rem:asssecond} and Proposition \ref{prop:asssecondkind}(5), $R^\triangleleft(y_1,x_2) \dashv \rho_3 = R^\triangleleft(y_1,\{u,z,x_2\}_i)$ for any $i \in \{1,2,3\}$. Hence, we can conclude that $R^\triangleleft(y_1,x_2) \dashv \rho_3 = R^\triangleleft(y_1, x_2 \prec \rho_3) = R^\triangleleft(y_1, x_2\succ \rho_3)$ for all $\rho_3 \in \mathfrak{R}(A)$.
            
            \item $(\rho_1 \dashv \rho_2 ) \dashv \rho_3 = \rho_1 \dashv (\rho_2 \dashv \rho_3) = \rho_1 \dashv (\rho_2 \vdash \rho_3)$. Because $\mathfrak{R}(A)$ is an associative dialgebra.
        \end{enumerate}
    \end{enumerate}

    The next two axioms we prove are $(X \dashv Y ) \vdash Z = (X \vdash Y) \vdash Z = X \vdash (Y \vdash Z)$.
    \begin{align*}
        &\left(\begin{pmatrix} \lambda_1 & x_1 \\ \bar{y}_1 & \rho_1\end{pmatrix} \dashv \begin{pmatrix} \lambda_2 & x_2 \\ \bar{y}_2 & \rho_2\end{pmatrix}\right) \vdash\begin{pmatrix} \lambda_3 & x_3 \\ \bar{y}_3 & \rho_3\end{pmatrix}
        \\
        &=\begin{pmatrix} \lambda_1  \dashv \lambda_2 + L^\triangleleft (x_1,y_2)& \lambda_1 \prec x_2 + x_1 \prec \rho_2 \\ \overline{ \lambda_2^* \succ y_1} + \overline{y_2 \succ \rho_1^* }& R^\triangleleft(y_1,x_2) + \rho_1 \dashv \rho_2\end{pmatrix} \vdash \begin{pmatrix} \lambda_3 & x_3 \\ \bar{y}_3 & \rho_3\end{pmatrix} = \begin{pmatrix}
            \alpha_{11} & \alpha_{12} \\ \alpha_{21} & \alpha_{22}
        \end{pmatrix}
    \end{align*}
    where
    \begin{align*}
        \alpha_{11}&= (\lambda_1 \dashv \lambda_2) \vdash \lambda_3 + L^\triangleleft(x_1,y_2) \vdash \lambda_3 + L^\triangleright(\lambda_1 \prec x_2,y_3) + L^\triangleright(x_1 \prec \rho_2,y_3)
        \\
        \alpha_{12}&= (\lambda_1 \dashv \lambda_2) \succ x_3 + L^\triangleleft(x_1,y_2) \succ x_3 + (\lambda_1 \prec x_2) \succ \rho_3+ (x_1 \prec \rho_2 ) \succ \rho_3
        \\
        \alpha_{21}&= \overline{\lambda_3^* \prec(\lambda_2^* \succ y_1)} + \overline{\lambda_3^* \prec (y_2 \succ \rho_1^*)}+ \overline{y_3 \prec R^\triangleright(x_2,y_1)}+\overline{y_3 \prec (\rho_2^* \vdash \rho_1^*)}
        \\
        \alpha_{22}&=  R^\triangleright(\lambda_2^*\succ y_1,x_3) +R^\triangleright(y_2 \succ \rho_1^*,x_3) +R^\triangleleft(y_1,x_2) \vdash \rho_3 + (\rho_1 \dashv \rho_2) \vdash \rho_3
    \end{align*}
    
    \begin{align*}
        &\left( \begin{pmatrix} \lambda_1 & x_1 \\ \bar{y}_1 & \rho_1\end{pmatrix} \vdash \begin{pmatrix} \lambda_2 & x_2 \\ \bar{y}_2 & \rho_2\end{pmatrix}\right) \vdash \begin{pmatrix} \lambda_3 & x_3 \\ \bar{y}_3 & \rho_3\end{pmatrix} 
        \\
        &= \begin{pmatrix}\lambda_1 \vdash \lambda_2 + L^\triangleright(x_1,y_2) & \lambda_1 \succ x_2 + x_1 \succ \rho_2\\ \overline{\lambda_2^*\prec y_1} + \overline{y_2 \prec \rho_1^*} & R^\triangleright(y_1,x_2) + \rho_1 \vdash \rho_2 \end{pmatrix}\vdash \begin{pmatrix} \lambda_3 & x_3 \\ \bar{y}_3 & \rho_3\end{pmatrix} = \begin{pmatrix} \beta_{11} & \beta_{12} \\ \beta_{21} & \beta_{22}\end{pmatrix}
    \end{align*}
    where
    \begin{align*}
        \beta_{11}&= (\lambda_1 \vdash \lambda_2)\vdash \lambda_3 + L^\triangleright(x_1,y_2)\vdash \lambda_3 + L^\triangleright(\lambda_1\succ x_2,y_3) + L^\triangleright(x_1 \succ \rho_2, y_3)
        \\
        \beta_{12}&= (\lambda_1 \vdash \lambda_2) \succ x_3 + L^\triangleright(x_1,y_2) \succ x_3 + (\lambda_1 \succ x_2) \succ \rho_3 + (x_1 \succ \rho_2) \succ \rho_3
        \\
        \beta_{21}&=  \overline{\lambda_3^* \prec (\lambda_2^* \prec y_1)} + \overline{\lambda_3^* \prec (y_2 \prec \rho_1^*)}+\overline{y_3 \prec R^\triangleleft(x_2,y_1)} + \overline{y_3 \prec(\rho_2^* \dashv \rho_1^*)}
        \\
        \beta_{22}&= R^\triangleright(\lambda_2^* \prec y_1,x_3) + R^\triangleright(y_2 \prec \rho_1^*,x_3) + R^\triangleright(y_1,x_2) \vdash \rho_3 + (\rho_1 \vdash \rho_2) \vdash \rho_3
    \end{align*}
    
    \begin{align*}
        &\begin{pmatrix} \lambda_1 & x_1 \\ \bar{y}_1 & \rho_1\end{pmatrix} \vdash \left(\begin{pmatrix} \lambda_2 & x_2 \\ \bar{y}_2 & \rho_2\end{pmatrix} \vdash \begin{pmatrix} \lambda_3 & x_3 \\ \bar{y}_3 & \rho_3\end{pmatrix}\right)
        \\
        &=\begin{pmatrix} \lambda_1 & x_1 \\ \bar{y}_1 & \rho_1\end{pmatrix} \vdash \begin{pmatrix} \lambda_2  \vdash \lambda_3 + L^\triangleright (x_2,y_3)& \lambda_2 \succ x_3 + x_2 \succ \rho_3 \\ \overline{\lambda_3^* \prec y_2} + \overline{y_3 \prec \rho_2^*}& R^\triangleright(y_2,x_3) + \rho_2 \vdash \rho_3\end{pmatrix} = \begin{pmatrix} \gamma_{11} & \gamma_{12} \\ \gamma_{21} & \gamma_{22}\end{pmatrix}
    \end{align*}
    where
    \begin{align*}
        \gamma_{11}&= \lambda_1 \vdash (\lambda_2 \vdash \lambda_3)+ \lambda_1 \vdash L^\triangleright(x_2,y_3) + L^\triangleright(x_1,\lambda_3^*\prec y_2) + L^\triangleright(x_1,y_3 \prec \rho_2^*)
        \\
        \gamma_{12}&= \lambda_1 \succ(\lambda_2 \succ x_3) + \lambda_1 \succ (x_2 \succ \rho_3) + x_1 \succ R^\triangleright(y_2,x_3) + x_1 \succ (\rho_2 \vdash \rho_3)
        \\
        \gamma_{21}&=  \overline{(\lambda_3^* \dashv \lambda_2^*)\prec y_1} + \overline{L^\triangleleft(y_3,x_2) \prec y_1} + \overline{(\lambda_3^* \prec y_2)\prec \rho_1^*} + \overline{(y_3 \prec \rho_2^*) \prec \rho_1^*}
        \\
        \gamma_{22}&= R^\triangleright(y_1,\lambda_2 \succ x_3) + R^\triangleright(y_1,x_2 \succ \rho_3) + \rho_1\vdash R^\triangleright(y_2,x_3) + \rho_1 \vdash (\rho_2 \vdash \rho_3)
    \end{align*}

    \begin{enumerate}
        \item Let us prove that $\alpha_{11} = \beta_{11} = \gamma_{11}$:
        \begin{enumerate}
            \item $(\lambda_1 \dashv \lambda_2) \vdash \lambda_3 = ( \lambda_1 \vdash \lambda_2) \vdash \lambda_3 = \lambda_1 \vdash (\lambda_2 \vdash \lambda_3)$. Because $\mathfrak{L}(A)$ is an associative dialgebra.

            \item $L^\triangleleft(x_1,y_2) \vdash \lambda_3= L^\triangleright(x_1,y_2) \vdash \lambda_3 = L^\triangleright(x_1, \lambda_3^* \prec y_2)$. By Remark \ref{rem:asssecond}, $L^\triangleleft(x_1,y_2) \vdash \lambda_3= L^\triangleright(x_1,y_2) \vdash \lambda_3$. If $\lambda_3 = L^\triangleleft(z,u)$ with $z,u \in A$, then, by Proposition \ref{prop:asssecondkind}(3), $ L^\triangleleft(x_1,y_2) \vdash \lambda_3 = L^\triangleleft(x_1,\{u,z,y_2\}_2) = L^\triangleleft ( x_1, \lambda_3^*\prec y_2)$. Otherwise, if $\lambda_3 = L^\triangleright(z,u)$ then, by Proposition \ref{prop:asssecondkind}(4), $ L^\triangleleft(x_1,y_2) \vdash \lambda_3 = L^\triangleleft(x_1,\{u,z,y_2\}_1) = L^\triangleleft ( x_1, \lambda_3^*\prec y_2)$. Therefore, for all $\lambda_1 \in \mathfrak{L}(A)$ we can assure that $L^\triangleleft(x_1,y_2) \vdash \lambda_3=  L^\triangleright(x_1, \lambda_3^* \prec y_2)$.

            \item $L^\triangleright(\lambda_1 \prec x_2, y_3) = L^\triangleright(\lambda_1 \succ x_2, y_3) = \lambda_1 \vdash L^\triangleright (x_2,y_3)$. If $\lambda_1 = L^\mu(z,u)$ with $z,u \in A$ and $\mu \in \{\triangleleft, \triangleright\}$, then, by Remark \ref{rem:asssecond} and Proposition \ref{prop:asssecondkind}(4), $\lambda_1 \vdash L^\triangleright (x_2,y_3) = L^\triangleright(\{z,u,x_2\}_i,y_3)$ for any $i \in \{1,2,3\}$. Hence, we can conclude that $L^\triangleright(\lambda_1 \prec x_2, y_3) = L^\triangleright(\lambda_1 \succ x_2, y_3) = \lambda_1 \vdash L^\triangleright (x_2,y_3)$ for all $\lambda_1 \in \mathfrak{L}(A)$.

            \item $L^\triangleright(x_1 \prec \rho_2, y_3) = L^\triangleright(x_1 \succ \rho_2, y_3) = L^\triangleright(x_1,y_3 \prec \rho_2^*)$. For both cases $\rho_2 = R^\triangleleft(z,u)$ or $\rho_2 = R^\triangleright(z,u)$ with $z,u \in A$, we have that $y_3 \prec \rho_2^* = \{y_3,u,z\}_1$. Thus, by Proposition \ref{prop:asssecondkind}(4), for any $i \in \{1,2,3\}$, we have $L^\triangleright(x_1,y_3 \prec \rho_2^*) = L^\triangleright(\{x_1,z,u \}_i,y_3) =L^\triangleright(x_1 \prec \rho_2, y_3) = L^\triangleright(x_1 \succ \rho_2, y_3)$.
        \end{enumerate}
        
        \item Let us prove that $\alpha_{12} = \beta_{12} = \gamma_{12}$:
        \begin{enumerate}
            \item $(\lambda_1 \dashv \lambda_2) \succ x_3 = (\lambda_1 \vdash \lambda_2) \succ x_3 = \lambda_1 \succ ( \lambda_2 \succ x_3)$. By Lemma \ref{lemma:diendomorphism}(1).
        
            \item $L^\triangleleft(x_1,y_2) \succ x_3 = L^\triangleright(x_1,y_2) \succ x_3 = x_1 \succ R^\triangleright(y_2,x_3)$. By the definition of the di-endomorphisms, all of them are equal to $\{x_1,y_2,x_3\}_3$.

            \item $(\lambda_1 \prec x_2) \succ \rho_3 = (\lambda_1 \succ x_2) \succ \rho_3 = \lambda_1 \succ (x_2 \succ \rho_3)$. Is the previous case (3b).

            \item $(x_1 \prec \rho_2 ) \succ \rho_3 = (x_1 \succ \rho_2) \succ \rho_3 = x_1 \succ ( \rho_2 \vdash \rho_3)$. Is the previous case (3d).
        \end{enumerate}
        \item Let us prove that $\alpha_{21} = \beta_{21} = \gamma_{21}$:
        \begin{enumerate}
            \item $\lambda_3^* \prec (\lambda_2^* \succ y_1) = \lambda_3^*\prec ( \lambda_2^* \prec y_1) = ( \lambda_3^* \dashv \lambda_2^*) \prec y_1$. Is the previous case (2a).
            
            \item $\lambda_3^* \prec (y_2 \succ \rho_1^*) = \lambda_3^* \prec (y_2 \prec \rho_1^*) = (\lambda_3^* \prec y_2) \prec \rho_1^*$. Is the previous case (2c).
            
            \item $y_3 \prec R^\triangleright(x_2,y_1) = y_3 \prec R^\triangleleft(x_2,y_1) = L^\triangleleft(y_3,x_2)\prec y_1$. By the definition of the di-endomorphisms, all of them are equal to $\{y_3,x_2,y_1\}_1$.
            
            \item $ y_3 \prec (\rho_2^* \vdash \rho_1^*) = y_3 \prec (\rho_2^* \dashv \rho_1^*) = (y_3 \prec \rho_2^*) \prec \rho_1^*$. Is the previous case (2d).
        \end{enumerate}

        \item Let us prove that $\alpha_{22} = \beta_{22} = \gamma_{22}$:

        \begin{enumerate}
            \item $R^\triangleright(\lambda_2^* \succ y_1, x_3) = R^\triangleright(\lambda_2^* \prec y_1, x_3) = R^\triangleright(y_1, \lambda_2 \succ x_3)$. For both cases $\lambda_2 = L^\triangleleft(z,u)$ or $\lambda_2 = L^\triangleright(z,u)$ with $z,u \in A$, we have that $\lambda_2 \succ x_3 = \{z,u,x_3\}_3$. Thus, by Proposition \ref{prop:asssecondkind}(8), for any $i \in \{1,2,3\}$, we have $R^\triangleright (y_1, \lambda_2 \succ x_3) = R^\triangleright(\{ u,z,y_1\}_i, x_3) = R^\triangleright(\lambda_2* \prec y_1,x_3) = R^\triangleright(\lambda_2^* \succ y_1,x_3)$.
            
            \item $R^\triangleright(y_2 \succ \rho_1^*,x_3) = R^\triangleright(y_2 \prec \rho_1^*,x_3) = \rho_1 \vdash R^\triangleright(y_2,x_3)$. If $\rho_1 = R^\mu(z,u)$ with $z,u \in A$ and $\mu \in \{\triangleleft, \triangleright\}$, then, by Remark \ref{rem:asssecond} and Proposition \ref{prop:asssecondkind}(8), $\rho_1 \vdash R^\triangleright(y_2,x_3) = R^\triangleright( \{y_2,u,z\}_i,x_3)$ for any $i \in \{1,2,3\}$. Hence, we can conclude that $\rho_1 \vdash R^\triangleright(y_2,x_3) =R^\triangleright(y_2 \succ \rho_1^*,x_3) = R^\triangleright(y_2 \prec \rho_1^*,x_3)$ for all $\rho_1 \in \mathfrak{R}(A)$.

            \item $R^\triangleleft(y_1,x_2) \vdash \rho_3 = R^\triangleright(y_1,x_2) \vdash \rho_3 = R^\triangleright(y_1,x_2 \succ \rho_3)$. By Remark \ref{rem:asssecond}, $R^\triangleleft(y_1,x_2) \vdash \rho_3 = R^\triangleright(y_1,x_2) \vdash \rho_3$. If $\rho_3 = R^\triangleleft(z,u)$ with $z,u \in A$, then, by Proposition \ref{prop:asssecondkind}(7), $R^\triangleleft(y_1,x_2) \vdash \rho_3 =R^\triangleright(y_1,\{x_2,z,u\}_2) = R^\triangleright(y_1, x_2 \succ \rho_3)$. Otherwise, if $\rho_3 = R^\triangleright(z,u)$, then, by Proposition \ref{prop:asssecondkind}(8), $R^\triangleleft(y_1,x_2) \vdash \rho_3= R^\triangleright(y_1,\{x_2,z,u\}_3)=R^\triangleright(y_1,x_2 \succ \rho_3)$. Therefore, for all $\rho_3 \in \mathfrak{R}(A)$ we can assure that $R^\triangleleft(y_1,x_2) \vdash \rho_3 =R^\triangleright(y_1,x_2 \succ \rho_3) $.

            \item $(\rho_1 \dashv \rho_2) \vdash \rho_3 = (\rho_1 \vdash \rho_2) \vdash \rho_3 = \rho_1 \vdash (\rho_2 \vdash \rho_3)$. Because $\mathfrak{R}(A)$ is an associative dialgebra.
        \end{enumerate}
    \end{enumerate}

    Finally, we prove the axiom $(X \vdash Y ) \dashv Z = X \vdash ( Y \dashv Z)$.

    \begin{align*}
        &\left( \begin{pmatrix} \lambda_1 & x_1 \\ \bar{y}_1 & \rho_1\end{pmatrix} \vdash \begin{pmatrix} \lambda_2 & x_2 \\ \bar{y}_2 & \rho_2\end{pmatrix}\right) \dashv \begin{pmatrix} \lambda_3 & x_3 \\ \bar{y}_3 & \rho_3\end{pmatrix} 
        \\
        &= \begin{pmatrix}\lambda_1 \vdash \lambda_2 + L^\triangleright(x_1,y_2) & \lambda_1 \succ x_2 + x_1 \succ \rho_2\\ \overline{\lambda_2^*\prec y_1} + \overline{y_2 \prec \rho_1^*} & R^\triangleright(y_1,x_2) + \rho_1 \vdash \rho_2 \end{pmatrix}\dashv \begin{pmatrix} \lambda_3 & x_3 \\ \bar{y}_3 & \rho_3\end{pmatrix} = \begin{pmatrix} \alpha_{11} & \alpha_{12} \\ \alpha_{21} & \alpha_{22} \end{pmatrix}
    \end{align*}
    where
    \begin{align*}
        \alpha_{11}&= (\lambda_1 \vdash \lambda_2) \dashv \lambda_3 + L^\triangleright(x_1,y_2) \dashv \lambda_3 + L^\triangleleft(\lambda_1 \succ x_2, y_3) + L^\triangleleft(x_1 \succ \rho_2,y_3)
        \\
        \alpha_{12}&= (\lambda_1 \vdash \lambda_2) \prec x_3 + L^\triangleright(x_1,y_2) \prec x_3 + (\lambda_1 \succ x_2) \prec \rho_3 + (x_1 \succ \rho_2) \prec \rho_3
        \\
        \alpha_{21}&=  \overline{\lambda_3^* \succ(\lambda_2^* \prec y_1)} + \overline{\lambda_3^* \succ (y_2 \prec \rho_1^*)} + \overline{y_3 \succ R^\triangleleft(x_2,y_1)} + \overline{y_3 \succ (\rho_2^* \dashv \rho_1^*)}
        \\
        \alpha_{22}&= R^\triangleleft(\lambda_2^*\prec y_1 ,x_3) + R^\triangleleft(y_2 \prec \rho_1^*,x_3) +R^\triangleright(y_1,x_2) \dashv \rho_3 + (\rho_1 \vdash \rho_2) \dashv \rho_3
    \end{align*}

    \begin{align*}
        &\begin{pmatrix} \lambda_1 & x_1 \\ \bar{y}_1 & \rho_1\end{pmatrix} \vdash \left(\begin{pmatrix} \lambda_2 & x_2 \\ \bar{y}_2 & \rho_2\end{pmatrix} \dashv \begin{pmatrix} \lambda_3 & x_3 \\ \bar{y}_3 & \rho_3\end{pmatrix}\right)
        \\
        &=\begin{pmatrix} \lambda_1 & x_1 \\ \bar{y}_1 & \rho_1\end{pmatrix} \vdash \begin{pmatrix} \lambda_2  \dashv \lambda_3 + L^\triangleleft (x_2,y_3)& \lambda_2 \prec x_3 + x_2 \prec \rho_3 \\ \overline{\lambda_3^* \succ y_2} + \overline{y_3 \succ \rho_2^*}& R^\triangleleft(y_2,x_3) + \rho_2 \dashv \rho_3\end{pmatrix} = \begin{pmatrix} \beta_{11} & \beta_{12} \\ \beta_{21} & \beta_{22} \end{pmatrix}
    \end{align*}
    where
    \begin{align*}
        \beta_{11}&= \lambda_1\vdash (\lambda_2 \dashv \lambda_3) + \lambda_1 \vdash L^\triangleleft(x_2,y_3) +L^\triangleright(x_1,\lambda_3^*\succ y_2) + L^\triangleright(x_1,y_3 \succ \rho_2^*)
        \\
        \beta_{12}&= \lambda_1 \succ (\lambda_2 \prec x_3) + \lambda_1 \succ (x_2 \prec \rho_3)+ x_1 \succ R^\triangleleft(y_2,x_3) + x_1 \succ (\rho_2 \dashv \rho_3)
        \\
        \beta_{21}&=  \overline{(\lambda_3^* \vdash \lambda_2^*) \prec y_1} + \overline{L^\triangleright(y_3,x_2)\prec y_1} + \overline{(\lambda_3^* \succ y_2) \prec \rho_1^*} + \overline{(y_3 \succ \rho_2^*) \prec \rho_1^*}
        \\
        \beta_{22}&= R^\triangleright(y_1,\lambda_2 \prec x_3) + R^\triangleright(y_1, x_2 \prec \rho_3) + \rho_1 \vdash R^\triangleleft(y_2,x_3) + \rho_1 \vdash (\rho_2 \dashv \rho_3)
    \end{align*}
    Since the arguments are similars as in the previous cases, we just give an indication.
   \begin{enumerate}
       \item Let us prove that $\alpha_{11} = \beta_{11}$:
       \begin{enumerate}
           \item $(\lambda_1 \vdash \lambda_2) \dashv \lambda_3 = \lambda_1 \vdash ( \lambda_2 \dashv \lambda_3)$. Because $\mathfrak{L}(A)$ is an associative dialgebra.
           
           \item $L^\triangleright(x_1,y_2) \dashv \lambda_3 = L^\triangleright(x_1, \lambda_3^* \succ y_2)$. By Proposition \ref{prop:asssecondkind}(2).
           
           \item $L^\triangleleft(\lambda_1 \succ x_2, y_3) = \lambda_1 \vdash L^\triangleleft(x_2,y_3)$. By Proposition \ref{prop:asssecondkind}(3).
           
           \item  $L^\triangleleft(x_1 \succ \rho_2,y_3) = L^\triangleright(x_1, y_3 \succ \rho_2^*)$. By Proposition \ref{prop:asssecondkind}(2) and (3).
       \end{enumerate}

       \item Let us prove that $\alpha_{12} = \beta_{12}$:
       \begin{enumerate}
           \item $(\lambda_1 \vdash \lambda_2) \prec x_3 = \lambda_1 \succ ( \lambda_2 \prec x_3)$. By Lemma \ref{lemma:diendomorphism}(2).
           
           \item $L^\triangleright(x_1,y_2) \prec x_3 = x_1 \succ R^\triangleleft(y_2,x_3)$. By definition of the di-endomorphims, both are equal to $\{x_1,y_2,x_3\}_2$.
           
           \item $(\lambda_1 \succ x_2) \prec \rho_3 = \lambda_1 \succ (x_2 \prec \rho_3)$.
           
           \item $(x_1 \succ \rho_2) \prec \rho_3 = x_1 \succ (\rho_2 \dashv \rho_3)$. By Lemma \ref{lemma:diendomorphism}(2).
       \end{enumerate}

       \item Let us prove that $\alpha_{21} = \beta_{21}$:
       \begin{enumerate}
           \item $\lambda_3^* \succ (\lambda_2^* \prec y_1) =(\lambda_3^* \vdash \lambda_2^*) \prec y_1$. By Lemma \ref{lemma:diendomorphism}(2).
           
           \item $\lambda_3^* \succ (y_2 \prec \rho_1^*) = (\lambda_3^* \succ y_2) \prec \rho_1^*$. Is the previous case (2c).
           
           \item $y_3 \succ R^\triangleleft(x_2,y_1) = L^\triangleleft (y_3,x_2 \prec y_1)$. By definition of the di-endomorphims, both are equal to $\{y_3,x_2,y_1\}_2$.
           
           \item $y_3 \succ (\rho_2^* \dashv \rho_1^*) = (y_3 \succ \rho_2^*)\prec \rho_1^*$. By Lemma \ref{lemma:diendomorphism}(2).
       \end{enumerate}

       \item Let us prove that $\alpha_{22} = \beta_{22}$:
       \begin{enumerate}
           \item $R^\triangleleft(\lambda_2^* \prec y_1,x_3) = R^\triangleright(y_1, \lambda_2 \prec x_3)$. By Propoposition \ref{prop:asssecondkind}(6) and (7).
           
           \item $R^\triangleleft(y_2 \prec \rho_1^*, x_3) = \rho_1 \vdash R^\triangleleft(y_2,x_3)$. By Proposition \ref{prop:asssecondkind}(7).
           
           \item $ R^\triangleright(y_1,x_2) \dashv \rho_3 = R^\triangleright(y_1,x_2 \prec \rho_3)$. By Proposition \ref{prop:asssecondkind}(6).
           
           \item $(\rho_1 \vdash \rho_2) \dashv \rho_3 = \rho_1 \vdash (\rho_2 \dashv \rho_3)$. Because $\mathfrak{R}(A)$ is an associative dialgebra.
       \end{enumerate}
   \end{enumerate} 
\end{proof}

\section*{Acknowledgment}

The first author would like to thank the Department of Analysis and Applied Mathematics at Complutense University for their hospitality during my sabbatical stay. The first author also thanks the algebra group at the Rey Juan Carlos University of Madrid for their support. Additionally, the first author acknowledges the financial support from CONAHCYT through its call for sabbatical year applications, first phase, 2024.

The second author was supported by Proyecto puente URJC “Lie algebras, Jordan systems and related structures” [grant number M3329]. 

\section*{Conflict of interest}
Authors state no conflict of interest.

\end{document}